\documentclass[11pt]{amsart}
\usepackage{mathrsfs}
\usepackage{amssymb}
\usepackage{amsmath,amscd}
\usepackage{amsfonts,amssymb}
\usepackage{latexsym}
\usepackage{amsbsy}

\usepackage{graphics,graphicx}
\usepackage[all]{xy}

\newtheorem{theorem}{Theorem}[subsection]
\newtheorem*{theorem*}{Theorem}

\newtheorem{lemma}{Lemma}[subsection]
\newtheorem{proposition}{Proposition}[subsection]

\newtheorem*{conjecture*}{Conjecture}

\newtheorem*{remark*}{Remark}
\newtheorem{definition}{Definition}[subsection]

\newcommand{\bcen}{\begin{center}}      \newcommand{\ecen}{\end{center}}
\newcommand{\bay}{\begin{array}}      \newcommand{\eay}{\end{array}}
\newcommand{\beq}{\begin{eqnarray*}}      \newcommand{\eeq}{\end{eqnarray*}}
\def\lr#1{\langle #1\rangle}
\newcommand{\cal}{\mathcal}
\def\az{\alpha}    
\def\bz{\beta}     
\def\gz{\gamma}
   
\def\sz{\sigma}

\def\dz{\delta}  

\def\llz{\Lambda}      
\def\lz{\lambda}    
     \def\cp{{\cal P}} \def\cc{{\cal C}}
\def\ch{{\cal H}}   \def\ct{{\cal T}} \def\co{{\cal O}}

\def\ca{{\cal A}}
\def\cb{{\cal B}}
\def\cz{{\cal Z}}

\def\ct{{\cal T}}

\def\cn{{\cal N}}
\def\bbn{{\mathbb N}}  \def\bbz{{\mathbb Z}}  \def\bbq{{\mathbb Q}}
\def\bbf{{\mathbb F}}   \def\bbe{{\mathbb E}} \def\bbf{{\mathbb F}}
\def\fq{{\mathbb F}_q}
  \def\bbc{{\mathbb C}}
\def\bba{{\mathbb A}}

\def\te{{\tilde E}}
\def\bi{{\bf 1}}
   
\def\fk#1{{\mathfrak #1}}

\def\lra{\longrightarrow}  

\def\ra{\rightarrow}
\def\hom{\operatorname{Hom}}
 \def\ind{\operatorname{ind}}

\def\Im{\operatorname{Im}}
\def\ext{\operatorname{Ext}}
\def\aut{\operatorname{Aut}}

\def\dim{\operatorname{dim}}
\def\cdim{\operatorname{codim}}
\def\min{\operatorname{min}}
\def\udim{\operatorname{\underline {dim}}}
\def\ed{\operatorname{End}}
\def\Hom{\operatorname{Hom}}
\def\Ext{\operatorname{Ext}}
\def\ol{\overline}
\def\mod{\operatorname{mod}}  
\def\Im{\operatorname{Im}}

\def\uq2{U_q(\hat{sl}_2)}
\def\fq{\bbf_q}

\def\wm{{\bf m}}

\def\nd{{\noindent}}
\def\mk{{\medskip}}
\def\hs{{\hskip 1cm}}

\def\det{\operatorname{det}}

\title[Canonical bases]{Canonical bases for the quantum extended Kac-Moody algebras and Hall polynomials}

\thanks{The research was
supported in part by  NSCF grant 10931006 and and by STCSM(12ZR1413200)  \\
${\quad}$The research was also supported in part by Research
Institute for Mathematical Sciences, Kyoto University, Kyoto, California Institute of Technology, Pasadena, California and Chern Institute of Mathematics, Nankai University, Tianjin,China}
\author{Guanglian Zhang}
\address{Department of Mathematics\\
Shanghai Jiao Tong University\\
Shanghai 200240, China}
 \email{g.l.zhang@sjtu.edu.cn}

\date{\today}

\begin{document}

\begin{abstract}
In this paper, the singular Ringel-Hall algebra for a tame quiver is
introduced and shown to be isomorphic to the positive part of the
quantum extended Kac-Moody algebra. A PBW basis is constructed and a
new class of perverse sheaves is shown to have purity property. This
allows to construct the canonical bases of the positive part of the
quantum extended Kac-Moody algebra. As an application, the existence
of Hall polynomials for tame quiver algebras is proved.

\end{abstract}

\maketitle

\bigskip
\bigskip

\centerline{\bf 0. Introduction}

\bigskip
The early nineties of the last century witnessed the invention of
canonical bases independently by Lusztig~\cite{L1,L2,L3} and by
Kashiwara~\cite{K}. Since then, canonical bases have been playing a
vital role in representation theory of quantum groups, Hecke
algebras, and quantized Schur algebras ~\cite{A,BN,LXZ,F1,F2,
Ro,VV}.

Around the same time Ringel--Hall algebras (of abelian categories)
were introduced by Ringel~\cite{R2} in order to obtain a categorical
version of Gabriel's theorem~\cite{BGP}. Ringel~\cite{R1} showed
that the Ringel--Hall algebra of the category of finite-dimensional
representations over a finite field of a Dynkin quiver, after a
twist by the Euler form, is isomorphic to the positive part of the
quantum group associated to the underlying graph of the quiver. This
was generalised by Green~\cite{G} to all quivers, with the whole
Ringel--Hall algebra replaced by its subalgebra generated by the
1-dimensional representations (the composition algebra).

The above relationship between Ringel--Hall algebras and quantum
groups enabled Lusztig to give geometric construction of canonical
bases by using quiver varieties~\cite{L1,L2,L3}. For finite types he
also has an algebraic construction ~\cite{L4}, and this was
generalised to affine types by Lin, Xiao and the author
in~\cite{LXZ}. The approach of Kashiwara is combinatorial~\cite{K}.

In this paper we study canonical bases for quantum extended affine
Kac-Moody algebras. We follow Lusztig's geometric approach. For this
purpose, we introduce the \emph{singular Ringel--Hall algebra}.
Thanks to a theorem of Sevenhant and Van den Bergh~\cite{SV}, the
whole Ringel--Hall algebra $\ch(Q)$ of a tame quiver $Q$ can be
obtained from the composition algebra $\cc(Q)$ by adding some
imaginary simple roots. The singular Ringel--Hall algebra $\ch^s(Q)$
will be defined as a subalgebra of $\ch(Q)$, obtained from $\cc(Q)$
by adding certain imaginary simple roots. It will be shown that
$\ch^s(Q)$ is isomorphic to the positive part of the quantum
extended Kac-Moody algebra of the underlying graph of $Q$. To
construct the canonical bases, we start with constructing a PBW
basis and, as in ~\cite{L1}, proving that certain perverse sheaves
over the associated quiver varieties have the purity property. It
then follows that certain elements in the singular Ringel--Hall
algebra associated to these sheaves form the canonical bases. In
\cite {KS}, Kang and Schiffmann obtained canonical bases of quantum
generalized Kac-Moody algebras via construct a class quiver with
edge loops. An important difference between our construction and
Kang and Schiffmann's construction is that our canonical bases are
made of simple perverse sheaves rather than semisimple perverse
sheaves, that is, our canonical bases coincide with classical
canonical bases of Lusztig's version.

This work is another motivated by the application of canonical bases
of the Fock space over the quantum affine $\mathfrak{sl}_n$ in the
study of decomposition numbers: the Lascoux--Leclerc--Thibon
conjecture ~\cite{LLT} which was solved by Ariki~\cite{A}, and the
Lusztig conjecture for $q$-Schur algebras by Varagnolo and
Vasserot~\cite{VV}. These Fock spaces are of type $A$ and lever $1$.
Those of type $A$ and of higher levers are studied in~\cite{A,KS}
and other papers. Our aim is to solve Lusztig's conjecture for
$q$-Schur algebras of other types. The first step is to find the
appropriate Fock spaces. We propose that the quantum extended affine
Kac--Moody algebra modulo certain elements corresponding to unstable
orbits in Nakajima's sense , considered as a module over the quantum
affine algebra, is the Fock space. For type $A$, this is verified.
 Our forthcoming work are canonical bases for generalized Kac-Moody Lie algebra,
 and understand quantum geometric
 Langlands duality by using our results about canonical bases.

  The paper is organized as follows. In \S 1 we
 give a quick review of the definitions of Ringel-Hall algebras and Double Ringel-Hall algebras.
In \S 2 we define the singular Ringel-Hall algebras $\ch^s(\llz)$
and study the relations between the singular Ringel-Hall algebras
and the quantized universal enveloping algebras. We prove
Proposition~2.2.1 and from this, we point out that the quotient of
the singular Ringel-Hall algebras modulo unstable orbits is
isomorphic to $q-$ Fock space in the case of type $A$. In \S 3 we
show that  the singular Ringel-Hall algebra $\ch^s(\llz)$ is
isomorphic to the positive part  $\mathbb{U}^+$ of the quantum
extended Kac-Moody algebra. In \S 4 we construct the $PBW$ type
basis for the positive part  $\mathbb{U}^+$ of the quantum extended
Kac-Moody algebras. In \S 5 we give the main theorem of this paper, that is, a description of the canonical basis of $\ch^s(\llz)\cong\mathbb{U}^+$. In $\S 6$ we prove that the closure of
semi-simple objects in $\ct_i$ have purity property.
In \S 7 we
study the fibres of $p_3$. We also give a new class of perverse
sheaves with purity property. In $\S 8$ we give the proof of Theorem 5.1.1.
In \cite{H}, A.Hubery proved the existence
of Hall polynomials on the tame quivers for Segre classes. In \S
9, by using the extension algebras of singular Ringel-Hall
algebras, we give a simple and direct proof for the existence of
Hall polynomials for the tame quivers.

\mk\nd{\bf Acknowledgments.} I would like to express my sincere
gratitude to Hiraku Nakajima, Xinwen Zhu and Jie Xiao for a number of interesting discussions.


\section{The Ringel--Hall algebra}\label{s:hall-alg}

Throughout the paper, let $\bbf_q$ denote a finite field with $q$
elements, and $k=\overline{\bbf}_q$ be the algebraic closure of $\bbf_q$.

\subsection{}\label{ss:1.1} A \emph{quiver} $Q=(I,H,s,t)$ consists of  a
vertex set $I$,  an arrow set $H$, and two maps $s,t:H\rightarrow I$
such that an arrow $\rho\in H$ starts at $s(\rho)$ and terminates at
$t(\rho).$ The two maps $s$ and $t$ extends naturally to the set of paths.
A \emph{representation} $V$ of $Q$ over a field $F$ is a collection $\{V_i:i\in I\}$ of $F$-vector spaces
and a collection $\{V(\rho):V_{s(\rho)}\rightarrow V_{t(\rho)}:\rho\in H\}$ of $F$-linear maps.

Let $\llz=\bbf_q Q$ be the \emph{path algebra} of $Q$ over the field $\bbf_q$.
Precisely, $\llz$ is a $\bbf_q$-vector space with basis all paths of $Q$ (including the trivial paths attached to all vertices), and the product $pq$ of two paths $p$ and $q$ is the concantenation of $p$ and $q$ if $s(p)=t(q)$ and is $0$ otherwise.
The
category of all finite-dimensional left $\llz$-modules, or
equivalently finite left $\llz$-modules, is equivalent to the category of finite-dimensional representations of
$Q$ over $\bbf_q.$ We shall simply identify $\llz$-modules with
representations of $Q$, and call a $\llz$-module \emph{nilpotent} if the corresponding representation is nilpotent. The category of nilpotent $\llz$-modules will be denoted by $\mod\llz$.

The set of isomorphism classes of nilpotent simple $\llz$-modules
is naturally indexed by the set $I$ of vertices of $Q.$ Hence the
Grothendieck group $G(\llz)$ of $\mod\llz$ is the free Abelian group
$\bbz I$. For each nilpotent $\llz$-module $M$, the dimension vector
$\udim M=\sum_{i\in I}(\dim M_i)i$ is an element of $ G(\llz)$.

The Euler form $\lr{-,-}$ on $G(\llz)=\bbz I$ is defined by
\beq\hs \lr{\az,\bz}=\sum_{i\in I}a_ib_i-\sum_{\rho\in
H}a_{s(\rho)}b_{t( \rho)}\eeq for $\az=\sum_{i\in I}a_i i$ and
$\bz=\sum_{i\in I}b_i i$ in $\bbz I.$ For any nilpotent $\llz$-modules $M$ and $N$ one has
$$\lr{\udim M, \udim
N}=\dim_{\fq}\hom_{\llz}(M,N)-\dim_{\fq}\ext_{\llz}(M,N).$$ The
symmetric Euler form is defined as
$$(\az,\bz)=\lr{\az,\bz}+\lr{\bz,\az}\ \ \text{for}\ \
\az,\bz\in\bbz I.$$ This gives rise to a symmetric generalized
Cartan matrix $C=(a_{ij})_{i,j\in I}$ with $a_{ij}=(i,j).$ It is
easy to see that $C$ is independent of the field $\fq$ and the
orientation of $Q.$


\subsection{The Ringel--Hall algebra}\label{ss:hall-alg}
Let $Q$ be a quiver, and $\bbf_qQ$ be the path algebra of $Q$ over the finite field $\bbf_q$.

Let $v=v_q=\sqrt q\in \bbc$ and $\cp$
 be the set of isomorphism classes of finite-dimensional nilpotent $\llz$-modules.
The \emph{(twisted) Ringel--Hall algebra} $\ch^*(\llz)$ is defined as the $\bbq(v)$-vector space with basis $\{u_{[M]}:[M]\in\cp\}$ and with multiplication  given by
$$u_{[M]}\ast
u_{[N]}=v^{\lr{\udim M, \udim
N}}\sum_{[L]\in\cp}g^{L}_{MN}u_{[L]}.$$
where $g_{MN}^{L}$ is the number of $\llz$-submodules
$W$ of $L$ such that $W\simeq N$ and $L/W\simeq M$ in $\mod \llz$. Its subalgebra generated by $\{u_{[M]}:M \text{ is a simple nilpotent } \Lambda \text{ module}\}$ is called the \emph{composition subalgebra} of $\ch^*(\llz)$ or the \emph{composition Ringel--Hall algebra} of $\llz$, and denoted by $\cc^*(\llz)$.
The Ringel--Hall algebra $\ch^*(\llz)$ and the composition algebra $\cc^*(\llz)$ is graded by $\bbn I$, namely, by dimension vectors of modules, since $g_{MN}^L\neq 0$ if and only if there is a short exact sequence $0\rightarrow N\rightarrow L\rightarrow M\rightarrow 0$, which implies that $\udim L =\udim M+\udim N$.
Following [R3], for any nilpotent
$\llz$-module $M$, we denote $$\lr{M}=v^{-\dim
 M+\dim\ed_{\llz}(M)}u_{[M]}.$$  Note that $\{ \lr{M} \;|\; M \in
\cp\}$ is  a $\bbq(v)$-basis of $\ch^*(\llz)$.

The $\bbq(v)$-algebra $ \ch^*(\llz)$ depends on $q(=v^2)$. We will use
$ \ch_q^*(\llz)$ to indicate the dependence on $q$ when such a need
arises.

\subsection{A construction by Lusztig}\label{ss:construction-lusztig}
Let  $V=\bigoplus_{i\in I}V_i$ be a finite-dimensional
$I$-graded $k$-vector space with a given
$\bbf_q$-rational structure by the Frobenius map $F:k\rightarrow k, a\mapsto a^q$. Let $\bbe_V$ be
the subset of $\bigoplus_{\rho\in H}\hom(V_{s(\rho)},V_{t(\rho)})$
consisting of elements which define nilpotent representations of $Q$. Note that
$\bbe_V=\bigoplus_{\rho\in H}\hom(V_{s(\rho)},V_{t(\rho)})$ when $Q$
has no oriented cycles. The space of $\bbf_q$-rational points of
$\bbe_V$ is the fixed-point set $\bbe^F_V.$

Let $G_V=\prod_{i\in I}GL(V_i)$, and its subgroup of
$\bbf_q-$rational points be $G^F_V$. Then the group $G_V=\prod_{i\in
I}GL(V_i)$ acts naturally on  $\bbe_V$ by
$$(g,x)\mapsto g\bullet x=x'\ \ \text{where}\ \
x'_{\rho}=g_{t(\rho)}x_{\rho}g^{-1}_{s(\rho)}\ \ \text{for all}\ \
\rho\in H.$$ For $x\in\bbe_V$, we denote by $\co_x$ the orbit of $x$.
This action restricts to an action of the finite group $G^F_V$ on $\bbe^F_V.$

For $\gz\in\bbn I,$ we fix a $I$-graded $k$-vector space $V_{\gz}$
with $\udim V_{\gz}=\gz.$ We set $\bbe_{\gz}=\bbe_{V_{\gz}}$ and
$G_{\gz}=G_{V_{\gz}}.$ For $\az,\bz\in\bbn I$ and $\gz=\az +\bz,$ we
consider the diagram
$$\xymatrix{\bbe_{\az}\times\bbe_{\bz} & \bbe' \ar[l]_(0.35){p_1}\ar[r]^{p_2}&\bbe''\ar[r]^{p_3}& \bbe_{\gz}.}$$ Here $\bbe''$ is the set of all pairs $(x,W)$,
consisting of $x\in\bbe_{\gz}$ and an $x$-stable $I$-graded subspace
$W$ of $V_{\gz}$ with $\udim W=\bz$, and $\bbe'$ is the set of all
quadruples $(x,W,R',R'')$, consisting of $(x,W)\in\bbe''$ and two
invertible linear maps $R':k^{\bz}\ra W$ and $R'':k^{\az}\ra
k^{\gz}/W.$  The maps are defined in an obvious way:
$p_2(x,W,R',R'')=(x,W),$ $p_3(x,W)=x,$ and
$p_1(x,W,R',R'')=(x',x''),$ where
$x_{\rho}R'_{s(\rho)}=R'_{t(\rho)}x'_{\rho}$ and
$x_{\rho}R''_{s(\rho)}=R''_{t(\rho)}x''_{\rho}$ for all $\rho\in H.$
By Lang's Lemma, the varieties and morphisms in this diagram are
naturally defined over $\bbf_q.$ So we have
$$\xymatrix{\bbe^F_{\az}\times\bbe^F_{\bz} & \bbe'^F \ar[l]_(0.35){p_1}\ar[r]^{p_2}&\bbe''^F\ar[r]^{p_3}& \bbe^F_{\gz}.}$$
For $M\in\bbe_{\az}, N\in\bbe_{\bz}$ and $L\in\bbe_{\az+\bz},$ we
define
$${\bf Z}=p_2p_1^{-1}(\co_M\times \co_N)\subseteq E'',\qquad {\bf Z}_{L,M,N}={\bf Z}\cap p_3^{-1}( \co_L).$$

For any map $p: X\rightarrow Y$ of finite sets,
$p^*:\bbc(Y)\rightarrow \bbc(X)$ is defined by $p^*(f)(x)=f(p(x))$
and $p_!: \bbc(X)\rightarrow \bbc(Y)$ is defined by $
p_!(h)(y)=\sum_{x \in p^{-1}(y)}h(x),$ on the integration along the
fibers. Let $\bbc_{G^F}(\bbe^F_V)$ be the space of $G^F_V$-invariant
functions $\bbe^F_V\ra \bbc ( \text{ or } \overline{Q_l}).$ Given
$f\in\bbc_{G^F}(\bbe^F_\az)$ and $g\in\bbc_{G^F}(\bbe^F_\bz)$, there
is a unique $h\in \bbc_G(\bbe''^F)$ such that
$p_2^*(h)=p_1^*(f\times g).$ Then define $f\circ g$ by
$$f\circ g=(p_3)_!(h)\in\bbc_{G^F}(\bbe^F_{\gz}).$$

 Let $${\wm}(\az,\bz)= \sum_{i\in
I}a_ib_i+\sum_{\rho\in H}a_{s(\rho)}b_{t( \rho)}.$$ We again define
the multiplication in the $\bbc$-space ${\bf K}=\bigoplus_{\az\in\bbn
I}\bbc_{G^F}(\bbe^F_{\az})$ by
$$f\ast g=v_q^{-{\wm}(\az,\bz)}f\circ g$$ for all $f\in\bbc_{G^F}(\bbe^F_\az)$ and
$g\in\bbc_{G^F}(\bbe^F_\bz).$ Then $({\bf K},\ast)$ becomes an
associative $\bbc$-algebra.

For $M\in\bbe^F_{\az}$, let  $\co_M\subset\bbe_{\az}$ be the
$G_{\az}$-orbit of $M.$ We take
$\bi_{[M]}\in\bbc_{G^F}(\bbe^F_{\az})$ to be the characteristic
function of $\co^F_M,$ and set $f_{[M]}=v_q^{-\dim\co_M}\bi_{[M]}.$
We consider the subalgebra $({\bf L},\ast)$ of $({\bf K},\ast)$
generated by $f_{[M]}$ over $\bbq(v_q),$ for all $M\in\bbe^F_{\az}$
and all $\az\in\bbn I.$ In fact ${\bf L}$ has a $\bbq(v_q)$-basis
$\{f_{[M]}|M\in\bbe^F_{\az}, \az\in\bbn I\}.$ Since  $\bi_{[M]}\circ
\bi_{[N]}(W)=g^W_{M N}$ for any $W\in\bbe^F_{\gz},$ we have


\begin{proposition}\cite{LXZ}\label{p:1.3.1}
The linear map
$\varphi:({\bf L},\ast)\lra \ch^*(\llz)$ defined by
$$\varphi(f_{[M]})=\lr{M},\ \ \ \text{for all}\ [M]\in\cp$$
is an isomorphism of associative $\bbq(v_q)$-algebras.
\end{proposition}

\subsection{The double Ringel--Hall algebra
$\mathscr{D}(\llz)$}\label{ss:double-hall-alg} First, we define a
Hopf algebra $\bar\ch^+(\llz)$ which is a $\bbq(v)$-vector space with
the basis $\{K_{\mu}u_{\alpha}^+|\mu\in\bbz[I],\az\in\cp\},$ whose
Hopf algebra structure is given as

(a) Multiplication (\cite{R1})
\begin{eqnarray*}u_{\az}^+*u_{\bz}^+&=&v^{\lr{\az,\bz}}\sum_{\lz\in\cp}g^{\lz}_{\az\bz}u^+_{\lz},
\text{ for all } \az,\bz\in\cp,\\
K_{\mu}*u_{\az}^+&=&v^{(\mu,\az,)}u_{\az}^+*K_{\mu},\text{ for all
} \az\in\cp,
\mu\in\bbn[I],\\
K_\mu*K_\nu&=&K_\nu*K_\mu\hspace{7pt}=\hspace{7pt}K_{\mu+\nu}, \text{ for all }
\mu,\nu\in\bbn[I].
\end{eqnarray*}

(b)Comultiplication (\cite{G})
\begin{eqnarray*}
\bigtriangleup(u_{\lz}^+)&=&\sum_{\az,\bz\in\cp}v^{\lr{\az,\bz}}
\frac{a_{\az}a_{\bz}}{a_{\lz}}g^{\lz}_{\az\bz}u^+_{\az}K_{\bz}\otimes
u_{\bz}^+, \text{ for all } \lz\in\cp,\\
\bigtriangleup(K_{\mu})&=&K_\mu\otimes K_\mu,\text{ for all } \mu\in\bbn[I]
\end{eqnarray*}
with counit $\epsilon(u_{\lz}^+)=0,$ for all $0\neq \lz\in\cp$, and
$\epsilon(K_\mu)=1.$ Here $a_\lz$ denotes the cardinality of the
finite set $\aut_{\llz}(M)$ with $[M]=\lz.$

(c)Antipode(\cite{X})
\begin{eqnarray*}
   S(u^+_{\lz})&=&\dz_{\lz0}+\sum(-1)^m\sum_{\pi\in\cp,\lz_1,\cdots,\lz_m\in\cp\setminus \{0\}}\times
\\&&\mbox{}v^{2\sum_{i<j}\lr{\lz_i,\lz_j}}\frac{a_{\lz_1}\cdots a_{\lz_m}}{a_{\lz}}g^\lz_{\lz_1\cdots\lz_m}
g^\pi_{\lz_1\cdots\lz_m}K_{-\lz}u^+_\pi, \text{ for all } \lz\in\cp,\\
S(K_\mu)&=&K_{-\mu} \text{  for all } \mu\in\bbz[I].
\end{eqnarray*}
 The subalgebra $\ch^+$ of $\bar\ch^+(\llz)$  generated by $\{u_\lz|\lz\in\cp\}$ is isomorphic to $\ch^*(\llz).$ Moreover, we have an isomorphism of vector spaces $\bar\ch^+(\llz)\cong \ct\otimes \ch^+$ , where $\mathcal {T}$ denotes the torus subalgebra generated by
$\{K_\mu:\mu\in\bbz[I]\}.$

  Dually, we can define a Hopf algebra $\bar\ch^-(\llz)$. Following
  Ringel, we have a bilinear form
  $\varphi:\bar\ch^+(\llz)\times\bar\ch^-(\llz){\lra}\bbq(v)$ defined by
  $$\varphi(K_\mu u_\az^+,K_\nu u^-_{\bz})=v^{-(\mu,\nu)-(\az,\nu)+(\mu,\bz)}\frac{|V_\az|}{a_\az}\dz_{\az\bz}$$
for all $\mu,\nu\in\bbz[I]$ and all $\az,\bz\in \cp.$ Thanks to
\cite{X}, we can form the reduced Drinfeld double
$\mathscr{D}(\llz)$ of the Ringel--Hall algebra of $\llz$, which admits a
triangular decomposition
$$\mathscr{D}(\llz)=\ch^-\otimes\mathcal {T}\otimes\ch^+.$$
It is a Hopf algebra, and the restriction of this structure on
$\bar\ch^-(\llz)=\ch^-\otimes\ct$ and $\bar\ch^+(\llz)=\ct\otimes\ch^+$ are given as above.

The subalgebra of $\mathscr{D}(\llz)$ generated by
$\{u_i^{\pm},K_{\mu}|i\in I, \mu\in\bbz[I]\}$ is also called the
\emph{composition algebra} of $\llz$ and denoted by $\mathcal {C}(\llz).$ It is
 a Hopf subalgebra of $\mathscr{D}(\llz)$ and admits a triangular decomposition
$$\mathcal {C}(\llz)=\mathcal {C}^{-}(\llz)\otimes\mathcal {T}\otimes\mathcal {C}^+(\llz),$$
where $\mathcal {C}^+(\llz)$ is the subalgebra generated by
$\{u_i^+:i\in I\}$ and $\mathcal {C}^{-}(\llz)$ is defined dually.
Moreover, the restriction $\varphi:\mathcal
{C}^+(\llz)\times\mathcal {C}^{-}(\llz){\lra}\bbq(v)$ is
non-degenerate (see\cite{HX}).

 In addition, $\mathscr{D}(\llz)$ admits an involution $\omega$
 defined by
\begin{eqnarray*}
\omega(u_{\lz}^+)=u_{\lz}^{-},~~ \omega(u_{\lz}^{-})=u_{\lz}^{+},
\text{ ~for all ~} \lz\in\cp;
\\
\omega(K_{\mu})=K_{-\mu},\text{ for all } \mu\in\bbn[I].\qquad\qquad
\end{eqnarray*}
We have $\varphi(x,y)=(\omega(x),\omega(y)).$ Obviously, $\omega$
induces an involution of $\mathcal {C}(\llz).$

\subsection{Tame quivers}\label{ss:tame-quiver}
Let $Q=(I,H,s,t)$ be a connected tame quiver, that is, a quiver whose underlying graph is an extended Dynkin diagram
of type $\tilde{A}_n$, $\tilde{D}_n$, $\tilde{E}_6$, $\tilde{E}_7$ or $\tilde{E}_8$, and which has no oriented cycles,
i.e. in $Q$ there are no paths $p$ with $s(p)=t(p)$.
We say that $Q$ is of type $\tilde{A}_{p,q}$ if the underlying graph of $Q$ is of type $\tilde{A}_{p+q-1}$ and there are $p$ clockwise oriented arrows and $q$ anti-clockwise oriented arrows.

Let $\Lambda=\fq Q$ be the path algebra of $Q$ over $\fq$.
Any (nilpotent) finite-dimensional $\Lambda$-module is a direct sum of modules of three types: preprojective, regular, and preinjective. The set of isomorphism classes of preprojective modules (respectively, preinjective modules) will be denoted by $\cp_{prep}$ (respectively, $\cp_{prei}$).  The regular part consists of a family of homogeneous tubes (i.e. tubes of period 1) and a finite number of non-homogeneous tubes, say $\ct_1,\ldots,\ct_l$ respecitvely of periods $r_1,\ldots,r_l$.
We have the following well-known results.


\begin{lemma}\label{l:2.1.1}
\begin{itemize}
\item[(a)] We have $l\leq 3$ and $\sum_{i=1}^l(r_i-1)=|I|-2$.
\item[(b)] Let $P$ be preprojective, $R$ be regular and $I$ be preinjective. Then
\[\Hom(R,P)=\Hom(I,R)=\Hom(I,P)=0,\]
\[\Ext(P,R)=\Ext(R,I)=\Ext(P,I)=0.\]
\item[(c)] Let $R$ and $R'$ be indecomposable regular in different tubes. Then
\[\Hom(R,R')=0=\Ext(R,R').\]
\item[(d)] Let $M\in\ct_i$ for some $1\leqslant i\leqslant l$, and
$$0{\lra}M_2{\lra}M{\lra}M_1{\lra}0$$
be a short exact sequence. Then $M_1\cong I_1\oplus N_1,M_2\cong
P_2\oplus N_2$, where $P_2$ is preprojective, $N_1,N_2\in\ct_i$, and
$I_1$ is preinjective.
\end{itemize}
\end{lemma}

Each tube is an abelian subcategory of $\mod\llz$ and a simple object in a tube is called a \emph{regular simple module} in $\mod\llz$. For a tube of period $r$, there are precisely $r$ simple objects (up to isomorphism), and the sum of dimension vectors of these regular simples is independent of the tube, and this sum will be denoted by $\delta$. Recall that  in Section~\ref{ss:1.1} we have defined the symmetric Euler form $(-,-)$ on $G(\llz)=\bbz I$. In our case, this form is positive semi-definite and its radical is free of rank $1$ generated by $\delta$.  We collect information on the invariants $l$, $r_1,\ldots,r_l$ and $\delta$ of connected tame quivers in the following table.

\[\begin{tabular}{|c|c|c|c|}\hline
type & ~~$\ell$~~ & ~~$r_1,\ldots,r_{\ell}$~~ & $\delta$\\ \hline
\raisebox{-15pt}[0pt]{$\tilde{A}_{n-1,1}$} & \raisebox{-15pt}[0pt]{1} & \raisebox{-15pt}[0pt]{$n-1$} & \xymatrix@!=0.2pc{&&1&&\\ 1 \ar@{-}[urr] \ar@{-}[r]&1\ar@{.}[rr]&&1\ar@{-}[r] &1\ar@{-}[ull] }\\ \hline
\raisebox{-15pt}[0pt]{$\tilde{A}_{p,q}(p,q\geq 2)$} & \raisebox{-15pt}[0pt]{2} & \raisebox{-15pt}[0pt]{$p$, $q$} & \xymatrix@!=0.2pc{&&1&&\\ 1 \ar@{-}[urr] \ar@{-}[r]&1\ar@{.}[rr]&&1\ar@{-}[r] &1\ar@{-}[ull] }\\ \hline
\raisebox{-30pt}[0pt]{$\tilde{D}_n (n\geq 3)$} & \raisebox{-30pt}[0pt]{3} & \raisebox{-30pt}[0pt]{$2,2,n-2$} & \xymatrix@!=0.2pc{1 &&&&&1\\ &2\ar@{-}[ul]\ar@{-}[dl]\ar@{-}[r]&2\ar@{.}[r]&2\ar@{-}[r]&2\ar@{-}[ur]\ar@{-}[dr] &\\ 1 &&&&&1} \\ \hline
\raisebox{-30pt}[0pt]{$\tilde{E}_6$} & \raisebox{-30pt}[0pt]{3} & \raisebox{-30pt}[0pt]{2,3,3} & \xymatrix@!=0.2pc{&&1&&\\ &&2\ar@{-}[u]\ar@{-}[d]&&\\ 1\ar@{-}[r]&2\ar@{-}[r]&3\ar@{-}[r]&2\ar@{-}[r]&1}\\ \hline
\raisebox{-15pt}[0pt]{$\tilde{E}_7$} & \raisebox{-15pt}[0pt]{3} & \raisebox{-15pt}[0pt]{2,3,4} & \xymatrix@!=0.2pc{&&&2\ar@{-}[d]&&&\\ 1\ar@{-}[r]&2\ar@{-}[r]&3\ar@{-}[r]&4\ar@{-}[r]&3\ar@{-}[r]&2\ar@{-}[r]&1}\\ \hline
\raisebox{-15pt}[0pt]{$\tilde{E}_8$} & \raisebox{-15pt}[0pt]{3} & \raisebox{-15pt}[0pt]{2,3,5} & \xymatrix@!=0.2pc{&&3\ar@{-}[d]&&&&&\\ 2\ar@{-}[r]&4\ar@{-}[r]&6\ar@{-}[r]&5\ar@{-}[r]& 4\ar@{-}[r]& 3\ar@{-}[r]& 2\ar@{-}[r]&1}\\ \hline
\end{tabular}\]
\bigskip

Let $K$ be the path algebra over $\bbf_q$ of the Kronecker quiver
$\xymatrix{\cdot\ar@<.5ex>[r]\ar@<-.5ex>[r] & \cdot}$. Then there is
an embedding $\mod K\hookrightarrow\mod\llz$. Precisely, let $e$ be
an extending vertex of $Q$ (i.e. the $e$-th entry of $\delta$ equals
$1$), let $P=P(e)$ be the indecomposable projective $\llz$-module
corresponding to $e$, and let $L$ be the unique indecomposable
preprojective $\llz$-module with dimension vector $\delta+\udim P$.
Let $\fk{C}(P,L)$ be the smallest full subcategory of $\mod\llz$
which contains $P$ and $L$ and is closed under taking extensions,
kernels of epimorphisms, and cokernels of monomorphisms in the
category of $\llz$-modules. Then $\fk{C}(P,L)$ is equivalent to
$\mod K$ (see for example~\cite{LXZ} Section 6.1). This embedding is
essentially independent of the field $\bbf_q$.


\section{Singular Ringel--Hall
algebras}\label{s:singular-hall-alg}

Let $Q$ be a tame quiver with vertex set $I$, and $\Lambda=\bbf_q Q$ be the path algebra of
$Q$ over the finite field $\bbf_q$. In this section we will define the singular Ringel--Hall
algebra of $\Lambda$ and give a description in terms of a set of generators and relations.

\subsection{The singular Ringel--Hall algebra $\ch^{s}(\llz)$}
Let $\ct_1,\ldots,\ct_l$ be the non-homogeneous tubes of $\mod\llz$ and assume that they are
respectively of period $r_1,\ldots,r_l$ (see Section~\ref{ss:tame-quiver}).

Let $\ch^*(\llz)$ be the twisted Ringel--Hall algebra of $\Lambda$,
which has a basis $\{u_{[M]}|M\in \mod\Lambda\}$ with structure
constants given by Hall numbers, see Section~\ref{ss:hall-alg}. The
\emph{singular Ringel--Hall algebra} of $\llz$, denoted by
$\ch^s(\llz)$, is defined as the subalgebra of $\ch^*(\llz)$
generated by $\{u_i,u_{[M]}:i\in I, M \in \ct_j, 1\leqslant
j\leqslant l\}$. Later we will prove the existence of Hall
polynomials, so that we have a generic version of the singular
Ringel--Hall algebra $\ch^s(\llz)$. We now set $\mathscr{D}^s(\llz)$
to be the subalgebra of $\mathscr{D}(\llz)$ generated by
$\{u_i^{\pm},u_{[M]}^{\pm}, K_\mu:i\in I, M \in \ct_j, 1\leqslant
j\leqslant l, \mu\in\bbn[I]\}.$ Namely, it is the reduced Drinfeld Double of
$\ch^s(\llz)$.

\begin{lemma}\label{l:2.1.2}
The $\mathscr{D}^s(\llz)$ is a Hopf algebra over $\bbq(v)$.
\end{lemma}

\begin{proof} We prove that $\mathscr{D}^s(\llz)$ is a Hopf subalgebra of $\mathscr{D}(\llz)$, i.e.
$\mathscr{D}^s(\llz)$ is closed under comultiplication and is closed
under antipode $S$. The former statement follows easily from
Lemma~\ref{l:2.1.1}. For the latter, it is sufficient to check on the
generators of $\mathscr{D}^s(\llz)$, since $S$ is an algebra
homomorphism. That $S(u_i^{\pm}),i\in I$, and
$S(K_\mu),\mu\in\mathbb{N}[I]$ belongs to $\mathscr{D}^s(\llz)$
follows immediately from the definition of $S$. Let $M\in \ct_j,
1\leq j\leq l$. We will prove that
$S(u^+_{[M]})\in\mathscr{D}^s(\llz),$ for $M\in\ct_i, 1\leqslant
i\leqslant l$ by induction on $\udim M$. The proof for $u^-$ is the
same.

Applying the equality $\mu(S \otimes
1)\bigtriangleup=\eta\epsilon$ to $u^+_{[M]}$, we obtain
\begin{eqnarray}\label{f:2.1.1}
S(u^+_{[M]})&=&-S(K_{\udim M})u^+_{[M]}\\
&&\mbox{}-\sum_{M_1,M_2\neq 0}v^{\lr{\udim M_1,\udim M_2}}
\frac{a_{M_1}a_{M_1}}{a_{M}}g^M_{M_1M_2}S(u^+_{[M_1]}K_{\udim
M_2})u^+_{[M_2]}.\nonumber\end{eqnarray} Assume $M_1\cong I_1\oplus
N_1,M_2\cong P_2\oplus N_2$, as in Lemma~\ref{l:2.1.1}. Note that
$\mathrm{Ext}^1(N_1,I_1)=0$. Therefore, we have
$u^+_{[N_1]}u^+_{[I_1]}=v^{\lr{\udim N_1 ,\udim I_1}}u^+_{[M_1]}$.
By induction hypothesis, we have
$S(u^+_{[N_1]})\in\mathscr{D}^s(\llz)$.  By Lemma~6.1 and 6.2 in
\cite{LXZ}, we have  $S(u^+_{[I_1]})\in\mathscr{D}^s(\llz)$. So
$S(u^+_{[M_1]})\in\mathscr{D}^s(\llz)$. Similarly,
$S(u^+_{[M_2]})\in\mathscr{D}^s(\llz)$. Therefore $(\ref{f:2.1.1})$
implies that $S(u^+_{[M]})\in\mathscr{D}^s(\llz)$.  
\end{proof}

\subsection{Decomposition of $\ch^{s}(\llz)$}\label{ss:decomposition} In this subsection, we
follow an idea of Sevenhant and Van den Bergh to obtain subalgebras
of $\ch^s(\llz)$ and $\mathscr{D}^s(\llz).$ (See also \cite{HX}.)

The twisted Ringel--Hall algebra $\ch^*(\Lambda)$ is naturally
$\mathbb{N}[I]$-graded, and so are its subalgebras $\ch^s(\llz)$ and
$\cc(\llz)$. We define a partial order on $\mathbb{N}[I]$: for
$\alpha,\beta\in\mathbb{N}[I]$, $\alpha\leq\beta$ if and only if
$\beta-\alpha\in\mathbb{N}[I]$. Clearly,
$\mathcal{C}(\llz)_\bz=\ch^s(\llz)_\bz$ if $\bz<\dz$.

Recall from Section~\ref{ss:double-hall-alg} that
$\varphi:\ch^+(\llz)\times\ch^-(\llz){\lra}\bbq(v)$ is a
non-degenerate bilinear form. It is easy to see that the restriction
of $\varphi$ on $\ch^{s,+}_\az\times\ch^{s,-}_\az$ is also
non-degenerate for all $\az\in \bbn[I]$. We now define
 $$\mathcal {L}^\pm_\dz=\{x^\pm\in\ch^{s,\pm}_{\dz}(\llz)|~\varphi(x^\pm,\mathcal {C}^{\mp}(\llz))=0\}.$$
 The non-degeneracy of $\varphi$ implies
\[\ch^{s,\pm}(\llz)_\dz=\mathcal {C}^{\pm}(\llz)_\dz\oplus\mathcal
{L}^\pm_\dz.\]
Let $\mathscr{D}^s(1)$ be the subalgebra of $\mathscr{D}^s$
generated by $\mathcal {C}^{\pm}(\llz)$ and $\mathcal{L}^{\pm}_\dz.$
Then we have a triangular decomposition
$$\mathscr{D}^s(1)=\mathscr{D}^s(1)^-\otimes\mathcal
{T}\otimes\mathscr{D}^s(1)^+.$$

Suppose $\mathcal {L}^{\pm}_{(m-1)\dz}$ and
$\mathscr{D}^s(m-1)^{\pm}$ have been defined, we inductively define
$\mathcal {L}^{\pm}_{m\dz}$ as follows:
\[
\mathcal{L}^{\pm}_{m\dz}=\{x^{+}\in\ch^{s,\pm}_{m\dz}|~\varphi(x^{\pm},\mathscr{D}^s(m-1)^\mp)=0\}.
\]
Let $\mathscr{D}^s(m)$ be the subalgebra of $\mathscr{D}^s$
generated by $\mathscr{D}^s(m-1)^{\pm}$ and $\mathcal{L}^{\pm}_{m\dz}.$ As
in the $m=1$ case, we have a triangular decomposition
$$\mathscr{D}^s(m)=\mathscr{D}^s(m)^-\otimes\mathcal
{T}\otimes\mathscr{D}^s(m)^+.$$
In this way we obtain a chain of subalgebras of $\mathscr{D}^s$
\[\cc(\llz)\subset \mathscr{D}^s(1)\subset\mathscr{D}^s(2)\subset\ldots\subset\mathscr{D}^s(m)\subset\ldots\subset\mathscr{D}^s\]
such that $\mathscr{D}^s=\bigcup_{m\in\bbn}\mathscr{D}^s(m)$.


\begin{lemma}\label{l:2.2.1}
Let $\eta_{n\dz}=\dim_{\bbq(v)}\mathcal
{L}^+_{n\dz}=\dim_{\bbq(v)}\mathcal {L}^-_{n\dz}$.  Then
$\eta_{n\dz}=l$.
\end{lemma}

\begin{proof} By proposition~6.5 in~\cite{LXZ},
we know that $\mathcal {C}^+(\llz)/(v-1)$ has a basis consiting of the following elements
\begin{itemize}
\item[a)] $\Psi(u_{[M(\az)]})$ for
 $\az\in \Phi^+_{Prep}$;
\item[b)] $\Psi(u_{\az,i})$ for $\az\in\ct_i, $ the  real
 roots,  $i=1,\ldots,l;$
\item[c)] $\Psi(u_{j,m\dz,i}-u_{j+1,m\dz,i}),$ $m\geq 1,$
 $1\leq j\leq r_i-1,$ $i=1,\ldots,l;$
\item[d)] $\Psi(\te_{n\dz}), n\geq 1$
\item[e)] $\Psi(u_{[M(\bz)]})$ for $\bz\in \Phi^+_{Prei}$.
\end{itemize}
While according to the definition $\ch^s$ the space  $\ch^s/(v-1)$ has a basis consiting of elements in
a) b) d) e) and all $\Psi(u_{j,m\dz,i}),$ $m\geq 1,$
 $1\leq j\leq r_i,$ $i=1,\ldots,l$. Thus we have $\eta_{n\dz}=l$ by the
 construction of $\mathcal
{L}^+_{n\dz}$ for all $n.$ 
\end{proof}

For each $n\dz,$ there exists a basis $\{x_n^1,\ldots,x_n^l\}$ of
$\mathcal {L}^+_{n\dz}$ and a basis $\{y_n^1,\ldots,y_n^l\}$ of
$\mathcal {L}^-_{n\dz}$ such that
$$\varphi(x_n^p,y_n^q)=\frac{1}{v-v^{-1}}\dz_{pq}.$$
Here $\dz_{pq}$ denotes the Kronecker's delta. Then we have
$$x_n^py_m^q-y_m^qx_n^p=\frac{K_{n\dz}-K_{-n\dz}}{v-v^{-1}}\dz_{pq}\dz_{mn},$$
for all $m,n\in\mathbb{N},1\leqslant p,q\leqslant l$
(see~\cite{HX}).

Let $J=\{(n\dz,p):n\in\mathbb{N},1\leqslant p\leqslant l\}.$ Define
\begin{eqnarray*}
\theta_i&=&\begin{cases} i & \text{ if } i\in I,\\ n\dz & \text{ if
} i=(n\dz,p)\in J,\end{cases}\\
x_i&=&\begin{cases} u_i^+ & \text{ if } i\in I,\\ x_n^p & \text{ if
}
i=(n\dz,p)\in J,\end{cases}\\
y_i&=&\begin{cases} -v^{-1}u_i^- & \text{ if } i\in I,\\-v^{-1}y_n^p
& \text{ if } i=(n\dz,p)\in J.\end{cases}\end{eqnarray*} By a
theorem of Sevenhant and Van den Bergh, we obtain that
$\mathscr{D}^s$ is generated by
$$\{x_i,y_i|~i\in I\cup J\}\cup\{K_\mu:\mu\in\bbn[I]\}$$
with the defining relations
\begin{eqnarray}
&K_0 =1,K_\mu K_\nu=K_{\mu+\nu} \text{ for all } \mu,\nu\in\bbn[I]&
\label{f:2.2.1}\\
&K_\mu x_i = v^{(\mu,\theta_i)}x_i K_\mu, \text{ and } K_\mu y_i =
v^{-(\mu,\theta_i)}y_i K_\mu \text{ for all } i\in I\cup J , \mu\in
\bbn[I]&\label{f:2.2.2}\\
&x_iy_j-y_jx_i=\frac{K_{\theta_i}-K_{-\theta_i}}{v-v^{-1}}\dz_{ij}
\text{ for all }i,j\in I\cup J& \label{f:2.2.3}\\
&\sum_{p+p'=1-a_{ij}}(-1)^px_i^{(p)}x_jx_i^{(p')}=0 \text{ and}
\sum_{p+p'=1-a_{ij}}(-1)^py_i^{(p)}y_jy_i^{(p')}=0& \label{f:2.2.4}\\
&x_ix_j =x_jx_i \text{ and } y_iy_j=y_jy_i \text{  for all } i,j\in
I\cup J \text{ with } (\theta_i,\theta_i)=0.& \label{f:2.2.5}
\end{eqnarray}
Here for an element $x$ and a positive integer $p$ the symbol $x^{(p)}$ denotes the divided power $\frac{x^p}{[p]!}$, where  $[p]=\frac{v^p-v^{-p}}{v-v^{-1}}$ and $[p]!=[1][2]\cdots[p]$.


Applying the relations above, the next statement is clear.

\begin{proposition}\label{p:2.2.2}
\begin{itemize}
\item[(a)] The elements $\{x_n^i|~n\in\bbn, 1\leqslant i\leqslant l\}$ are central in $\ch^{s,+}$.
Dually, the elments $\{y_n^i|~n\in\bbn, 1\leqslant i\leqslant l\}$ are central in $\ch^{s,-}$.
\item[(b)] We have an algebra isomorphism $\ch^{s,+}\cong \mathcal {C}^+\otimes \bbq(v)[x_n^i|~n\in\bbn, 1\leqslant i\leqslant l$.
\item[(c)] The element $x_n^j$ commutes with $\ch^{s,-}$, and the element $y_n^j$ commutes with $\ch^{s,+}$.
\end{itemize}
\end{proposition}


In particular, when the quiver $Q$ is of type
$\widetilde{A}_{n-1,1}$ for some $n\geq 3$, we have $l=1$ and we
have an algebra isomorphism $\ch^{s,+}\cong
U_q(\widetilde{sl}_n)^+\otimes \bbq(v)[x_1,\ldots,x_n,\ldots]$. This
means that $\ch^{s,+}$ is isomorphic to the $q$-Fock space of type
$A$.  Also this proposition implies that the algebra structure
$\ch^s$ not only  depends on the type, but also on the orientation.




\section{The quantum extended Kac--Moody algebra}

\subsection{} Let $(I,(,))$ be a datum in the sense of Green [G,
3.1]. Let $J$ be an index set and let
$\theta:J\rightarrow\bbn[I]\backslash\{0\},j\mapsto \theta_j$ be a
map with finite fibers. We denote by $\widetilde{G}$ the triple $(I,
(,),\{\theta_j: j\in J\})$, and call it an \emph{extended Green
datum}. We extend $\theta$ to a map $\theta:I\cup J\rightarrow\bbn
[I]\backslash\{0\}$ by setting $\theta_i=i$ for $i\in I$.

Let $Q$ be a tame quiver with vertex set $I$, and let $\llz$ be the path algebra. Recall that a Hopf subalgebra $\mathscr{D}^s(\llz)$ of $\mathscr{D}(\llz)$ is defined in Section~\ref{s:singular-hall-alg}.
Let $\widetilde{G}$ be the extended Green datum
corresponding to $\mathscr{D}^s(\Lambda)$. Precisely, $J=\{(n\dz,p):1\leqslant p\leqslant l\}$ and
$\theta_{(n\delta,p)}=n\delta$. From $\widetilde{G}$ we define a new
datum $\widetilde{G}'=(I\cup J, (,)'),$ where
$(i,j)'=(\theta_i,\theta_j)$ for all $i,j\in I\cup J.$ Let
$\mathscr{D}'(\Lambda)$ be the reduced Drinfeld double corresponding
to $\widetilde{G}'$. We have the following proposition similar to~\cite[Proposition 3.8]{DX}.

\begin{proposition}\label{p:3.1} There exists a surjective Hopf algebra
homomorphism $F:\mathscr{D}'(\Lambda)\rightarrow
 \mathscr{D}^{s}(\Lambda)$ such that $F(x_i^{\pm})=x_i^{\pm}$ and
 $F(K_i)=K_{\theta_i}$ for $i\in I\cup J,$ and $\ker F$ is the ideal
 generated by  $\{K_j-K_{\theta_j}| j\in J\}.$
 \end{proposition}

In particular, we have ${\mathscr{D}'}^{> 0}\cong \ch^{s,+}. $

\subsection{} Let
$\mathbb{U}=\mathbb{U}^-\otimes\mathbb{U}^0\otimes\mathbb{U}^+$ be
the quantized enveloping algebra in the sense of Drinfeld and Jimbo
with the generators $\{E_i,F_i,K_i,K_{-i}| i\in I\cup J\}$ subject
to relations similar to (\ref{f:2.2.1})--(\ref{f:2.2.5}) as in Section~\ref{ss:decomposition}. $\mathbb{U}$
is called the \emph{quantum extended Kac--Moody algebra}. Similar to~\cite[Corollary 8.3]{DX}, we have
$$\mathbb{U}\cong \mathscr{D}'.$$ So we obtain

\begin{proposition}\label{p:3.2}
 There exists an  algebra
 isomorphism $G:\mathbb{U}^+\stackrel{\sim}{\rightarrow} \ch^{s,+}$ such that $G(E_i)=x_i^{+}$
 for $i\in I\cup J.$
\end{proposition}




\section{PBW-basis of $\mathbb{U}^+$}\label{s:4}

\subsection{}\label{ss:4.1}
Let $Q$ be a connected tame quiver without oriented cycles with vertex set $I$, and let $\llz=\bbf_q Q$ be the path algebra of $Q$ over the finite field $\bbf_q$. Recall that the regular part of $\mod\llz$ is the direct sum of a family of homogeneous tubes and finitely many non-homogeneous tubes $\ct_1,\ldots,\ct_l$.

\begin{lemma}\label{l:regular-commutative}
 Let $M$ be a module in the direct sum of the homogeneous tubes, and $N$ be a regular module. Then $u_{[M]}*u_{[N]}=u_{[N]}*u_{[M]}$. If in addition no direct summands of $N$ belongs to the same tube as any indecomposable dirct summand of $M$, then $u_{[M]}*u_{[N]}=u_{[N]}*u_{[M]}=u_{[M\oplus N]}$.
\end{lemma}
\begin{proof}
 If no direct summands of $N$ belongs to the same tube as any indecomposable dirct summand of $M$, the statement follows from Lemma~\ref{l:2.1.1} (c) and the fact that the dimension vector of $M$ is a multiple of $\delta$, which lies in the radical of the Euler form. If $N$ is also in the direct sum of the homogeneous tubes, the statement is proved in~\cite{Zh}. Generally we write $u_{[N]}=u_{[N_1]}*u_{[N_2]}$ with $N_1\oplus N_2\cong N$ such that $N_1$ belongs to the direct sum of the homogeneous tubes and $N_2$ belongs to the direct sum of the non-homogeneous tubes. Then
\begin{eqnarray*}
u_{[M]}*u_{[N]}&=&u_{[M]}*u_{[N_1]}*u_{[N_2]}\hspace{7pt}=\hspace{7pt}u_{[N_1]}*u_{[M]}*u_{[N_2]}\\
&=&u_{[N_1]}*u_{[N_2]}*u_{[M]}\hspace{7pt}=\hspace{7pt}u_{[N]}*u_{[M]}.
\end{eqnarray*}
\end{proof}

Let $P$, $L$, $K$ and $\fk{C}(P,L)$ be as in Section~\ref{ss:tame-quiver}, and let $F:\mod K\cong\fk{C}(P,L)\hookrightarrow\mod\llz$ be the exact embedding.
It gives rise to an injective homomorphism of algebras, still denoted
by $F:\ch^*(K)\hookrightarrow\ch^*(\llz).$ In $\ch^*(K)$ a distinguished element $E_{m\dz_K}$ is defined for any $m\geq 1$, and we set $E_{m\dz}=F(E_{m\dz_K}),$ see~\cite{LXZ} for more details. Since $E_{m\dz_K}\in
\cc^*(K),$ and $\lr{L}, \lr{P}\in \cc^*(\llz)$, it follows that $E_{m\dz}$ is in
$\cc^*(\llz)$ and even  in $\cc^*(\llz)_{\cz}$, where $\cz=\bbq[v,v^{-1}]$.  Let $\cal K$ be the
subalgebra of $\cc^*(\llz)$ generated by $E_{m\dz}$ for $m\in\bbn,$
it is a polynomial ring on infinitely many variables
$\{E_{m\dz}|m\geq1\},$ and its integral form is the polynomial ring
on variables $\{E_{m\dz}|m\geq1\}$ over $\cz.$

We denote by $\fk{C}_0$ (respectively, $\fk{C}_1$) the full subcategory of
$\fk{C}(P,L)$ consisting of the $\Lambda$-modules which belong to
homogeneous (respectively, non-homogeneous) tubes of $\mod\Lambda.$

We now decompose $E_{n\dz}$ as
   $$ E_{n\delta}=E_{n\delta,1}+E_{n\delta,2}+E_{n\delta,3} , $$
 where
\begin{eqnarray}
   E_{n\delta,1}&=&v^{-n|\delta|}
   \sum_{[M]:M\in\fk{C}_1,\udim M=n\delta}u_{[M]}\label{f:4.1.1}\\
    E_{n\delta,2}&=&v^{-n|\delta|}
   \sum_{\substack{[M]:M\in\fk{C},\udim M=n\delta\\
   M=M_1\oplus M_2, 0\neq M_1\in \fk{C}_1, 0\neq M_2\in
   \fk{C}_0}}u_{[M]}\label{f:4.1.2}\\
   E_{n\delta,3}&=&
   v^{-n|\delta|}\sum_{[M]:M\in\fk{C}_0, \udim M=n\delta}u_{[M]},\label{f:4.1.3}
   \end{eqnarray}
where $|\delta|$ is the sum of all entries of $\delta$.

\begin{lemma}\label{l:E-n-delta} Let $n,n'$ be two positive integers. Then we have
\begin{itemize}
 \item[(a)] $E_{n\delta,1}*E_{n'\delta,3}=E_{n'\delta,3}*E_{n\delta,1}$;
 \item[(b)] $E_{n\delta,2}=\sum_{m=1}^{n-1} E_{m\delta,1}*E_{(n-m)\delta,3}$;
 \item[(c)] $E_{n\delta,3}*E_{n'\delta,3}=E_{n'\delta,3}*E_{n\delta,3}$.
\end{itemize}
\end{lemma}
\begin{proof} (a) and (c) follows from Lemma~\ref{l:regular-commutative}. (b) holds because
 \begin{eqnarray*}
E_{n\delta,2}
&=&v^{-n|\delta|}
   \sum_{\substack{[M_1\oplus M_2]:\udim (M_1\oplus M_2)=n\delta,\\
   0\neq M_1\in \fk{C}_1, 0\neq M_2\in
   \fk{C}_0}}u_{[M_1\oplus M_2]}\\
&=&v^{-n|\delta|}
   \sum_{\substack{[M_1\oplus M_2]:\udim (M_1\oplus M_2)=n\delta,\\
   0\neq M_1\in \fk{C}_1, 0\neq M_2\in
   \fk{C}_0}}u_{[M_1]}*u_{[ M_2]}\\
&=&v^{-n|\delta|}\sum_{m=1}^{n-1}
   \sum_{\substack{[M_1]:M_1\in\fk{C}_1,\\ \udim M_1=m\delta}}
   \sum_{\substack{[M_2]:M_2\in\fk{C}_0,\\ \udim M_2=(n-m)\delta}}u_{[M_1]}*u_{[ M_2]}\\
&=&\sum_{m=1}^{n-1}
   (\sum_{\substack{[M_1]:M_1\in\fk{C}_1,\\ \udim M_1=m\delta}}v^{-m|\delta|}u_{[M_1]})
   *(\sum_{\substack{[M_2]:M_2\in\fk{C}_0,\\ \udim M_2=(n-m)\delta}}v^{-(n-m)|\delta|}u_{[M_2]})\\
&=&\sum_{m=1}^{n-1} E_{m\delta,1}*E_{(n-m)\delta,3}.
\end{eqnarray*}
\end{proof}

For a partition
$\textbf{w}=(w_1,\ldots,w_t)$ of $n,$ we define
\begin{eqnarray*}
E_{\textbf{w}\dz}&=&E_{w_1\dz}*\cdots*E_{w_t\dz}\\
E_{\textbf{w}\dz,3}&=&E_{w_1\dz,3}*\cdots*E_{w_t\dz,3}.\end{eqnarray*}
Let ${\bf P}(n)$ be the set of all partitions of $n,$ and
$\lr{N}=v^{-\dim N+\dim \ed(N)}u_{[N]}.$ Set
$${\bf
B}=\{\lr{P}*\lr{M}*E_{\textbf{w}\dz,3}*\lr{I}|~ P\in\cp_{prep},
M\in\bigoplus_{i=1}^l\ct_i, I\in\cp_{prei},{\bf w}\in{\bf P}(n),
n\in\bbn\},$$
where recall that $\cp_{prep}$ respectively $\cp_{prei}$ is the set of isomorphism classes of preprojective respectively preinjective $\llz$-modules, and $\ct_1,\ldots,\ct_l$ are the non-homogeneous tubes of $\mod\llz$.
Then we have the following:


\begin{theorem}\label{t:4.1.1}  The set ${\bf B}$ is a
$\bbq(v)-$basis of $\ch^{s,+}$ ($\cong\mathbb{U}^+$).
\end{theorem}

\begin{proof}
We have by definition that $E_{n\delta}$ and $E_{n\delta,1}$ belong
to $\ch^{s,+}$. Then it follows from Lemma~\ref{l:E-n-delta} (b) by
induction on $n$ that both $E_{n\delta,2}$ and $E_{n\delta,3}$
belong to $\ch^{s,+}$. As a consequence, the set ${\bf B}$ is contained in
$\ch^{s,+}$. Because ${\bf B}$ is linear independent over $\bbq(v)$,
it remains to show that ${\bf B}$ linearly spans $\ch^{s,+}$.

Let $\Pi_i^a$ be the set of aperiodic $r_i$-tuples
of partitions, for all $1\leqslant i\leqslant l.$ Set
\begin{eqnarray*}
  {\bf B}^c&=&\{\lr{P}\ast E_{\pi_1
}\ast\cdots\ast  E_{\pi_l }\ast E_{{\bf
w}\dz}\ast\lr{I}:\\
&&P\in\cp_{prep},I\in\cp_{prei},\pi_i\in\Pi_i^a,1\leqslant i\leqslant l, {\bf
w}\in{\bf P}(n),n\in\bbn\}.
 \end{eqnarray*}
By  Proposition~7.2 in \cite{LXZ}, we know that ${\bf B}^c$
 is a
$\bbq(v)$-basis of $\cc^*(\llz).$  By definition $\ch^{s,+}$ is
generated by ${\bf B}^c$ and $\{u_{[M]}:M\in\bigoplus_{i=1}^l
\ct_i\}$. The latter elements belong to ${\bf B}$, so now it remains
to prove that each element in ${\bf B}^c$ ,$u_{[M]}\ast\lr{P}$ and
$(E_{{\bf w}\dz}\ast\lr{I})\ast u_{[M]}$ can be linearly spanned by
${\bf B}$.

For ${\bf w}=(w_1,\ldots,w_t)\in {\bf P}(n)$, we have by Lemma~\ref{l:E-n-delta} (b)
\begin{eqnarray*}
E_{{\bf w}\delta}&=&E_{w_1\delta}*\cdots*E_{w_t\delta}\\
&=&\prod_{j=1}^t(E_{w_j\delta,1}+\sum_{m_j=1}^{w_j-1}E_{m_j\delta,1}*E_{(w_j-m_j)\delta,3}+E_{w_j\delta,3})
\end{eqnarray*}
By Lemma~\ref{l:E-n-delta} (a) (c), we can write $E_{{\bf w}\delta}$ as a linear combination of elements of the form
$E_{m_1\delta,1}*\cdots *E_{m_r\delta,1}*E_{m'_1\delta,3}*\cdots*E_{m'_{r'}\delta,3}$,
which itself is a linear combination of elements of the form $\langle M\rangle*E_{{\bf w'}\delta,3}$,
where $M\in\bigoplus_{i=1}^l \ct_i$ and ${\bf w'}$ is a partition.

Similarly, we can prove that each element of form
$u_{[M]}\ast\lr{P}$ and $(E_{{\bf w}\dz}\ast\lr{I})\ast u_{[M]}$ can
be linearly spanned by ${\bf B}$. This completes the proof.
\end{proof}

\section{Canonical bases of $\mathbb{U}^+$ ($\cong\ch^s(\llz)$)}

In this section we give the main theorem of this paper, that is, a description of the canonical basis of $\ch^s(\llz)\cong\mathbb{U}^+$.

\subsection{}
Let $Q$ be a tame quiver with vertex set $I$, and $\Lambda=\bbf_q Q$
be the path algebra of $Q$ over the finite field $\bbf_q$. We denote
by $M(x),x\in\bbe_\az,$ the $\llz$-module of dimension vector $\az$
corresponding to $x.$ For subsets $\ca\subset\bbe_{\az}$ and
$\cb\subset\bbe_{\bz},$ we define the extension set $\ca\star\cb$ of
$\ca$ by $\cb$  to be
\begin{eqnarray*} \ca\star\cb&=&\{z\in\bbe_{\az+\bz}|\ \text{there exists an exact
sequence}\ \\
&& \quad 0\ra M(x)\ra M(z)\ra M(y)\ra 0 \ \text{with}\
x\in\cb,\ y\in\ca\}.
\end{eqnarray*}
It follows from the definition that
$\ca\star\cb=p_3p_2(p_1^{-1}(\ca\times \cb))$, see
Section~\ref{ss:construction-lusztig} for the definitions of
$p_1,p_2$ and $p_3$. Because $p_1$ is a locally trivial fibration ,
we have $\overline{\ca}\star \overline{\cb} \subseteq
\overline{\ca\star\cb}$(see Lemma~2.3 in \cite{LXZ}) . In
particular, $\ol{\co}_M\star\ol{\co}_N=\ol{\co}_{M\oplus N}$ if
$\ext(M,N)=0,$ i.e. $\co_{M\oplus N}$ is open and dense  in
$\ol{\co}_M\star\ol{\co}_N.$

Set $\cdim\ca=\dim\bbe_{\az}-\dim\ca.$ We will need the following:


\begin{lemma}[\cite{Re}]\label{l:5.1.1}
Let  $\az, \bz\in\bbn I$. If $\ca\subset\bbe_{\az}$ and $\cb\subset\bbe_{\bz}$  are
irreducible algebraic varieties and are stable under the action of
$G_{\az}$ and $G_{\bz}$ respectively, then $\ca\star\cb$ is
irreducible and stable under the action of $G_{\az+\bz},$ too.
Moreover,
$$\cdim \ca\star\cb =\cdim\ca+\cdim\cb-\lr{\bz,\az}+r,$$
where $0\leq r\leq\min\{\dim_k\hom(M(y),M(x))|y\in\cb,x\in\ca\}$.
\end{lemma}

Recall that for $x\in \bbe_\az$ the symbol $\co_x$ (or $\co_{M(x)}$) denotes the $G_\az$-orbits of $x$.
We now introduce two orders in $\llz-$mod as follows:
\begin{itemize}
\item $N\leq_{deg} M$  if $\co_N\subseteq\overline{\co}_M$.
\item $N\leq_{ext} M$  if there exist $M_i,U_i,V_i$ and short exact sequence
\[0\lra U_i\lra M_i\lra V_i\lra 0\]
such that $M=M_1,M_{i+1}=U_i\oplus V_i, 1\leqslant
i\leqslant p$, and $N=M_{p+1}$ for $p\in\mathbb{N}$.
\end{itemize}
%

\begin{proposition}\label{p:5.1.2}\cite{Z}
 The orders $\leq_{deg},\leq_{ext}$ are equivalent in $\llz-mod.$
\end{proposition}


We denote by $\cp_X$ the category of perverse sheaves on an
algebraic variety $X.$  Let $f$ be a locally closed embedding from
$X$ to algebraic variety $Y$. One has the intermediate extension
functor
$$f_{!*}:\cp_X\lra \cp_Y, P\longmapsto \Im\{ ^p\mathcal{H}^0(f_!P)\lra {^p\mathcal{H}}^0(f_*P)\}.$$

Let $V$ be a locally closed, smooth, irreducible subvariety of $X$,
of dimension $d$ and let $\mathcal{L}$ be an irreducible
$\mathbb{Q}_l-$ local system on $V.$ Then $\mathcal{L}[d]$ is an
irreducible perverse sheaf on $V$ and there is a unique irreducible
perverse sheaf $\widetilde{\mathcal{L}}[d],$ whose restriction to
$V$ is $\mathcal{L}[d]$, we have
$\widetilde{\mathcal{L}}[d]=IC_{\overline{V}}(\mathcal{L}),$ where
$IC_{\overline{V}}(\mathcal{L})$ is the intersection cohomolgy
complex of Deligne -Goresky-Macpherson of $\overline{V}$ with
coefficients in $\mathcal{L}.$ The extension of
$\widetilde{\mathcal{L}}[d]$ to $X$ (by $0$ outside $\overline{V}$)
is an irreducible perverse sheaf on $X.$

 In particular, suppose $\mathcal{L}$ is a local system on a
nonsingular Zariski dense open subset $j:U\lra Y$ of the irreducible
$n-$dimensional $Y.$ Then
$IC_Y(\mathcal{L}):=j_{!*}\mathcal{L}[n]\in\cp_Y.$

\begin{definition}\label{d:5.1.3} \cite{KW} Let $\mathcal{K}$ be a Weil
complex. Then $\mathcal{K}$ is said to \emph{have the purity
property} if on all stalks of the semisimplication
$\mathscr{H}^v(\mathcal{K})_x^{ss}(v/2)$ of cohomology sheaves the
Frobenius $F_x$ acts trivially.
\end{definition}




 \begin{lemma}\label{l:5.1.4}  Let $Y$ be an irreducible algebraic variety of dimension $m$, and let $p$ be a smooth morphism $p:X\rightarrow Y$ of relative dimension $d$. Suppose $U_0$  is a nonsingular
Zariski dense open subset of $Y,$ $j_0:U_0\rightarrow Y$ is an open embedding. Then we have the following
  a cartesian square
\[
 \xymatrix{V=p^{-1}(U_0)\ar[d]^{p|_U}\ar[r]^(0.65)j &  X\ar[d]^p\\
 U \ar[r]^{j_0} &Y.}
\]
If $j_{0 !*}\overline{\bbq}_l[m]$ has the purity property, then $j_{
!*}\overline{\bbq}_l[d+m]$ has the purity property.
\end{lemma}

\begin{proof} By the definition of $j_{*}$ and $j_!,$ we have a
natural morphism
$$\varphi_Y:~~^p\mathcal{H}^0(j_{0!}\overline{\bbq}_l[m])\lra
{^p\mathcal{H}}^0(j_{0*}\overline{\bbq}_l[m]).$$ It induces an
intermediate extension functor
$$ j_{0!*}: \cp_U\lra\cp_Y$$
such that $j_{0!*}\overline{\bbq}_l[m]=\Im\{
^p\mathcal{H}^0(j_{0!}\overline{\bbq}_l[m])\lra
{^p\mathcal{H}}^0(j_{0*}\overline{\bbq}_l[m])\}. $ Furthermore,
$$p^*[d]\circ\varphi_Y: p^*[d](^p\mathcal{H}^0(j_{0!}\overline{\bbq}_l[m]))
\lra p^*[d](^p\mathcal{H}^0(j_{0*}\overline{\bbq}_l[m])).$$
 Since $p: X\lra Y$ is smooth of relative dimension
$d,$ $p^*[d]=p^{!}[-d]$ is $t$-exact , we have
$$\begin{array}{l}
\varphi_X: ~~^p\mathcal{H}^0(p^*j_{0!}\overline{\bbq}_l[m+d])\lra
{^p\mathcal{H}}^0(p^*j_{0*}\overline{\bbq}_l[m+d]),\\[3mm]
\varphi_X: ~~^p\mathcal{H}^0(j_{!}\overline{\bbq}_l[m+d])\lra
{^p\mathcal{H}}^0(j_{*}\overline{\bbq}_l[m+d]).
\end{array}$$

So $p^*[d](j_{0 !*}\overline{\bbq}_l[m])=j_{
!*}\overline{\bbq}_l[d+m].$ By the same argument, we get
$p^*[d]\circ F=F\circ p^*[d].$ Since $j_{0 !*}\overline{\bbq}_l[m]$
has the purity property,  the statement of the lemma is true.
\end{proof}

 Let $\cn_{w_i}$ and $\cn_{w_i,3}$ be respectively the union
of orbits of regular modules of $\fk{C}(P,L)$ and $\fk{C}_0(P,L)$
with dimension vector $w_i\dz.$ Set $\cn_{{\bf
w}}=\cn_{w_1}\star\cdots\star\cn_{w_t}$ and $\cn_{{\bf
w},3}=\cn_{w_1,3}\star\cdots\star\cn_{w_t,3}.$
For any $P\in\cp_{prep}, M\in\bigoplus_{i=1}^l\ct_i,
I\in\cp_{prei},\pi_i\in\Pi_i^a,1\leqslant i\leqslant l,{\bf
w}\in{\bf P}(n), n\in\bbn,$ we define the varieties
\begin{eqnarray*}
\co_{P,\pi_1,\ldots,\pi_l,{\bf w},I}&=&\co_P\star\co_{\pi_1}\star\cdots\star\co_{\pi_l}
\star\cn_{{\bf w}}\star\co_I ,\\
\co_{P,M,{\bf
w},I}&=&\co_P\star\co_M\star\cn_{{\bf w},3}\star\co_I.\end{eqnarray*}

According to \cite{L4}£¬ \cite{L5} and \cite{LXZ}, we know that
$IC_{\overline{\co}_{P,\pi_1,\ldots,\pi_l,{\bf
w},I}}(\overline{\bbq}_l)$ has the purity property. In order to
construct the canonical basis of $\ch^s(\llz)$
($\cong\mathbb{U}^+$), we need to study the purity property of
$\overline{\co}_{P,M,{\bf w},I}.$


\begin{theorem}\label{t:5.1.5} Assume $Q$ is not of type $\widetilde{E}_8$.
Let $X=\overline{\co}_{P,M,{\bf
w},I}$, for $P\in\cp_{prep}, M\in\bigoplus_{i=1}^l\ct_i,
I\in\cp_{prei}, {\bf w}\in {\bf P}(n), n\in\bbn.$ Then
$IC_X(\overline{\bbq}_l)$ has the purity property.
\end{theorem}

We postpone the proof of Theorem~\ref{t:5.1.5} to later sections.

\subsection{}\label{ss:canonical-basis}
Let
\begin{eqnarray}
b_{\overline{\co}_{P,M,{\bf w},I}}&=&
\sum_{i,N\in{\overline{\co}_{P,M,{\bf
w},I}}^F}(-1)^iv^{i+\dim\co_N}\dim
IH^i((IC_{\overline{\co}_{P,M,{\bf
w},I}}(\overline{\bbq}_l))_{N})\lr{N}.\label{f:7.2.1}
\end{eqnarray} Set
\begin{eqnarray*}{\bf CB}&=&\{b_{\overline{\co}_{P,M,{\bf w},I}}|P\in\mathcal {P}_{prep},
M\in\oplus_{i=1}^l\ct_i,I\in\mathcal {P}_{prei},{\bf w}\in {\bf
P}(n),n\in\bbn\}.\end{eqnarray*}
Then we have the main theorem of this paper

\begin{theorem}\label{t:7.2.1} For $Q$ not of type $\widetilde{E}_8$, the set ${\bf CB}$ is the
canonical basis of $\ch^s(\llz)$ ($\cong\mathbb{U}^+$).
\end{theorem}


\begin{proof}
By Proposition~\ref{p:1.3.1}, Theorem~\ref{t:4.1.1} and
Theorem~\ref{t:5.1.5}, the proof is complete.
\end{proof}




\section{Purity Properties of Perverse sheaves of closure
of semi-simple objects in $\ct_i$}

We start from this section the proof of Theorem~\ref{t:5.1.5}.
In this section we deal with the special case $\co_M$ for $M$ a semi-simple object in any non-homogeneous tube.
We proceed type by type.






\subsection{Type $\widetilde{A}_{n,1}$}

Let $Q$ be the quiver
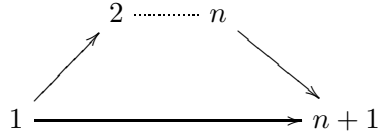
\[\xymatrix{&2\ar@{.}[r]&n\ar[dr]\\
1\ar[ur]\ar[rrr]&&&n+1}\]
There is only one non-homogeneous tube, which is of period $n$, and the regular simples $E_1,\ldots,E_n$ respectively have dimension vectors
$(1,0,0,\ldots,0,1),(0,1,0,\ldots,0,0),\ldots,(0,0,0,\ldots,1,0).$




\begin{lemma}\label{l:5.2.1}Let $X=\overline{\co}_{\bigoplus_{i=1}^nE_i}.$ Then
$IC_X(\overline{\bbq}_l)$ has the purity property.
\end{lemma}

In the following , the set of $m\times n$ matrices over $k$ is
denoted by $M_{mn}(k),$   the $k-$vector space generated by column
vectors of a matrix $A$ is denoted by $C(A)$, the  projective space
of a $k-$vector space $V$ is denoted by $P(V).$

\begin{proof}[Proof of Lemma~\ref{l:5.2.1}]
Let $V$ be the $I$-graded vector space
$V=\bigoplus_{i=1}^{n+1}V_i,$ with $V_i=k.$ Set
$\az=(1,1,\ldots,1)\in\bbn [I],$ then
$$\begin{array}{l}
\bbe_{\az}=\bbe_{V}=\{x=(x_{12},x_{23},\ldots,x_{n,n+1},x_{1,n+1})\in
k^{n+1}\},\\[5pt]
\co_{\bigoplus_{i=1}^nE_i}=\{x\in\bbe_{\az}|~x_{12}=x_{23}=\ldots=x_{n,n+1}=0,x_{1,n+1}\neq
0\},\\[5pt]
X=\overline{\co}_{\bigoplus_{i=1}^nE_i}=\{x\in\bbe_{\az}|~x_{12}=x_{23}=\ldots=x_{n,n+1}=0\}=\bba^1.
\end{array}$$
 It is clear that
$IC_{X}(\overline{\bbq}_l)$ has the purity property. The statement is proved.
\end{proof}

\subsection{Type $\widetilde{D}_n$}

Let $Q$ be the quiver
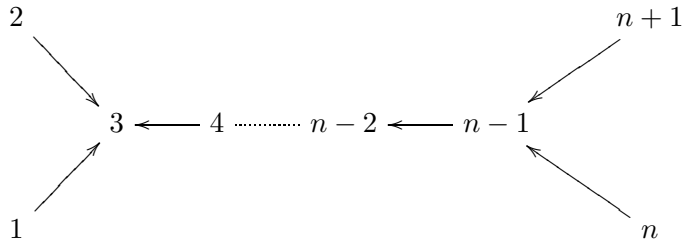
\[\xymatrix{ 2\ar[dr] &&&&& n+1\ar[dl]\\
&3&4\ar[l]\ar@{.}[r]& n-2 &n-1\ar[l]  &\\
 1\ar[ur] &&&&&n\ar[ul]}\]
There are three non-homogeneous tubes $\ct_1,\ct_2,\ct_3$, respectively of periods $2,2,n-2$. Let $E_1,E_2$ be the regular simples in $\ct_1$, let $E'_1,E'_2$ be the regular simples in $\ct_2$ and let $E''_1,\ldots,E''_{n-2}$ be the regular simples in $\ct_3$. Their dimension vectors are given as follows
\begin{itemize}
\item[$\ct_1$:] $(1,0,1,1\ldots,1,1,0)$, $(0,1,1,1,\ldots,1,0,1)$;
\item[$\ct_2$:] $(1,0,1,1\ldots,1,0,1)$, $(0,1,1,1,\ldots,1,1,0)$;
\item[$\ct_3$:] $(1,1,1,0\ldots,0,0,0)$, $(0,0,1,1,\ldots,1,1,1)$, $(0,0,0,1,0,\ldots,,0,0,0)$, $\ldots$, $(0,0,0,0,\ldots,0,1,0,0)$.
\end{itemize}


\begin{lemma}\label{l:5.2.2}
 Set
$X_1=\overline{\co}_{E_1\oplus E_2}, X_2=\overline{\co}_{E'_1\oplus
E'_2},$ and $X_3=\overline{\co}_{\bigoplus_{i=1}^{n-2} E_i''}.$
Then $IC_{X_i}(\overline{\bbq}_l)$ has
the purity property for any $i=1,2,3.$
\end{lemma}

\begin{proof}We only prove for $i=1.$ The other cases can be proved
similarly. Let $V$ be the $I$-graded vector space
$V=\bigoplus_{i=1}^{n+1}V_i,$ with $V_1=V_2=V_{n}=V_{n+1}=k$ and $V_i=k^2$
for all $3\leqslant i\leqslant n-1.$ Set
$\az=(1,1,2,\ldots,2,1,1)\in\bbn [I]$. Then
\begin{eqnarray*}
\bbe_{\az}=\bbe_{V}&\hspace{-8pt}=\hspace{-8pt}&\{x|~x=(x_{13},x_{23},x_{43},\ldots,x_{n-1,n-2},x_{n,n-1},x_{n+1,n-1})\\
                    && \in M_{21}(k)\times M_{21}(k)\times  M_{22}(k)\times\cdots\times M_{22}(k)\times M_{21}(k)\times M_{21}(k) \}.
\end{eqnarray*}
Thanks to \cite{[DR]}, we have $E_1\oplus E_2=M(y)$ for
$y\in\bbe_{\az}$ given by
\[y_{13}=y_{n,n-1}=\left(\begin{array}{l}
1\\
0\end{array}\right),~~y_{23}=y_{n+1,n-1}=\left(\begin{array}{l}
0\\
1\end{array}\right),~~y_{43}=\ldots=y_{n-1,n-2}=\left(\begin{array}{ll}
1&0\\
0&1\end{array}\right).\]


We claim that
\begin{eqnarray*}
\co_{E_1\oplus E_2}&=&\{x\in\bbe_{\az}|(x_{13}~~~x_{23})=x_{43}\cdots x_{n-1,n-2}\cdot(x_{n,n-1}~~~x_{n+1,n-1})\cdot
\left(\begin{array}{ll}
\xi_1&0\\
0&\xi_2\end{array}\right),\\
&&\xi_1\xi_2\neq 0,\det(x_{i,i-1})\neq 0, 4\leqslant i\leqslant
n-1,\det(x_{n,n-1}~~~x_{n+1,n-1})\neq 0\}
\\
&=&k^*\times
k^*\times\{x=(x_{43},\cdots,x_{n-1,n-2},x_{n,n-1},x_{n+1,n-1})|\\
&&\det(x_{i,i-1})\neq 0, 4\leqslant i\leqslant
n-1\,
\det(x_{n,n-1}~~~x_{n+1,n-1})\neq 0\}.
\end{eqnarray*}
Let $S$ be the set on the right hand side of the equality. For any $x\in\co_{E_1\oplus E_2}$ there exists $g=(g_i)_{i\in I}\in
GL_{\az}$ such that $x=g\bullet y$, i.e.
$$x_{13}=g_3\left(\begin{array}{l}
1\\
0\end{array}\right) g_1^{-1},~~x_{23}=g_3\left(\begin{array}{l}
0\\
1\end{array}\right) g_2^{-1},~~x_{i,i-1}=g_{i-1}g_i^{-1},4\leqslant
i\leqslant n-1,$$ $$x_{n,n-1}=g_{n-1}\left(\begin{array}{l}
1\\
0\end{array}\right) g_n^{-1},~~
  x_{n+1,n-1}=g_{n-1}\left(\begin{array}{l}
0\\
1\end{array}\right) g_{n+1}^{-1}.$$
The inclusion $\co_{E_1\oplus E_2}\subseteq S$ follows immediately, with $\xi_1=g_n g_1^{-1}$ and $\xi_2=g_{n+1}g_2^{-1}$.
%
Conversely, for an element in $S$, we have $x=g\bullet y$ for $g=(g_i)_{i\in I}\in GL_{\alpha}$ with $g_1=\xi_1^{-1}$, $g_2=\xi_2^{-1}$, $g_i=x_{i+1,i}\cdots x_{n-1,n-2} (x_{n,n-1}~~~x_{n+1,n-1})$ for $3\leq i\leq n-2$, $g_{n-1}=(x_{n,n-1}~~~x_{n+1,n-1})$ and $g_n=1$, $g_{n+1}=1$.

Let
\[
X'=\{(x_{43},\ldots,x_{n-1,n-2},x_{n,n-1},x_{n+1,n-1})\in M_{22}(k)\times \cdots\times M_{22}(k)\times M_{21}(k)\times M_{21}(k)\}.
\]
Then
$$
X=\overline{\co}_{E_1\oplus E_2} =k\times k\times
X'=\{x|x=(\xi_1,\xi_2,x'), \xi_1,\xi_2\in k, x'\in
X'\}=\bba^{4n-10}.
$$
The result is clear now.
\end{proof}

\subsection{Type $\tilde{E}_6$}

Let $Q$ be the quiver
\[\xymatrix{&&7\ar[d]&&\\ &&6\ar[d]&&\\ 1\ar[r]&2\ar[r]&3&4\ar[l]&5\ar[l]}\]
There are three non-homogeneous tubes $\ct_1,\ct_2,\ct_3$, respectively of periods $2,3,3$. Let $E_1,E_2$ be the regular simples in $\ct_1$, let $E'_1,E'_2,E'_3$ be the regualr simples in $\ct_2$, and let $E''_1,E''_2,E''_3$ be the regular simples in $\ct_3$. Their dimension vectors are given as follows
\begin{itemize}
\item[$\ct_1$:] $(1,1,2,1,1,1,1)$, $(0,1,1,1,0,1,0)$;
\item[$\ct_2$:] $(1,1,1,1,0,0,0)$, $(0,1,1,0,0,1,1)$, $(0,0,1,1,1,1,0)$;
\item[$\ct_3$:] $(1,1,1,0,0,1,0)$, $(0,1,1,1,1,0,0)$, $(0,0,1,1,0,1,1)$.
\end{itemize}



\begin{lemma}\label{l:5.2.3}
Set $X_1=\overline{\co}_{E_1\oplus E_2},
X_2=\overline{\co}_{E'_1\oplus E'_2\oplus E'_3},$ and
$X_3=\overline{\co}_{E''_1\oplus E''_2\oplus E''_3}.$
 Then $IC_{X_i}(\overline{\bbq}_l)$ has the
purity property for any $i=1,2,3$.
\end{lemma}

\begin{proof} We only prove for $i=1.$ The other cases can be proved
similarly. Let $V$ be the $I$-graded vector space
$V=\bigoplus_{i=1}^{7}V_i$ with $V_1=V_5=V_7=k,
V_2=V_4=V_6=k^2$ and $V_3=k^3$. Set $\az=(1,2,3,2,1,2,1)\in\bbn [I],$ then
\begin{eqnarray*}
\bbe_{\az}=\bbe_{V}&\hspace{-8pt}=\hspace{-8pt}&\{x|~x=(x_{12},x_{23},x_{43},x_{54},x_{63},x_{76}),\\
                    && x_{12},x_{54},x_{76}\in M_{21}(k),
x_{23},x_{43},x_{63}\in  M_{32}(k) \}.
\end{eqnarray*}
Thanks to \cite{[DR]}, we have $E_1\oplus E_2=M(y)$ for
$y\in\bbe_{\az}$ given by
\[y_{12}=y_{54}=y_{76}=\left(\begin{array}{l}
1\\
0\end{array}\right),~~ y_{23}=\left(\begin{array}{ll}
1&0\\
0&0\\
0&1\end{array}\right),~~y_{43}=\left(\begin{array}{ll}
0&0\\
1&0\\
0&1\end{array}\right),~~y_{63}=\left(\begin{array}{ll}
1&0\\
1&0\\
0&1\end{array}\right).\]
%
%
By definition any $x\in\co_{E_1\oplus E_2}$ is of the form $x=g\bullet y$ for some $g=(g_i)_{i\in I}\in
GL_{\az}$, i.e.
\begin{eqnarray} \label{f:5.2.6}& x_{12}=g_2\left(\begin{array}{l} 1\\0
\end{array}\right)g_1^{-1},~~
x_{23}=g_3\left(\begin{array}{ll}
1&0\\
0&0\\
0&1\end{array}\right)g_2^{-1},~~ x_{54}=g_4\left(\begin{array}{l} 1\\0
\end{array}\right)g_5^{-1},& \\
&x_{43}=g_3\left(\begin{array}{ll}
0&0\\
1&0\\
0&1\end{array}\right)g_4^{-1},~~ x_{76}=g_6\left(\begin{array}{l} 1\\0
\end{array}\right)g_7^{-1},~~
x_{63}=g_3\left(\begin{array}{ll}
1&0\\
1&0\\
0&1\end{array}\right)g_6^{-1}.& \nonumber
\end{eqnarray}

Let $Z$ be the subvariety of $k^2\times {\mathbb{P}}^1\times k\times
P(k^2\times k^2\times k^2)\times \bbe_{\az}$ consisting of those
elements $z=(z_0,[\lz_1:\lz_2],z_2,[z_3,z_4,z_5],x)$ such that
\begin{eqnarray*}
 x_{12}\lz_1=z_0\lz_2, x_{23}z_0+x_{43}x_{54}z_2=x_{63}x_{76},~~x_{23}z_3=x_{43}z_4=x_{63}z_5, z_2\neq 0&\\
\det (z_0~~~z_3)\neq 0,~~  \det(x_{12}~~~ z_3)\neq 0,~~\det(x_{54}~~~z_4)\neq 0,~~\det(x_{76}~~~z_5)\neq 0,&\\
\det(x_{23}x_{12}~~~x_{43}x_{54}~~~x_{63}z_5)\neq 0,~~
\lz_1\lz_2\neq 0.&
\end{eqnarray*}
Define the polynomial map $$\varphi: Z\longrightarrow \co_{E_1\oplus
E_2}$$ by $\varphi(z)=x$. We will show that $\varphi$ is an
isomorphism. For $z\in Z,$ we take
\begin{eqnarray*}
&g_1=z_1=\frac{\lz_1}{\lz_2},~~g_2=(x_{12}z_1~~~z_3),~~g_3=(x_{23}x_{12}z_1~~~x_{43}x_{54}z_2~~~x_{63}z_5),\\
&g_4=(x_{54}z_2~~~z_4),~~g_5=z_2,~~g_6=(x_{76}~~~z_5),~~g_7=1.&
\end{eqnarray*}
Then it is straightforward to check that $x=g\bullet y$ for
$g=(g_i)_{i\in I}$. So $\varphi$ is well-defined.

The inverse map $\psi$ of $\varphi$
$$\psi:  \co_{E_1\oplus E_2}\longrightarrow Z,\psi(x)=z,$$ is defined as follows.
Let $x=g\bullet y$ be an element in $\co_{E_1\oplus E_2}$, where $g=(g_i)_{i\in I}\in
GL_{\az}$. Then the coefficients of $x$ satisfy the equations (\ref{f:5.2.6}), and this implies that they satisfy the following equations
\begin{eqnarray}\label{f:5.2.7}
&x_{23}x_{12}g_1=g_3\left(\begin{array}{l}
1\\0\\0\end{array}\right),~~x_{43}x_{54}g_5=g_3\left(\begin{array}{l} 0\\1\\0\end{array}\right),~~
x_{63}x_{76}g_7=g_3\left(\begin{array}{l}
 1\\1\\0\end{array}\right),&\\
&x_{23}g_2=g_3\left(\begin{array}{ll}1&0\\0&0\\0&1\end{array}\right),~~
x_{43}g_4=g_3\left(\begin{array}{lll}
0&0\\1&0\\0&1\end{array}\right),
~~x_{63}g_6=g_3\left(\begin{array}{ll}
1&0\\1&0\\0&1\end{array}\right).&\nonumber
\end{eqnarray}
Set
\begin{eqnarray*}&\frac{\lz_1}{\lz_2}=z_1=\frac{g_1}{g_7},~~ z_0=x_{12}z_1,~~ z_2=\frac{g_5}{g_7},
 ~~z_3=g_2\left(\begin{array}{l} 0\\1\end{array}\right),
~~z_4=g_4\left(\begin{array}{l} 0\\1\end{array}\right),
~~z_5=g_6\left(\begin{array}{l} 0\\1\end{array}\right).&\end{eqnarray*}
Then
\begin{eqnarray*}x_{23}x_{12}z_1+x_{43}x_{54}z_2=x_{63}x_{76},{\quad} x_{23}z_3=x_{43}z_4=x_{63}z_5.\end{eqnarray*}
By (\ref{f:5.2.7}), we know that $x_{23}x_{12}$ and $x_{43}x_{54}$ are
linearly independent over $k$. Thus $z_1$ and $z_2$ are uniquely
determined by $x$.
Now suppose $x_{23}z'_3=x_{43}z'_4=x_{63}z'_5$, and suppose that the
four matrices $(x_{12}~~~z'_3)$, $(x_{54}~~~z'_4)$,
$(x_{76}~~~z'_5)$ and $(x_{23}x_{12}~~~x_{43}x_{54}~~~x_{63}z'_5)$
are all invertible. We denote by C(M) the vector space generated by
columns of matrix M. By

$$\dim_k C(x_{23})\cap C(x_{43})=\dim_k C(x_{43})\cap C(x_{63})=\dim_k C(x_{63})\cap C(x_{23})=2+2-3=1,$$

we have
\begin{eqnarray}\label{f:5.2.8}x_{23}z_3'=ax_{23}z_3,x_{43}z_4'=ax_{43}z_4,x_{63}z_5'=ax_{63}z_5\end{eqnarray}
for some $a\in k^*.$
Since $x_{23}$, $x_{43}$ and $x_{63}$ have full column ranks, (\ref{f:5.2.8}) implies
$$z_3'=az_3,z_4'=az_4,z_5'=az_5.$$
Hence $\psi$ is well-defined, and $\varphi$ is an isomorphism.

Furthermore, let $X$ be the subvariety of $k^2\times
{\mathbb{P}}^1\times k\times P(k^2\times k^2\times k^2)\times
\bbe_{\az}$ consisting of those elements
$z=(z_0,[\lz_1:\lz_2],z_2,[z_3,z_4,z_5],x)$ such that
\[ x_{12}\lz_1=z_0\lz_2,~~ x_{23}z_0+x_{43}x_{54}z_2=x_{63}x_{76},~~ x_{23}z_3=x_{43}z_4=x_{63}z_5.\]
Then it is easy to
see that $X=\overline{\co}_{E_1\oplus E_2}.$

We denote by $P(i)$ (respectively, $I(i)$) the indecomposable projective (respectively injective)
module corresponding to $i$ for any $i\in I.$  Set
$Y=\overline{\co}_{\tau^{-2}P(1)\oplus E_2\oplus I(1)}.$ Since
$$ \Hom_{\llz}(\tau^{-2}P(1),I(1))=\Hom_{\llz}(E_2,I(1))=0,$$ and
$$\lr{\udim \tau^{-2}P(1),\udim E_2}=0=\dim_k\Hom( \tau^{-2}P(1),E_2),$$
we have
$$\dim_k \ed(\tau^{-2}P(1)\oplus E_2\oplus I(1))=3, ~\dim_kX=\dim_kY+1. $$
Furthermore, up to the isomorphism $\psi,$ we have
$$Y=\{z|z\in X, \lz_2=0\}.$$
%
%
Let $p:X\rightarrow Y$ be the canonical projection from $X$ to $Y$.
 Since
$\tau^{-1}P(1)\oplus E_2\oplus I(1)$ is aperiodic,
$IC_{Y}(\overline{\bbq}_l)$ has the purity property by Theorem~5.4
in \cite{L4} . It is clear that $p$ is  smooth with relative
dimension $1.$ Therefore, by Lemma \ref{l:5.1.4}, the statement of
the lemma is true.
\end{proof}

\subsection{Type $\tilde{E}_7$}
Let $Q$ be the quiver
\[\xymatrix{&&&8\ar[d]&&&\\ 1\ar[r]&2\ar[r]&3\ar[r]&4&5\ar[l]&6\ar[l]&7\ar[l]}\]
There are three non-homogeneous tubes $\ct_1,\ct_2,\ct_3$, respectively of periods $2,3,4$. Let $E_1,E_2$ be the regular simples in $\ct_1$, let $E'_1,E'_2,E'_3$ be the regualr simples in $\ct_2$, and let $E''_1,E''_2,E''_3.E''_4$ be the regular simples in $\ct_3$. Their dimension vectors are given as follows
\begin{itemize}
\item[$\ct_1$:] $(1,1,2,2,1,1,0,1)$, $(0,1,1,2,2,1,1,1)$;
\item[$\ct_2$:] $(1,1,1,2,1,1,1,1)$, $(0,1,1,1,1,1,0,0)$, $(0,0,1,1,1,0,0,1)$;
\item[$\ct_3$:] $(1,1,1,1,1,0,0,0)$, $(0,1,1,1,0,0,0,1)$, $(0,0,1,1,1,1,1,0)$, $(0,0,0,1,1,1,0,1)$.
\end{itemize}



\begin{lemma}\label{l:5.2.4}
Set $X_1=\overline{\co}_{E_1\oplus E_2},
X_2=\overline{\co}_{E'_1\oplus E'_2\oplus E'_3},$ and
$X_3=\overline{\co}_{E''_1\oplus E''_2\oplus E''_3\oplus E''_4}.$
 Then $IC_{X_i}(\overline{\bbq}_l)$ has the
purity property for any $i=1,2,3.$
\end{lemma}

\begin{proof} We only prove for $i=3.$ The other cases can be proved
similarly. Let $V$ be the $I-$graded vector space
$V=\oplus_{i=1}^{8}V_i,$ and $V_1=V_7=k,
V_2=V_6=V_8=k^2,V_3=V_5=k^3,V_4=k^4.$ Set
$\az=(1,2,3,4,3,2,1,2)\in\bbn [I],$ then
\begin{eqnarray*}
\bbe_{\az}=\bbe_{V}&=&\{x|~x=(x_{12},x_{23},x_{34},x_{54},x_{65},x_{76},x_{84}),\\
                    && x_{12},x_{76}\in M_{21}(k),
x_{23},x_{65}\in  M_{32}(k), x_{34},x_{54}\in  M_{43}(k),x_{84}\in
M_{42}(k)\}.
\end{eqnarray*}
Thanks to \cite{[DR]}, we can get $E''_1\oplus E''_2\oplus
E''_3\oplus E''_4=M(x),$ where $x\in\bbe_{\az}$ and

$$\begin{array}{l}
x_{12}=x_{76}=\left[\begin{array}{l}
1\\
0\end{array}\right],~x_{23}=\left[\begin{array}{ll}
1&0\\
0&1\\
0&0\end{array}\right],~ x_{65}=\left[\begin{array}{ll}
0&0\\
1&0\\
0&1\end{array}\right], ~x_{84}=\left[\begin{array}{ll}
0&0\\
1&0\\
0&0\\
0&1\end{array}\right],\\
 x_{34}=\left[\begin{array}{lll}
1&0&0\\
0&1&0\\
0&0&1\\
0&0&0\end{array}\right],~x_{54}=\left[\begin{array}{lll}
1&0&0\\
0&0&0\\
0&1&0\\
0&0&1\end{array}\right].
\end{array}$$
For any $x\in\co_{E''_1\oplus E''_2\oplus E''_3\oplus E_4''},$ there
 exists $(g_i)_{i\in I}\in GL_\az$ such that
\begin{eqnarray}\label{f:5.2.9}{\qquad}x_{12}=g_2\left[\begin{array}{l} 1\\0
\end{array}\right]g_1^{-1},
x_{23}=g_3\left[\begin{array}{ll}
1&0\\
0&1\\
0&0\end{array}\right]g_2^{-1}, x_{34}=g_4\left[\begin{array}{lll}
1&0&0\\0&1&0\\0&0&1\\0&0&0
\end{array}\right]g_3^{-1},\end{eqnarray}
$${\quad}x_{54}=g_4\left[\begin{array}{lll}
1&0&0\\0&0&0\\0&1&0\\0&0&1
\end{array}\right]g_5^{-1}, x_{65}=g_5\left[\begin{array}{ll}
0&0\\1&0\\0&1
\end{array}\right]g_6^{-1},
x_{76}=g_6\left[\begin{array}{l}
1\\
0
\end{array}\right]g_7^{-1},$$
${\qquad\qquad\qquad}x_{84}=g_4\left[\begin{array}{ll}
0&0\\1&0\\0&0\\0&1
\end{array}\right]g_8^{-1}.$

 Let $Z$ be the subvariety of $k^2\times \mathbb{P}^1\times k^2\times k^3 \times k^2\times
k^2\times k^3\times k^2\times k\times \bbe_\az$ consisting of those
elements $z=(z_0,[\lz_1:\lz_2],z_2,z_3,z_{41},z_{42},z_5,z_6,z_7,x)$
such that
$$\begin{array}{l}
{\qquad\qquad}x_{12}\lz_1=z_0\lz_2,x_{34}x_{23}z_0=x_{54}z_5,~x_{34}x_{23}z_2=x_{84}z_{41},
x_{34}z_3=x_{54}x_{65}x_{76}z_7,~x_{84}z_{42}=x_{54}x_{65}z_6,\\
{\qquad\qquad}z_7\neq 0, \left[\begin{array}{lll} z_2&:&z_5
\end{array}\right]\in P(k^2\times k^3),~\left[\begin{array}{lll}
z_2&:&z_{41}
\end{array}\right]\in P(k^2\times k^2),~\left[\begin{array}{lll}
z_3&:&z_7
\end{array}\right]\in P(k^3\times k),\\
{\qquad\qquad}\left[\begin{array}{lll} z_{42}&:&z_6
\end{array}\right]\in P(k^2\times k^2)
 ,\left[\begin{array}{ll} z_0&z_2
\end{array}\right],~\left[\begin{array}{ll} z_{41}&z_{42}
\end{array}\right],~\left[\begin{array}{ll} x_{12}&z_2
\end{array}\right],~\left[\begin{array}{ll} x_{76}&z_6
\end{array}\right]\in GL_2,\\
{\qquad\qquad}\left[\begin{array}{lll} x_{23}&z_3
\end{array}\right],
~\left[\begin{array}{ll} z_5&x_{65}
\end{array}\right]\in GL_3, \lz_1\lz_2\neq 0\}
\end{array}$$

We defined the polynomial map
$$\varphi:Z\longrightarrow\co_{E''_1\oplus E''_2\oplus E''_3\oplus
E''_4}$$
  by $\varphi(z)=x.$ We will show that $\varphi$ is an
isomorphism.

On the one hand, we take
$$\begin{array}{l}
g_1=z_1=\frac{\lz_1}{\lz_2},~ g_2=\left[\begin{array}{ll}
x_{12}z_1&z_2\end{array}\right],~ g_3=\left[\begin{array}{lll}
x_{23}x_{12}z_1&x_{23}z_2&z_3\end{array}\right],~g_5=\left[\begin{array}{lll}
z_5&x_{65}x_{76}z_7&x_{65}z_6\end{array}\right],
\\
 g_6=\left[\begin{array}{ll}
x_{76}z_7&z_6\end{array}\right],~g_7=z_7,g_8=\left[\begin{array}{ll}
z_{41}&z_{42}\end{array}\right],~g_4=\left[\begin{array}{llll}
x_{34}x_{23}x_{12}z_1&x_{34}x_{23}z_2&x_{34}z_3&x_{54}x_{65}z_6\end{array}\right].
\end{array}$$

 Thus,
$$\begin{array}{l}
g_2\left[\begin{array}{l} 1\\0
\end{array}\right]g_1^{-1}=\left[\begin{array}{ll}
x_{12}z_1&z_2\end{array}\right]\left[\begin{array}{l} 1\\0
\end{array}\right]z_1^{-1}=x_{12},\\
g_3\left[\begin{array}{ll}
1&0\\0&1\\0&0\end{array}\right]g_2^{-1}=\left[\begin{array}{lll}
x_{23}x_{12}z_1&x_{23}z_2&z_3\end{array}\right]\left[\begin{array}{ll}
1&0\\0&1\\0&0\end{array}\right] \left[\begin{array}{ll}
x_{12}z_1&z_2\end{array}\right]^{-1}=x_{23},\\
g_4\left[\begin{array}{lll}
1&0&0\\0&1&0\\0&0&1\\0&0&0\end{array}\right]g_3^{-1}=\left[\begin{array}{llll}
x_{34}x_{23}x_{12}z_1&x_{34}x_{23}z_2&x_{34}z_3&x_{54}x_{65}z_6\end{array}\right]\left[\begin{array}{lll}
1&0&0\\0&1&0\\0&0&1\\0&0&0\end{array}\right]\left[\begin{array}{lll}
x_{23}x_{12}z_1&x_{23}z_2&z_3\end{array}\right]^{-1}\\
{\qquad\qquad\qquad\qquad}=x_{34},\\
g_4\left[\begin{array}{lll}
1&0&0\\0&0&0\\0&1&0\\0&0&1\end{array}\right]g_5^{-1}=\left[\begin{array}{llll}
x_{34}x_{23}x_{12}z_1&x_{34}x_{23}z_2&x_{34}z_3&x_{54}x_{65}z_6\end{array}\right]\left[\begin{array}{lll}
1&0&0\\0&0&0\\0&1&0\\0&0&1\end{array}\right]\left[\begin{array}{lll}
z_5&x_{65}x_{76}z_7&x_{65}z_6\end{array}\right]^{-1}\\
{\qquad\qquad\qquad\qquad}=\left[\begin{array}{lll}
x_{34}x_{23}x_{12}z_1&x_{34}z_3&x_{54}x_{65}z_6\end{array}\right]\left[\begin{array}{lll}
z_5&x_{65}x_{76}z_7&x_{65}z_6\end{array}\right]^{-1}=x_{54},\\
g_5\left[\begin{array}{ll}
0&0\\1&0\\0&1\end{array}\right]g_6^{-1}=\left[\begin{array}{lll}
z_5&x_{65}x_{76}z_7&x_{65}z_6\end{array}\right]\left[\begin{array}{ll}
0&0\\1&0\\0&1\end{array}\right]\left[\begin{array}{ll}
x_{76}z_7&z_6\end{array}\right]^{-1}=x_{65},\\\end{array}$$

$$\begin{array}{l}
g_6\left[\begin{array}{l}
1\\0\end{array}\right]g_7^{-1}=\left[\begin{array}{ll}
x_{76}z_7&z_6\end{array}\right]\left[\begin{array}{l}
1\\0\end{array}\right]=x_{76},\\
g_4\left[\begin{array}{ll}
0&0\\1&0\\0&0\\0&1\end{array}\right]g_8^{-1}=\left[\begin{array}{llll}
x_{34}x_{23}x_{12}z_1&x_{34}x_{23}z_2&x_{34}z_3&x_{54}x_{65}z_6\end{array}\right]\left[\begin{array}{ll}
0&0\\1&0\\0&0\\0&1\end{array}\right]\left[\begin{array}{ll}
z_{41}&z_{42}\end{array}\right]^{-1}=x_{84},\\
\end{array}$$
that is, $\varphi$ is well defined.

On the other hand, we get the inverse map $\psi$ of $\varphi$
$$\psi:\co_{E''_1\oplus E''_2\oplus E''_3\oplus E''_4}\longrightarrow Z, ~\psi(x)=z,$$
and $z$ is defined in the following.

Since $x\in\co_{E''_1\oplus E''_2\oplus E''_3\oplus E''_4},$ there
exists $(g_i)_{i\in I}\in GL_\az$ such that $x$ satisfies (\ref
{f:5.2.9}) .Therefore,
\begin{eqnarray}\label{f:
5.2.10}{\qquad\qquad}x_{34}x_{23}x_{12}g_1=x_{54}g_5\left[\begin{array}{l}
1\\0\\0\end{array}\right]=g_4\left[\begin{array}{l}
1\\0\\0\\0\end{array}\right],~x_{34}x_{23}g_2\left[\begin{array}{l}
0\\1\end{array}\right]=x_{84}g_8\left[\begin{array}{l}
1\\0\end{array}\right]=g_4\left[\begin{array}{l}
0\\1\\0\\0\end{array}\right],\end{eqnarray}
$${\qquad\quad}x_{34}g_3\left[\begin{array}{l}
0\\0\\1\end{array}\right]=x_{54}x_{65}x_{76}g_7=g_4\left[\begin{array}{l}
0\\0\\1\\0\end{array}\right],~x_{84}g_8\left[\begin{array}{l}
0\\1\end{array}\right]=x_{54}x_{65}g_6\left[\begin{array}{l}
0\\1\end{array}\right]=g_4\left[\begin{array}{l}
0\\0\\0\\1\end{array}\right].
$$
Set
$$\begin{array}{l}\frac{\lz_1}{\lz_2}=z_1=g_1,~,z_0=x_{12}z_1,~z_2=g_2\left[\begin{array}{l}
0\\1\end{array}\right],~z_3=g_3\left[\begin{array}{l}
0\\0\\1\end{array}\right],~z_{41}=g_8\left[\begin{array}{l}
1\\0\end{array}\right],~z_{42}=g_8\left[\begin{array}{l}
0\\1\end{array}\right],\\
z_5=g_5\left[\begin{array}{l} 1\\0\\0\end{array}\right],
~z_6=g_6\left[\begin{array}{l}
0\\1\end{array}\right],z_7=g_7.\end{array}$$ It is easy to see that
$z=(z_0,[\lz_1:\lz_2],z_2,z_3,z_{41},z_{42},z_5,z_6,z_7,x)\in Z.$

 Suppose $z'=(z'_0,[\lz_1':\lz_2'],z'_2,z'_3,z'_{41},z'_{42},z'_5,z'_6,z'_7,x)\in
 Z,$ we have
 \begin{eqnarray}\label{f:5.2.11}{\qquad\qquad}~x_{34}x_{23}z_0=x_{54}z_5,~x_{34}x_{23}z_2=x_{84}z_{41},
x_{34}z_3=x_{54}x_{65}x_{76}z_7,~x_{84}z_{42}=x_{54}x_{65}z_6,z_7\neq
0,\end{eqnarray}
$${\qquad\quad}~x_{34}x_{23}z'_0=x_{54}z'_5,~x_{34}x_{23}z'_2=x_{84}z_{41},
x_{34}z'_3=x_{54}x_{65}x_{76}z'_7,~x_{84}z'_{42}=x_{54}x_{65}z'_6,z'_7\neq
0.
 $$

 Since $dim_k C(x_{54})\cap C(x_{34}x_{23})=1,$ and $x_{54}z'_5,x_{54}z_5\in C(x_{54})\cap C(x_{34}x_{23}),$ we have $x_{54}z_5=\lz
 x_{54}z'_5$ for some $\lz.$ Because $rank(x_{54})=3,$
 there exists a matrix $y_{54}$
such that $y_{54}x_{54}=I_3.$ Thus, $z_5=\lz z'_5.$ Moreover,
$z_0=\lz z'_0$, that is ,$~[z_0:z_5]\in P(k^2\times k^3).$

Based on (\ref{f:5.2.9}),  we have
$$C(x_{34}x_{23})=
C(g_4\left[\begin{array}{ll} 1&0\\0&1\\0&0\\0&0\end{array}\right])
\text{ and } C(x_{84})= C(g_4\left[\begin{array}{ll}
0&0\\1&0\\0&0\\0&1\end{array}\right]).$$
 Since
 $$\begin{array}{l}C(x_{84})\cap
C(x_{34}x_{23})=C(g_4\left[\begin{array}{l}
0\\1\\0\\0\end{array}\right]),\end{array}$$ we obtain $z_2=\mu z'_2$
and $z_{41}=\mu z'_{41}$ for some $\mu\in k^*,$ that is
$[z_2:z_{41}]\in P(k^2\times k^2).$ In the same way, we get
$[z_3:z_7]\in P(k^3\times k), [z_{42}:z_6]\in P(k^2\times k^2).$

Hence, $\psi$ is well defined, and $\varphi$ is an isomorphism.

Let
$$\begin{array}{l}X=\{z=(z_0,[\lz_1:\lz_2],z_2,z_3,z_{41},z_{42},z_5,z_6,z_7,x)|z\in k^2\times \mathbb{P}^1\times k^2\times k^3 \times k^2\times k^2\times
k^3\times k^2\times k\times \bbe_\az,\\
{\qquad\qquad}~x_{12}\lz_1=z_0\lz_2,~x_{34}x_{23}z_0=x_{54}z_5,~x_{34}x_{23}z_2=x_{84}z_{41},
x_{34}z_3=x_{54}x_{65}x_{76}z_7,~x_{84}z_{42}=x_{54}x_{65}z_6,\\
{\qquad\qquad}\left[\begin{array}{lll} z_0&:&z_5
\end{array}\right]\in P(k^2\times k^3),~\left[\begin{array}{lll}
z_2&:&z_{41}
\end{array}\right]\in P(k^2\times k^2),~\left[\begin{array}{lll}
z_3&:&z_7
\end{array}\right]\in P(k^3\times k),\\
{\qquad\qquad}\left[\begin{array}{lll} z_{42}&:&z_6
\end{array}\right]\in P(k^2\times k^2)
 \}
\end{array}$$
Then $X=\overline{\co}_{E'_1\oplus E'_2\oplus E'_3\oplus E'_4}.$

Set $Y=\overline{\co}_{\tau^{-2}P(7)\oplus E''_2\oplus E''_3\oplus
E''_4\oplus I(1)}.$ Since
$$Hom_{\llz}(\tau^{-2}P(7),I(1))=Hom_{\llz}(E_i'',I(1))=0,i=2,3,4$$
$$\lr{\udim \tau^{-2}P(7), \udim E_i''}=dim_k Hom_{\llz}(\tau^{-2}P(7), E_i'')=0, i=2,3,4,$$
we have
$$dim_k End_{\llz}(\tau^{-2}P(7)\oplus E''_2\oplus E''_3\oplus
E''_4\oplus I(1))=5,~dim X=dim Y+1.$$ Furthermore, up to the
isomorphism $\varphi,$ we have
$$Y=\{z|z \in X , \lz_2=0\}$$
Let $p:X\rightarrow Y$ be the canonical projection from $X$ to $Y$.
Thus,  Lemma \ref{l:5.2.4} follows from the lemma \ref{l:5.1.4}.
 \end{proof}

\subsection{Type $\widetilde{E}_8$}
Let $Q$ be the quiver
\[ \xymatrix{&&9\ar[d]&&&&&\\ 1\ar[r]&2\ar[r]&3&4\ar[l]& 5\ar[l]& 6\ar[l]& 7\ar[l]&8\ar[l]}\]
There are three non-homogeneous tubes $\ct_1,\ct_2,\ct_3$, respectively of periods $2,3,5$. Let $E_1,E_2$ be the regular simples in $\ct_1$, let $E'_1,E'_2,E'_3$ be the regualr simples in $\ct_2$, and let $E''_1,E''_2,E''_3.E''_4,E''_5$ be the regular simples in $\ct_3$. Their dimension vectors are given as follows
\begin{itemize}
\item[$\ct_1$:] $(1,2,3,2,2,1,1,0,2)$, $(1,2,3,3,2,2,1,1,1)$;
\item[$\ct_2$:] $(1,2,2,1,1,1,0,0,1)$, $(0,1,2,2,2,1,1,1,1)$, $(1,1,2,2,1,1,1,0,1)$;
\item[$\ct_3$:] $(1,1,1,1,1,0,0,0,0)$, $(0,1,1,1,0,0,0,0,1)$, $(1,1,2,1,1,1,1,1,1)$,$(0,1,1,1,1,1,1,0,0)$,\\ $(0,0,1,1,1,1,0,0,1)$.
\end{itemize}



\begin{lemma}\label{l:5.2.5}
Set $X_1=\overline{\co}_{E_1\oplus E_2},
X_2=\overline{\co}_{E'_1\oplus E'_2\oplus E'_3},$ and
$X_3=\overline{\co}_{E''_1\oplus E''_2\oplus E''_3\oplus E''_4\oplus
E''_5}.$
Then $IC_{X_i}(\overline{\bbq}_l)$ has the purity property for any $i=1,2,3.$
\end{lemma}

\begin{proof} We only prove for $i=1.$ The other cases can be proved
similarly. Let $V$ be the $I-$graded vector space
$V=\oplus_{i=1}^{9}V_i,$ and $V_1=V_7=k^2,
V_2=V_5=k^4,V_3=k^6,V_4=k^5,V_6=V_9=k^3,V_8=k.$ Set
$\az=(2,4,6,5,4,3,2,1,3)\in\bbn [I],$ then
\begin{eqnarray*}
\bbe_{\az}=\bbe_{V}&=&\{x|~x=(x_{12},x_{23},x_{43},x_{54},x_{65},x_{76},x_{87},x_{93}),\\
                    && x_{12}\in M_{42}(k),x_{23}\in
M_{64}(k),x_{43}\in M_{65}(k) ,x_{54}\in M_{54}(k) ,x_{65}\in
M_{43}(k) ,x_{76}\in M_{32}(k),\\
&&x_{87}\in M_{21}(k),x_{93}\in M_{63}(k)\}.
\end{eqnarray*}
Thanks to \cite{[DR]}, we can get $E_1\oplus E_2=M(x),$ where
$x\in\bbe_{\az}$ and

$$\begin{array}{l}
x_{12}=\left[\begin{array}{ll} 0&0\\1&0\\0&0\\0&1
\end{array}\right],~x_{23}=\left[\begin{array}{llll}
0&0&0&0\\
1&0&0&0\\
0&1&0&0\\
0&0&0&0\\
0&0&1&0\\
0&0&0&1\end{array}\right],~ x_{43}=\left[\begin{array}{lllll}
1&0&0&0&0\\
0&1&0&0&0\\
0&0&0&0&0\\
0&0&1&0&0\\
0&0&0&1&0\\
0&0&0&0&1\end{array}\right], ~x_{54}=\left[\begin{array}{llll}
1&0&0&0\\
0&1&0&0\\
0&0&1&0\\
0&0&0&1\\
0&0&0&0\end{array}\right],\\
 x_{65}=\left[\begin{array}{lll}
1&0&0\\
0&0&0\\
0&1&0\\
0&0&1\end{array}\right],~x_{76}=\left[\begin{array}{ll}
1&0\\
0&1\\
0&0\end{array}\right],~x_{87}=\left[\begin{array}{l}
0\\
1\end{array}\right],~x_{93}=\left[\begin{array}{lll}
1&0&0\\
1&1&0\\
0&1&0\\
0&0&1\\
0&0&1\\
0&0&1.\end{array}\right].
\end{array}$$

For any $x\in\co_{E_1\oplus E_2},$ there exists $(g_i)_{i\in I}\in
GL_\az$ such that
\begin{eqnarray}\label{f:5.2.11}{\qquad\qquad}x_{12}=g_2\left[\begin{array}{ll}
0&0\\1&0\\0&0\\0&1
\end{array}\right]g_1^{-1},
x_{23}=g_3\left[\begin{array}{llll}
0&0&0&0\\
1&0&0&0\\
0&1&0&0\\
0&0&0&0\\
0&0&1&0\\
0&0&0&1
\end{array}\right]g_2^{-1},
x_{43}=g_3\left[\begin{array}{lllll}
1&0&0&0&0\\
0&1&0&0&0\\
0&0&0&0&0\\
0&0&1&0&0\\
0&0&0&1&0\\
0&0&0&0&1\end{array}\right]g_4^{-1},\end{eqnarray}
$$\begin{array}{l} x_{54}=g_4\left[\begin{array}{llll}
1&0&0&0\\
0&1&0&0\\
0&0&1&0\\
0&0&0&1\\
0&0&0&0\end{array}\right]g_5^{-1},
x_{65}=g_5\left[\begin{array}{lll}
1&0&0\\
0&0&0\\
0&1&0\\
0&0&1\end{array}\right]g_6^{-1}, x_{76}=g_6\left[\begin{array}{ll}
1&0\\
0&1\\
0&0\end{array}\right]g_7^{-1},\\
 x_{87}=g_7\left[\begin{array}{l} 0\\1
\end{array}\right]g_8^{-1},~x_{93}=g_3\left[\begin{array}{lll}
1&0&0\\
1&1&0\\
0&1&0\\
0&0&1\\
0&0&1\\
0&0&1\end{array}\right]g_9^{-1}.\end{array}$$

Set $$\begin{array}{l}
 Z=\{z=(z_{11},z_{12},z_2,z_3,z_4,z_5,z_6,z_7,z_8,[\lz_1:\lz_2],z_{91},z_{92},z_{93},x)|\\
 {\qquad\qquad}z\in k^2\times k^2\times k^4\times k^4
 \times k^5\times k^4\times k^3\times k^2\times k^2\times\mathbb{P}^1\times k^3\times k^3\times k^3\times
 \bbe_\az,~x_{87}\lz_1=z_8\lz_2,\\{\qquad\qquad}x_{23}z_2=x_{43}x_{54}z_5,
~x_{23}z_3=x_{43}x_{54}x_{65}z_6,~
~x_{93}z_{91}=x_{43}x_{54}x_{65}x_{76}z_7+x_{23}z_2,\\{\qquad\qquad}x_{93}z_{92}=x_{23}x_{12}z_{11}+x_{43}x_{54}z_5,
~x_{23}x_{12}z_{12}=x_{43}z_4,
\\{\qquad\qquad}
x_{93}z_{93}=x_{43}x_{54}x_{65}x_{76}z_8+x_{23}z_3+x_{23}x_{12}z_{12},\\{\qquad\qquad}[z_3:z_6]\in
P(k^4\times k^3),
 [z_{12}:z_4]\in P(k^2\times k^5),\left[\begin{array}{ll}
 z_{11}&z_{12}\end{array}\right],~\left[\begin{array}{ll}
 z_7&z_8\end{array}\right]\in GL_2,\\{\qquad\qquad}\left[\begin{array}{lll}
 x_{76}z_7&x_{76}z_8&z_6\end{array}\right]\in GL_3
 ~\left[\begin{array}{llll}
 z_2&x_{12}z_{11}&z_3&x_{12}z_{12}\end{array}\right]\in GL_4,\\{\qquad\qquad}\left[\begin{array}{llll}
 x_{65}x_{76}z_7&z_5&x_{65}x_{76}z_8&x_{65}z_6\end{array}\right]\in GL_4,\\
{\qquad\qquad} \left[\begin{array}{llllll}
 x_{43}x_{54}x_{65}x_{76}z_7&x_{43}x_{54}z_5&x_{23}x_{12}z_{11}&x_{43}x_{54}x_{65}x_{76}z_8
 &x_{43}x_{54}x_{65}z_6&x_{43}z_4\end{array}\right]\in GL_6,\\
 {\qquad\qquad}\left[\begin{array}{lll}
 z_{91}&z_{92}&z_{93}\end{array}\right]\in GL_3,~\left[\begin{array}{lllll}
 x_{54}x_{65}x_{76}z_7&x_{54}z_5&x_{54}x_{65}x_{76}z_8&x_{54}x_{65}z_6&z_4\end{array}\right]\in GL_5,~\lz_1\lz_2\neq 0\}
\end{array}$$

We define the polynomial map
$$\varphi: Z\longrightarrow \co_{E_1\oplus E_2}$$
 by $\varphi(z)=x$.

 If we take
$$\begin{array}{l}
g_1=\left[\begin{array}{ll}
 z_{11}&z_{12}\end{array}\right],~g_2=\left[\begin{array}{llll}
 z_2&x_{12}z_{11}&z_3&x_{12}z_{12}\end{array}\right],\\g_3=\left[\begin{array}{llllll}
 x_{43}x_{54}x_{65}x_{76}z_7&x_{43}x_{54}z_5&x_{23}x_{12}z_{11}&x_{43}x_{54}x_{65}x_{76}z_8
 &x_{43}x_{54}x_{65}z_6&x_{43}z_4\end{array}\right],\\
 g_4=\left[\begin{array}{lllll}
 x_{54}x_{65}x_{76}z_7&x_{54}z_5&x_{54}x_{65}x_{76}z_8&x_{54}x_{65}z_6&z_4\end{array}\right],
 ~g_5=\left[\begin{array}{llll}
 x_{65}x_{76}z_7&z_5&x_{65}x_{76}z_8&x_{65}z_6\end{array}\right],\\
 g_6=\left[\begin{array}{lll}
 x_{76}z_7&x_{76}z_8&z_6\end{array}\right],~g_7=\left[\begin{array}{ll}
 z_7&z_8\end{array}\right],
\end{array}$$
one may get $x\in\bbe_\az$ satisfies (\ref{f:5.2.11}), that is, $\varphi$
is well defined.

On the other hand, we have the inverse map $\psi$ of $\varphi$
$$\psi: \co_{E_1\oplus E_2}\longrightarrow Z, \psi(x)=z,$$
and $z$ is defined as in the following.

Since $x\in\co_{E_1\oplus E_2},$ there exists $(g_i)_i\in GL_{\udim
E_1\oplus E_2}$ such that $x$ satisfies (\ref{f:5.2.11}). Moreover,
\begin{eqnarray}
\label{f:5.2.12}x_{23}g_2\left[\begin{array}{l}
 1\\0\\0\\0\end{array}\right]=x_{43}x_{54}g_5\left[\begin{array}{l}
 0\\1\\0\\0\end{array}\right],~x_{23}g_2\left[\begin{array}{l}
 0\\0\\1\\0\end{array}\right]=x_{43}x_{54}x_{65}g_6\left[\begin{array}{l}
 0\\0\\1\end{array}\right],\end{eqnarray}
$$\begin{array}{l}
 {\qquad\qquad}x_{23}x_{12}g_1\left[\begin{array}{l}
 0\\1\end{array}\right]=x_{43}g_4\left[\begin{array}{l}
 0\\0\\0\\0\\1\end{array}\right],~x_{93}g_9\left[\begin{array}{l}
 1\\0\\0\end{array}\right]=x_{43}x_{54}x_{65}x_{76}g_7\left[\begin{array}{l}
 1\\0\end{array}\right]+x_{23}g_2\left[\begin{array}{l}
 1\\0\\0\\0\end{array}\right],\\
 {\qquad\qquad}x_{93}g_9\left[\begin{array}{l}
 0\\1\\0\end{array}\right]=x_{23}x_{12}g_1\left[\begin{array}{l}
 1\\0\end{array}\right]+x_{43}x_{54}g_5\left[\begin{array}{l}
 0\\1\\0\\0\end{array}\right],\\
 {\qquad\qquad}x_{93}g_9\left[\begin{array}{l}
 0\\0\\1\end{array}\right]
 =x_{43}x_{54}x_{65}x_{76}x_{87}g_8+x_{23}g_2\left[\begin{array}{l}
 0\\0\\1\\0\end{array}\right]+x_{23}g_2\left[\begin{array}{l}
 0\\0\\1\\0\end{array}\right]+x_{23}x_{12}g_1\left[\begin{array}{l}
 0\\1\end{array}\right].
\end{array}$$
Set
$$\begin{array}{l}
z_{11}=g_1\left[\begin{array}{l}
 1\\0\end{array}\right],~z_{12}=g_1\left[\begin{array}{l}
 0\\1\end{array}\right],~z_2=g_2\left[\begin{array}{l}
 1\\0\\0\\0\end{array}\right],~z_3=g_2\left[\begin{array}{l}
 0\\0\\1\\0\end{array}\right],~z_4=g_4\left[\begin{array}{l}
 0\\0\\0\\0\\1\end{array}\right],~z_5=g_5\left[\begin{array}{l}
 0\\1\\0\\0\end{array}\right],\\
 z_6=g_6\left[\begin{array}{l}
 0\\0\\1\end{array}\right],~z_7=g_7\left[\begin{array}{l}
 1\\0\end{array}\right],~z_8=x_{87}g_8,~ \frac{\lz_1}{\lz_2}=g_8,
\end{array}$$
then we have
$z=(z_{11},z_{12},z_2,z_3,z_4,z_5,z_6,z_7,z_8,[\lz_1:\lz_2],z_{91},z_{92},z_{93},x)\in
Z.$

Suppose
$z'=(z'_{11},z'_{12},z'_2,z'_3,z'_4,z'_5,z'_6,z'_7,z'_8,[\lz'_1:\lz'_2],z'_{91},z'_{92},z'_{93},x)\in
Z,$ then we also have
\begin{eqnarray}
\label{f:5.2.13}x_{23}z'_2=x_{43}x_{54}z'_5, ~x_{23}z'_3=x_{43}x_{54}x_{65}z'_6,~x_{87}\lz_1'=z_8'\lz_2',\end{eqnarray}
$$\begin{array}{l}
{\qquad}x_{93}z'_{91}=x_{43}x_{54}x_{65}x_{76}z'_7+x_{23}z'_2,~x_{93}z'_{92}
=x_{23}x_{12}z'_{11}+x_{43}x_{54}z'_5,\\
{\qquad}x_{23}x_{12}z'_{12}=x_{43}z'_4,
~x_{93}z'_{93}=x_{43}x_{54}x_{65}x_{76}z'_8+x_{23}z'_3+x_{23}x_{12}z'_{12}.
\end{array}$$
Based on (\ref{f:5.2.11}), we obtain
$$C(x_{23})=C(g_3\left[\begin{array}{llll}
0&0&0&0\\
1&0&0&0\\
0&1&0&0\\
0&0&0&0\\
0&0&1&0\\
0&0&0&1
\end{array}\right]),\text{ and } C(x_{43}x_{54}x_{65})=C(g_3\left[\begin{array}{lll}
1&0&0\\
0&0&0\\
0&0&0\\
0&1&0\\
0&0&1\\
0&0&0
\end{array}\right]),$$
$$C(x_{23}x_{12})=C(g_3\left[\begin{array}{ll}
0&0\\
0&0\\
1&0\\
0&0\\
0&0\\
0&1
\end{array}\right]) \text{ and } C(x_{43})=C(g_3\left[\begin{array}{lllll}
1&0&0&0&0\\
0&1&0&0&0\\
0&0&0&0&0\\
0&0&1&0&0\\
0&0&0&1&0\\
0&0&0&0&1
\end{array}\right])$$
Therefore,
$$C(x_{43}x_{54}x_{65})\cap
C(x_{23})=C(g_3\left[\begin{array}{l}
0\\
0\\
0\\
0\\
1\\
0
\end{array}\right] ),{\quad} C(x_{23}x_{12})\cap C(x_{43})=C(g_3\left[\begin{array}{l}
0\\
0\\
0\\
0\\
0\\
1
\end{array}\right]).$$ It implies that $z'_3=\mu z_3, z'_6=\mu z_6,z'_{12}=\nu z_{12}, \text{ and } z'_4=\nu z_4$
for some $\mu$ and $\nu,$ that is, $[z_3:z_6]\in P(k^4\times k^3)
\text{ and } [z_{12}:z_4]\in P(k^2\times k^5).$

Hence $\psi$ is well defined, and $\varphi$ is an isomorphism.

Let
$$\begin{array}{l}X=\{z=(z_{11},z_{12},z_2,z_3,z_4,z_5,z_6,z_7,z_8,[\lz_1:\lz_2],z_{91},z_{92},z_{93},x)|\\
 {\qquad\qquad}z\in k^2\times k^2\times k^4\times k^4
 \times k^5\times k^4\times k^3\times k^2\times k^2\times\mathbb{P}^1\times k^3\times k^3\times k^3\times
 \bbe_\az,~x_{23}z_2=x_{43}x_{54}z_5,\\
{\qquad\qquad}x_{23}z_3=x_{43}x_{54}x_{65}z_6,~
~x_{93}z_{91}=x_{43}x_{54}x_{65}x_{76}z_7+x_{23}z_2,~x_{93}z_{92}=x_{23}x_{12}z_{11}+x_{43}x_{54}z_5,\\
{\qquad\qquad}x_{23}x_{12}z_{12}=x_{43}z_4,
~x_{93}z_{93}=x_{43}x_{54}x_{65}x_{76}x_{87}z_8+x_{23}z_3+x_{23}x_{12}z_{12},~x_{87}\lz_1=z_8\lz_2,\\
 {\qquad\qquad}[z_3:z_6]\in P(k^4\times k^3),~[z_{12}:z_4]\in P(k^2\times k^5)\}\end{array},$$
then we get $X=\overline{\co}_{E_1\oplus E_2}.$

Let $P$ be the pre-projective with dimension vector $(123322101),$ and
set $Y=\overline{\co}_{P\oplus E_1\oplus I(8)}.$

Since
 $$\begin{array}{l}Hom_{\llz}(P,I(8))=Hom_{\llz}(E_1,I(8))=0,\\
\lr{\udim P,\udim E_1}=Hom_{\llz}(P,E_1)=0,
 \end{array}$$
we get
$$dim End_{\llz}(P\oplus E_1\oplus I(8))=3, dim X=dim Y+1.$$

Furthermore, up to the isomorphism $\varphi,$ we have
$$Y=\{z|z\in X, \lz_2=0\}.$$
Let $x'\in\bbe_{\az}$ such that $x=(0,x_{87},0)+x'$ and define
$$p:X\longrightarrow Y, p(z)=(z_{11},z_{12},z_2,z_3,z_4,z_5,z_6,z_7,z_8,[\lz_1:0],z_{91},z_{92},z_{93},x').$$
Applying Lemma \ref{l:5.1.4}, the proof is complete.

\end{proof}

Furthermore, in the same way, we can prove that the closure of
orbits of semi-simple objects has the  purity property.

In order to prove Theorem \ref{t:5.1.5}, we not only need to discuss the
closure of semi-simple objects in $\ct_i$ which has purity property,
 but also need to study the fibres of $p_3$.

\section{The Fibres of $p_3$}

 \subsection{} Let $\cp_{\gamma}$ be the set of
 $\llz$ modules of dimension vector $\gamma$, up to isomorphism.
 From the definition of $p_3$ in $1.3,$ it follows that ${p_3}^{-1}(L)=\cup_{N\in\cp_\bz,N\subseteq L}{\bf
 Z}_{L,M,N}.$   Thus we need to discuss some properties of the variety ${\bf
  Z}_{L,M,N}.$ In fact, we know that $g^L_{MN}$ is the number of
  rational points of the variety ${\bf
 Z}_{L,M,N}$. Set $Z_L=\{(M,N)|(M,N)\in\cp_\az\times\cp_\bz,N\subseteq L, \text{ and }L/N\cong M\}.$
 If $Z_L$ has only finitely many
 elements, then $p^{-1}(L)$ has a natural stratification with strata ${\bf
 Z}_{L,M,N}$ indexed by $(M,N)\in Z_L.$

For the convenience of discussions below, we need to introduce some
 notations about the BGP reflection functors (see[BGP] or[DR]).

We define $\sz_iQ$ to be the quiver obtained from $Q$ by reversing
the direction of every arrow connected to the vertex $i.$ If $i$ is
a sink of $Q,$  $\sz_i^+$ is defined as follows:
$$\sz_i^+:\mod\llz\lra\mod\sz_i\llz,$$
where $\llz=\fq(Q)$ (resp. $\llz=k(Q)$) and $\sz_i\llz=\fq(\sz_iQ)$
(resp. $\sz_i\llz=k(\sz_iQ)$) are path algebras. Therefore $\sz_i^+$
is an exact functor on the full subcategory $\mod\llz(i)$ of
$\mod\llz$ consisting of modules which does not have $S_i$ as a
direct summand, and induce quasi-inverse equivalence between
$\mod\llz(i)$ and the full subcategory $\mod\sz_i\llz(i)$ consisting
of modules which does not have direct summand
 isomorphic to the simple injective module $S_i.$

 According to the notation of Hall polynomial (see \cite{R3}), we use $|\cdot|$
to denote the ordinary cardinality of a finite set.

\begin{lemma}\label{l:6.1} {\sl Let $M,M_1,M_2$ and $N$ be  $\llz-$
 modules and let $N$ (resp. $M_2$) be a pre-projective ( resp. regular). If $Hom_{\llz}(M_2,M_1)=0,$ then

   $$g^{M\oplus M_2}_{M_1\oplus
   M_2,N}=g^{M}_{M_1,N}\frac{|Hom_{\llz}(M,M_2)|}{|Hom_{\llz}(M_1,M_2)|}.$$}
\end{lemma}

\begin{proof} Let $X,Y$ and $L$ be $\llz-$modules, and set
$$\begin{array}{l}
W(X,Y;L)=\{(f,g)\in Hom_{\llz}(X,L)\times Hom_{\llz}(L,Y)|\\
0\longrightarrow X\longrightarrow L\longrightarrow Y\longrightarrow
0\text{ is a short exact sequence }\}.
\end{array}$$
The action of $Aut(X)\times Aut(Y)$ on $W(X,Y;L)$ is defined by
$$(a,c)(f,g)=(fa,c^{-1}g).$$
It induces the orbit space
$$V(X,Y;L)=\{(f,g)^{\wedge}|(f,g)\in W(X,Y;L)\}.$$
 Then $|V(X,Y;L)|=g_{YX}^L.$

Note that the actions of $Aut(N)\times Aut(M_1\oplus M_2)$ and
$Aut(N)\times Aut( M_1)$ on $W(N,M_1\oplus M_2;M\oplus M_2)$  and
$W(N,M_1;M),$ respectively, are free. So we have \begin{eqnarray}
\label{f:6.1}
{\qquad\qquad\qquad}|V(N, M_1;M)|=\frac{|W(N,
M_1;M)|}{|Aut(N)||Aut( M_1)|}
\text{ ,~and }\end{eqnarray}
$$\begin{array}{l}
{\qquad\qquad\qquad\qquad}|V(N,M_1\oplus M_2;M\oplus
M_2)|=\frac{|W(N,M_1\oplus M_2;M\oplus M_2)|}{|Aut(N)||Aut(M_1\oplus
M_2)|}.\end{array}$$

If $(\left[\begin{array}{l}
f_1\\
f_2
\end{array}\right],\left[\begin{array}{ll}
g_{11}&0\\
g_{21}&g_{22}
\end{array}\right])\in W(N,M_1\oplus M_2;M\oplus M_2).$ By the regular part of $\mod\llz$ is an abelian subcategory, then we get $g_{22}$ is injective. So is invertible and $f_2=-g_{22}^{-1}g_{21}f_1.$

Consider the map
$$\begin{array}{l}\varphi:W(N,M_1\oplus M_2;M\oplus M_2)\longrightarrow
W(N, M_1;M)\times Hom_{\llz}(M,M_2)\times Aut(M_2),\\
(\left[\begin{array}{l}
f_1\\
-g_{22}^{-1}g_{21}f_1
\end{array}\right],\left[\begin{array}{ll}
g_{11}&0\\
g_{21}&g_{22}
\end{array}\right])\longmapsto (f_1,g_{11},g_{21},g_{22}).
\end{array}$$
It is easy to see that the inverse map $\psi$ of $\varphi$ is
$$\begin{array}{l}\psi:W(N, M_1;M)\times Hom_{\llz}(M,M_2)\times Aut(M_2)
\longrightarrow W(N,M_1\oplus M_2;M\oplus M_2),\\
(f_1,g_{11},g_{21},g_{22})\longmapsto (\left[\begin{array}{l}
f_1\\
-g_{22}^{-1}g_{21}f_1
\end{array}\right],\left[\begin{array}{ll}
g_{11}&0\\
g_{21}&g_{22}
\end{array}\right]).
\end{array}$$

 Thus, $\varphi$ is an isomorphism. Therefore, we have
\begin{eqnarray}\label{f:6.2}|W(N,M_1\oplus M_2;M\oplus M_2)|=|W(N, M_1;M)||Hom_{\llz}(M,M_2)\times Aut(M_2)|.\end{eqnarray}
 Hence the proof follows from (\ref{f:6.1}) and (\ref{f:6.2}). \end{proof}

\begin{proposition}\label{p:6.2} {\sl Let $S$ be a simple projective
module of $k Q$ ( except $Q=\widetilde{E}_8$) corresponding to a
unique sink point, $P$ be a preprojective, and let $M$ be a regular
semi-simple object in $\ct_i$ for some $i,1\leqslant i\leqslant l.$
Then ${p_3^{-1}|_{\overline{\co}_M\times\overline{\co}_S}}(P\oplus L)$ has the purity property.~}\end{proposition}

\begin{proof} If $g_{MS}^{P\oplus L}\neq 0,$ we have
$$0\longrightarrow S\stackrel{f}\longrightarrow P\oplus L\stackrel{\left[\begin{array}{l}
g_1\\g_2
\end{array}\right]^t}\longrightarrow M\longrightarrow 0.$$
According to the representation theory of quivers, we know that
$L\in\ct_i.$ First, we claim that $L$ is a regular semi-simple
submodule of $M.$ Obviously, $g_2$ is a monomorphism. Let,
otherwise, $rad_{\ct_i}(L)$ be the regular radical  of $L$ in the
full sub-category $\ct_i$ of $\mod\llz,$ then we have
$g_2(rad_{\ct_i}(L))=0$ and $rad_{\ct_i}(L)\subseteq Im (f)$. But
$Hom_{\llz}(rad_{\ct_i}(L),S)=0,$ this is a contradiction.

Because $M$ and $L$ are semi-simple in $\ct_i,$ by $g_2$ is a monomorphism, there exists a morphism $l_2$
from $M$ to $L$ such that $l_2g_2=id_L.$ Then
$M=M_1\oplus Im(g_2),$ where $M_1=\{x-g_2l_2(x)|x\in M\}.$

Since
$M\cong M_1\oplus Im (g_2),$ the above short exact sequence may be rewritten as following form:
$$\begin{array}{l}0\longrightarrow S\stackrel{\left[\begin{array}{l}
f_1\\-g_{22}^{-1}g_{21}f_1
\end{array}\right]}\longrightarrow P\oplus L\stackrel{\left[\begin{array}{ll}
g_{11}&0\\g_{21}&g_{22}
\end{array}\right]}\longrightarrow M_1\oplus Im (g_2)\longrightarrow
0\end{array}$$
where  $g_{22}$ is an isomorphism. Moreover,
$$\begin{array}{l} 0\longrightarrow S\stackrel{
f_1 }\longrightarrow P\stackrel{ g_{11}}\longrightarrow
M_1\longrightarrow 0.
\end{array}$$

 For any
  $$(\left[\begin{array}{l}
f_1\\
f_2
\end{array}\right],\left[\begin{array}{ll}
g_{11}&g_{12}\\
g_{21}&g_{22}
\end{array}\right])\in W(S, M_1\oplus L;P\oplus L).$$
Because  $\left[\begin{array}{l}
g_{12}\\
g_{22}
\end{array}\right]$ is a split monomorphism, there is a morphism
$\left[\begin{array}{ll} l_{21}& l_{22}
\end{array}\right]$
from $M_1\oplus L$ to $L$ such that $$\left[\begin{array}{ll}
l_{21}& l_{22}
\end{array}\right]\left[\begin{array}{l}
g_{12}\\
g_{22}
\end{array}\right]=id_L.$$

Without loss of generality, we may assume that $l_{22}$ is
invertible up to an automorphism in $Aut_{\mod\llz}(M_1\oplus L)$.

Since
$$\left[\begin{array}{ll}
id_{M_1}&-g_{12}\\
0&id_L
\end{array}\right]\left[\begin{array}{ll}
id_{M_1}& 0\\l_{21}&l_{22}
\end{array}\right]\left[\begin{array}{ll}
g_{11}&g_{12}\\
g_{21}&g_{22}
\end{array}\right]=\left[\begin{array}{ll}
g_{11}-g_{12}(l_{21}g_{11}+l_{22}g_{21})&0\\
l_{21}g_{11}+l_{22}g_{21}&id_L
\end{array}\right],$$
we have $c=\left[\begin{array}{ll}
id_{M_1}&-g_{12}\\
0&id_L
\end{array}\right]\left[\begin{array}{ll}
id_{M_1}& 0\\l_{21}&l_{22}
\end{array}\right]\in Aut_{\llz}(M_1\oplus L)$ such that
$$c\left[\begin{array}{ll}
g_{11}&g_{12}\\
g_{21}&g_{22}
\end{array}\right]=\left[\begin{array}{ll}
g'_{11}&0\\
g'_{21}&1
\end{array}\right].$$

Thus $$V(S,M_1\oplus L;P\oplus L)=\{(\left[\begin{array}{l}
f_1\\
-g_{21}f_1
\end{array}\right],\left[\begin{array}{ll}
g_{11}&0\\
g_{21}&1
\end{array}\right])^{\wedge}|(\left[\begin{array}{l}
f_1\\
-g_{21}f_1
\end{array}\right],\left[\begin{array}{ll}
g_{11}&0\\
g_{21}&1
\end{array}\right])\in W(S,M_1\oplus L;P\oplus L)\}.$$

Consider the map
$$\begin{array}{l}\varphi:V(S, M_1\oplus L;P\oplus L)\longrightarrow V(S, M_1;P)
\end{array}$$
sending $(\left[\begin{array}{l}
f_1\\
-g_{21}f_1
\end{array}\right],\left[\begin{array}{ll}
g_{11}&0\\
g_{21}&1
\end{array}\right])^{\wedge}$ to $(f_1,g_{11})^{\wedge}.$ Let $V(g_{11})=\{hg_{11}| h\in
Hom_{\llz}(M_1,L)\}.$ Then $V(g_{11})$ is a $k-$ vector subspace of
$Hom_{\llz}(P,L)$ . Because  $g_{11}$ is an epimorphism, we have
$V(g_{11})\cong Hom_{\llz}(M_1,L)$ as a $k-$spaces. Therefore
$$(A)\left\{{\qquad}\begin{array}{l} \varphi^{-1}((f_1,g_{11})^{\wedge})=Hom_{\llz}(P,L)/V(g_{11})=\bba^{dim_k
Hom_{\llz}(P,L)/V(g_{11})},\\
{\qquad} |V(S, M_1\oplus L;P\oplus L)|=|V(S, M_1;P)|\times
|\bba^{dim_k Hom_{\llz}(P,L)/V(g_{11})}|.\end{array}\right.$$

  The rest of the proof can be divided into several cases:

Case $1: Q=\widetilde{A}_n.$

Let $E_{i}, 1\leq i\leq n$ be simple objects in the full
non-homgeneous subcategory of $\llz-\text{mod}$ corresponding to the
dimension vectors listed in 6.1. Since there is only one
simple object $E_n$ such that $Ext^1(E_{n},S)\neq 0,$ by Lemma~6.1,
we only need to discuss the case $N=\oplus mE_{n}$.

  Because $P(n)\bigoplus\oplus
(m-1)E_{n}$ is a unique non-trivial extension of $S$ by $N,$ it is
easy to see that $g_{mE_{n},S}^{P(n)\bigoplus\oplus (m-1)E_{n}}=1.$

Set $M=\oplus a_nE_{n}\bigoplus \oplus_{i\neq n}a_iE_{i}.$ Since

$$Hom_{\llz}(\oplus_{i\neq n}a_iE_{i},\oplus
a_nE_{n})=Hom_{\llz}(\oplus a_nE_{n},\oplus_{i\neq n}a_iE_{i})=0,$$

we obtain, by Lemma \ref{l:6.1},
$$\begin{array}{l}
g_{\oplus a_nE_{n}\bigoplus \oplus_{i\neq n}a_iE_{i},S}^{(P(n)\oplus
(a_n-1)E_{n})\bigoplus(\oplus_{i\neq n}a_iE_{i})}=g_{\oplus
a_nE_{n},S}^{P(n)\oplus (a_n-1)E_{n}}|Hom_{\llz}(P(n)\oplus
(a_n-1)E_{n},\oplus_{i\neq
n}a_iE_{i})|\\
\qquad\qquad\qquad\qquad\qquad\qquad\quad=|Hom_{\llz}(P(n),\oplus_{i\neq
n}a_iE_{i})|
\end{array}.$$

 Moreover,
$$\begin{array}{l}{\bf Z}_{({P(n)\oplus
(a_n-1)E_{n})\bigoplus\oplus_{i\neq
n}a_iE_{i}},M,S}=Hom_{\llz}(P(n),\oplus_{i\neq
n}a_iE_{i}),\\
({p_3^{-1}|_{\overline{\co}_M\times\overline{\co}_S}})(P(n)\oplus
(a_n-1)E_{n})\bigoplus(\oplus_{i\neq
n}a_iE_{i})=Hom_{\llz}(P(n),\oplus_{i\neq n}a_iE_{i}).\end{array}$$
Therefore ${\bf Z}_{{(P(n)\oplus
(a_n-1)E_{n})\bigoplus(\oplus_{i\neq n}a_iE_{i})},M,S}$ has the
purity property.

Case $2: Q=\widetilde{D}_n,n\geq 4.$

 We only consider $\ct_1$ in 6.2. The other cases can be proved
 similarly.
Let $E_{i}, i=1,2$ be simple objects in the full subcategory $\ct_1$
of $\llz-\text{mod}$ corresponding to the dimension vectors listed
in 6.2.

Since
$$dim_k Ext^1(E_{i},S)=-\lr{\udim E_{i}, \udim S}=1,$$
we have
$$dim_k Ext^1(\oplus mE_{i},S)=m.$$
Moreover,
$$\begin{array}{l}\xi:{\quad}0\longrightarrow S\stackrel{f_{\xi}=\left[\begin{array}{l}
1\\-1
\end{array}\right]}\longrightarrow P(1)\oplus P(n)\stackrel{g_{\xi}=\left[\begin{array}{ll}
g_1&g_2
\end{array}\right]}\longrightarrow E_{1}\longrightarrow 0\end{array}$$
and $$\begin{array}{l} \eta:{\quad}0\longrightarrow
S\stackrel{\left[\begin{array}{l} 1\\-1
\end{array}\right]}\longrightarrow P(2)\oplus P(n+1)\longrightarrow
E_{2}\longrightarrow 0\end{array}$$ are unique non-trivial
extensions
 of $S$ by $E_{1}$ and $E_{2},$ respectively.

Based on the claim above, let $M$ be a non-trivial extension of $S$
by $\oplus mE_{~1}\bigoplus\oplus rE_{~2},$ then it is of the form
$P\bigoplus(\oplus (m-a)E_{1}\oplus\oplus (r-b)E_{2})$ for some $
0\leq a\leq m,0\leq b\leq r,$ where $P$ is a pre-projective.

Obviously, the non-trivial extension of $S$ by $\oplus mE_{1}$
(resp.$\oplus rE_{2}$) has a unique form $P(1)\oplus P(n)\bigoplus
\oplus (m-1)E_{1}$ (resp. $P(2)\oplus P(n+1)\bigoplus \oplus
(r-1)E_{2}$).

Since
$$\begin{array}{l}|Ext^1_{\llz}(\oplus mE_{1},S)_{P(1)\oplus P(n)\bigoplus \oplus (m-1)E_{1}}|
=|Ext^1_{\llz}(\oplus mE_{2},S)_{P(2)\oplus P(n+1)\bigoplus \oplus
(m-1)E_{2}}|=q^m-1,\\
g^{L}_{MN}=\frac{|Ext^{1}(M,N)_L||Aut(L)|}{|Aut(M)||Aut(N)||Hom|(M,N)|},\end{array}$$
we obtain
$$g^{P(1)\oplus P(n)\bigoplus \oplus (m-1)E_{1}}_{\oplus mE_{1},S}=
g^{P(2)\oplus P(n+1)\bigoplus \oplus (m-1)E_{2}}_{\oplus
mE_{2},S}=q^{m-1}(q-1).$$ Applying Lemma \ref{l:6.1}, we get
$$\begin{array}{l}
g^{P(1)\oplus P(n)\oplus(m-1)E_{1}\bigoplus \oplus rE_{2}}_{\oplus
mE_{1}\bigoplus \oplus rE_{2},S}=q^{m-1}(q-1)|Hom_{\llz}(P(1)\oplus
P(n),\oplus rE_{2})|,\\[4mm]
g^{P(2)\oplus P(n+1)\oplus (m-1)E_{~2}\bigoplus \oplus
rE_{1}}_{\oplus mE_{2}\bigoplus\oplus
rE_{1},S}=q^{m-1}(q-1)|Hom_{\llz}(P(2)\oplus P(n+1),\oplus
rE_{1})|.
\end{array}$$

Moreover, by $(A)$ and Lemma 16.13 in \cite{KW}, we have that
$$\begin{array}{l}V(S,\oplus mE_{2}\bigoplus\oplus
rE_{1};P(1)\oplus P(n)\oplus(m-1)E_{1}\bigoplus \oplus rE_{2}),\\
V(S,\oplus mE_{2}\bigoplus\oplus rE_{1};P(2)\oplus P(n+1)\oplus
(m-1)E_{2}\bigoplus \oplus rE_{1})\end{array}$$ have the purity
property.

Because $Ext^1(\oplus mE_{1},\oplus mE_{1})=0,$ we get
$\overline{\co}_{\oplus mE_{1}}=\bbe_{\udim \oplus mE_{1}}.$

So,
$$({p_3^{-1}|_{\overline{\co}_{\oplus
mE_{1}}\times\overline{\co}_S}})(P(1)\oplus
P(n)\oplus(m-1)E_{1})=P^{dim Hom(S,P(1)\oplus
P(n)\oplus(m-1)E_{1})-1}.$$

Thus,
$$({p_3^{-1}|_{\overline{\co}_{\oplus
mE_{1}}\times\overline{\co}_S})(P(1)\oplus
P(n)\oplus(m-1)E_{1}})$$ has purity property.

We now consider the cases $m\neq 0, r\neq 0.$ By $(A)$ and
Mayer-Vietoris sequence, we may consider the case $a=m,b=r,$ only.

  Consider the short exact
sequence
\begin{eqnarray}\label{f:6.3}{\qquad}0\longrightarrow S\longrightarrow P\longrightarrow \oplus
mE_{1}\bigoplus \oplus rE_{2}\longrightarrow 0.\end{eqnarray}

 Applying $Hom(\oplus
mE_{1}\bigoplus \oplus rE_{2},~)$ to (\ref{f:6.3}) we have

$$\begin{array}{l}0\longrightarrow Hom(\oplus mE_{1}\bigoplus \oplus
rE_{2},S)\longrightarrow
Hom(\oplus mE_{1}\bigoplus \oplus rE_{2},P)\\
\longrightarrow Hom(\oplus mE_{1}\bigoplus \oplus rE_{2},\oplus
mE_{1}\bigoplus \oplus rE_{2}) \longrightarrow
Ext^1(\oplus mE_{1}\bigoplus \oplus rE_{2},S)\\
\longrightarrow Ext^1(\oplus mE_{1}\bigoplus \oplus
rE_{2},P)\longrightarrow Ext^1(\oplus mE_{1}\bigoplus \oplus
rE_{2},\oplus mE_{1}\bigoplus \oplus rE_{2})\longrightarrow
0\end{array}.$$

Since
$$ Hom(\oplus mE_{1}\bigoplus \oplus
rE_{2},S)=Hom(\oplus mE_{1}\bigoplus \oplus rE_{2},P)=0,$$ we get
$$m^2+r^2\leq m+r.$$

Thus
$$m=r=1.$$

We now have the short exact sequence
\begin{eqnarray}\label{f:6.4}{\qquad} 0\longrightarrow S\longrightarrow P\longrightarrow
E_{1}\oplus E_{2}\longrightarrow 0. \end{eqnarray}

If we also have the short exact sequence
\begin{eqnarray}\label{f:6.5}{\qquad\qquad} 0\longrightarrow S\longrightarrow P\stackrel{f}\longrightarrow
X\longrightarrow 0 \end{eqnarray} where $X\ncong E_{1}\oplus E_{2}.$

Suppose that $P$ is indecomposable. By the projectivity of simple module $S,$ we
obtain $Hom(P,S)=0.$ In addition, by AR-theory, we have $Hom(X,S)=0$ and
$Hom(X,P)=0.$

Since
$$\begin{array}{l}<\udim X, \udim S>=dim Hom(X,S)-dim Ext^1(X,S)=-dim Ext^1(X,S),\\
<\udim (E_1+E_2), \udim S>=dim Hom(E_1+E_2,S)-dim Ext^1(E_1+E_2,S)=-dim Ext^1(E_1+E_2,S),\\
\text{ and }<\udim X, \udim S>=<\udim (E_1+E_2), \udim S>,
\end{array}$$
we have $Ext^1(X,S)=k^2$

Applying $Hom(X,)$ to (\ref{f:6.5}) again, we have
$$\begin{array}{l}
0\longrightarrow Hom(X,S)\longrightarrow Hom(X,P)\longrightarrow
Hom(X,X)\\
\stackrel{f^*}\longrightarrow Ext^1(X,S)\longrightarrow
Ext^1(X,P)\longrightarrow Ext^1(X,X)\longrightarrow0.\end{array}$$

Because  $f^*$ is injective, we get $dim_kEnd(X)\leq 2.$ Thus
$X\notin\overline{\co}_{E_{1}\oplus E_{2}}.$

Suppose that $P$ is decomposable.
 Since $$\lr{\dz, \udim P}=\lr{\dz, \udim
S}=-2,$$ we may assume that
 $P=P_1\oplus P_2$ with $P_1$ and $P_2$ being pre-projective
indecomposable. For convenience, we may assume that $n=4.$ Let $Y_i,
i=1,2,4,5$ be pre-projective indecomposables of dimension vector
$(01211),(10211),(11201),(11210),$ respectively. Then $P\cong
P(1)\oplus Y_1$, or $P\cong P(2)\oplus Y_2$, or $P\cong P(4)\oplus
Y_4$, or $P\cong P(5)\oplus Y_5.$ Without loss of generality, we set
$P=P(1)\oplus Y_1.$

Consider
\begin{eqnarray}\label{f:6.6} 0\longrightarrow S\longrightarrow P\longrightarrow
E_{1}\oplus E_{2}\longrightarrow 0. \end{eqnarray}

Since
$$\begin{array}{l}dim Hom(P(1),E_{1})=\lr{\udim P(1),\udim E_{1}}=1,
 dim Hom(P(1),E_{2})=\lr{\udim P(1),\udim E_{2}}=0,\\
 dim Hom(Y_1,E_{1})=\lr{\udim Y_1,\udim E_{1}}=0,
 dim Hom(Y_1,E_{2})=\lr{\udim Y_1,\udim E_{2}}=1.\end{array}$$
 The short exact seqence (\ref{f:6.6}) turns into  the following form
 $$0\longrightarrow S\longrightarrow P(1)\oplus Y_1\stackrel{\left[\begin{array}{ll}
f_1&0\\0&f_2
\end{array}\right]}\longrightarrow E_{1}\oplus E_{2}
 \longrightarrow 0.$$

 But $dim Im(f_1)\leq  dim P(1)=2<3=dim E_{1}.$ This gives
 a contradiction. Thus $P$ must be indecomposable.

 We first point out that $P$ is the extension of $S$ by
 $E_{1}\oplus E_{2}.$ In fact, $P\cong
 M(x),$ and
$E_{1}\oplus E_{2}\cong M(y),$ where $x_{13}=\left[\begin{array}{l}
1\\0\\1
\end{array}\right],x_{23}=\left[\begin{array}{l}
1\\0\\0
\end{array}\right],x_{43}=\left[\begin{array}{l}
0\\1\\0
\end{array}\right],x_{53}=\left[\begin{array}{l}
0\\1\\1 \end{array}\right];y_{13}=\left[\begin{array}{l} 1\\0
\end{array}\right],y_{23}=\left[\begin{array}{l} 0\\1
\end{array}\right],y_{43}=\left[\begin{array}{l} 1\\0
\end{array}\right],y_{53}=\left[\begin{array}{l} 0\\1
\end{array}\right].$ It is easy to see that $(f_i)_i$ is an
epimorphism from $P$ to $E_1\oplus E_2,$ here
$f_1=f_5=-1,f_2=f_4=1,f_3=\left[\begin{array}{lll} 0&1&-1\\1&0&-1
\end{array}\right].$

From the discussions above, we obtain
$$(p_3^{-1}|_{\overline{\co}_{E_{1}\oplus E_{2}}\times\co_S})(P)=
V(S,E_{1}\oplus E_{2};P).$$

Let $X$ be the set of epimorphisms from $P$ to $E_{1}\oplus E_{2}.$
Then
$$V(S,E_{1}\oplus E_{2};P)=X/{Aut(E_{1}\oplus E_{2})}. $$
Here the action of $Aut_{\llz}(E_{1}\oplus E_{2})$ on $X$ is defined
by
$$h\circ f=hf.$$

Since
$$\begin{array}{l}
Hom_{\llz}(P,E_{1}\oplus E_{2})=\{(f_i)_i| f_1=a, f_2=b, f_4=-a,
f_5=-b, f_3=\left[\begin{array}{lll} 0&-a&a\\b&0&-b
\end{array}\right], a,b\in k\},\\
End_{\llz}(E_{1}\oplus E_{2})=\{(h_i)_i| h_1=c, h_2=d, h_4=c, h_5=d,
h_3=\left[\begin{array}{ll} c&0\\0&d
\end{array}\right], c,d\in k\},

\end{array}$$
we have
$$\begin{array}{l}
X=\{(f_i)_i|(f_i)_i\in Hom_{\llz}(P,E_{1}\oplus E_{2}), a\neq 0, b\neq 0 \},\\
Aut_{\llz}(E_{1}\oplus E_{2})=\{(h_i)_i|(h_i)_i\in
End_{\llz}(E_{1}\oplus E_{2}), c\neq 0, d\neq 0\}.
\end{array}$$

Then $X/Aut(E_{1}\oplus E_{2})$ has only one point correspond to the
orbit of $(1,1,\left[\begin{array}{lll} 0&-1&1\\1&0&-1
\end{array}\right],-1,-1),$ that is,
$$V(S,E_{1}\oplus E_{2};P)=\{1 \text{ point }\}.$$

Thus $({p_3^{-1}|_{\overline{\co}_{E_1\oplus E_2}\times \co_S}})(P)$
has the purity property.

 Case $3: Q=\widetilde{E}_6$

 We only consider $\ct_1$ in 6.3. Other cases can be proved
 similarly.
Let $E_{i}, i=1,2$ be simple objects in the full subcategory $\ct_1$
of $\llz-\text{mod}$ corresponding to the dimension vectors listed
in 6.3.

Since
$$dim Ext^1(E_{1}, S)=-\lr{\udim E_{1},\udim S}=1, \text{ and } dim
Ext^1(E_{2}, S)=-\lr{\udim E_{2},\udim S}=2,$$ we have
$$dim Ext^1(\oplus mE_{1}, S)=m, dim Ext^1(\oplus mE_{2}, S)=2m.$$
Let $X$ and $Y_1$ be indecomposable pre-projective modules of
dimension vectors $(0121010)$ and $(0111000),$ respectively.

Then
$$\begin{array}{l}
\zeta:{\qquad\qquad}0\longrightarrow
S\stackrel{\left[\begin{array}{l} 1\\1\\-1
\end{array}\right]}\longrightarrow P(1)\oplus P(5)\oplus
P(7)\longrightarrow E_{1}\longrightarrow 0
\end{array}$$
and
$$\begin{array}{l}
\zeta_1:{\qquad\qquad}0\longrightarrow S\stackrel{(f_i)_{i\in
I}}\longrightarrow
X\longrightarrow E_{2}\longrightarrow 0,\\
\zeta_1:{\qquad\qquad}0\longrightarrow
S\stackrel{\left[\begin{array}{l} 1\\-1
\end{array}\right]}\longrightarrow P(6)\oplus
Y_1\longrightarrow E_{2}\longrightarrow 0
\end{array}$$
is a basis of $Ext^1(E_{1}, S)$ and $Ext^1(E_{2}, S),$ respectively,
where $f_3=\left[\begin{array}{l} 1\\-1
\end{array}\right],f_i=0,i\neq 3.$

Based on the claim above, we may assume that $P\bigoplus(\oplus
(m-a)E_{1}\oplus\oplus(r-b)E_{2})$ is a non-trivial extension of $S$
by $\oplus mE_{1}\oplus\oplus rE_{2}$ for some $0\leq a\leq m,0\leq
b\leq r,$ where $P$ is a preprojective.

Up to isomorphism, $(P(1)\oplus P(5)\oplus P(7))\bigoplus(\oplus
(m-1)E_{1})$ is the unique non-trivial extension of $S$ by $\oplus
mE_{1}.$

Base on \cite{Z}, we known that for any $X\in \mathbb{E}_{\underline{dim}mE_1}$
we have
$X\in\overline{\co}_{mE_1}.$

By the proof of Lemma 16.13 in \cite{KW}, we get that
$$\begin{array}{l}{p_3|_{\overline{\co}_{\oplus
mE_{1}}\times \co_S}}^{-1}(P(1)\oplus P(5)\oplus P(7)\oplus
(m-1)E_{1})\\=(Hom(S,P(1)\oplus P(5)\oplus P(7)\oplus
(m-1)E_{1})-\{0\})/k^*=P^{dim_kHom(S,P(1)\oplus P(5)\oplus P(7)\oplus
(m-1)E_{1})-1}.\end{array}$$ has the purity property.

Because of $E_{1}=\tau(E_2)$ ,
 we can reduce the problem $p_3^{-1}(P(1)\oplus
P(5)\oplus P(7)\oplus (m-1)E_{1}\bigoplus\oplus rE_{2})$ into the
same problem of the pair $(\tau^{-1}(S),\tau^{-1}(E_{1}))$ .

 Since $E_{2}$ is
non-sincere, we can also reduce the problem of $Ext^1_{\llz}(\oplus
mE_{2},S)$ into the same problem of $Ext^1_{{\llz}'}(\oplus
mE_{2},S),$ where ${\llz}'$ is the path algebra of Dynkin quiver
$D_4:$

$$\xymatrix{&6\ar[dr]&&\\
            & &3&4\ar[l] \\
            &2\ar[ur]  }$$

\begin{center}
Fig.7.1\\
\end{center}

We know that the $AR-$ quiver of $D_4$ is
\begin{center}
\text{ where}
P(3)=S.\qquad\qquad\qquad\qquad\qquad\qquad\qquad\qquad\qquad\qquad\qquad\qquad\qquad\qquad\qquad
\qquad\qquad\qquad\\
\end{center}

{\tiny$$\xymatrix@C=0.3cm@R=0.3cm{&&P(6)={\begin{smallmatrix}
 1& & \\
  &1&0\\
 0& &
 \end{smallmatrix} } \ar[dr]&&Y_1={\begin{smallmatrix}
 0& & \\
  &1&1\\
 1& &
 \end{smallmatrix} } \ar[dr]&&N_1= {\begin{smallmatrix}
 1& & \\
  &0&0\\
 0& &
 \end{smallmatrix} }                \\   
 &P(3)={\begin{smallmatrix}
 0& & \\
  &1&0\\
 0& &
 \end{smallmatrix} }\ar[ur]\ar[dr] \ar[r]
 &P(2)={\begin{smallmatrix}
 0& & \\
  &1&0\\
 1& &
 \end{smallmatrix} } \ar[r]
 &X={\begin{smallmatrix}
 1& & \\
  &2&1\\
 1& &
 \end{smallmatrix} } \ar[ur]\ar[dr]\ar[r]
 &Y_2={\begin{smallmatrix}
 1& & \\
  &1&1\\
 0& &
 \end{smallmatrix} } \ar[r]
 &Z={\begin{smallmatrix}
 1& & \\
  &1&1\\
 1& &
 \end{smallmatrix} } \ar[ur]\ar[dr]\ar[r]
 &N_2={\begin{smallmatrix}
 0 & \\
  &0&0\\
 1& &
 \end{smallmatrix} }                   \\ 
 & &P(4)={\begin{smallmatrix}
 0& & \\
  &1&1\\
 0& &
 \end{smallmatrix} }  \ar[ur]&
 &Y_3={\begin{smallmatrix}
 1& & \\
  &1&0\\
 1& &
 \end{smallmatrix} } \ar[ur]&
 &N_3={\begin{smallmatrix}
 0& & \\
  &0&1\\
 0& &
 \end{smallmatrix} } }$$}
 \begin{center}
 Fig.7.2
\end{center}

 Let $Y,N$ be any $\llz-$ modules, and let $N_1=S_6, N_2=S_2,N_3=S_4.$
  We denote by $W_N$  the set of maps $w:S\rightarrow Y$ such that
either $w=0$ or $Hom_{\llz} (N, Cok w)\neq 0,$ and denote by $W_t$
the set of maps $w:S\rightarrow Y$ which factor through $P(t).$
Suppose $Hom_{\llz}(S,Y)=k^{n+1}$, then the coordinate ring of
$Hom_{\llz}(S,Y)$ is $k[T_0,\cdots,T_n].$

The action of $k^*=Aut(S)$ on $Hom_{\llz}(S,Y)=k^{n+1}$ is defined
by
$$a\circ (w_0,\cdots,w_n)=(w_0a,\cdots,w_na).$$

When $Y=X,$ we have $Hom_{\llz}(S,X)=k^2,$ and the coordinate ring
of $Hom_{\llz}(S,X)$ is $k[T_0,T_1].$ Let $r_{N_j}$ be the defining
relations of $W_{N_j},$ then
$r_{N_1}=T_1,r_{N_2}=T_0,r_{N_3}=T_0-T_1,$ that is,
$W_{N_j}=\mathscr{L}(r_{N_j})$ for $j=1,2,3.$

 Thanks to the Corollary on page 166 in~\cite{R3}, we get
$$
~ {\bf Z}_{E_2,S,X}=(\bba^2-\cup_{j=1}^3\mathscr{L}(r_{N_j}))/k^*
=P^1- {\{[1:0],[0:1],[1:1] \}}.$$

Since
$$dim_k Hom_{\llz}(S,E_{2})=\lr{\udim S,\udim E_{2}}=1,$$
we have
$$ Hom_{\llz}(S,X\bigoplus\oplus (m-1)E_{2})=k^{m+1}.$$

For any $0\neq w\in Hom_{\llz}(S,X\bigoplus\oplus (m-1)E_{2}),$
suppose
$$Cok w=a_1Y_1\oplus a_2Y_2\oplus a_3Y_3\oplus (m-i)E_{2}\oplus
 b_1N_1\oplus b_2N_2\oplus b_3N_3.$$  Then we have
 \begin{eqnarray}\label{f:6.7}{\qquad\qquad}\left\{\begin{array}{lll}
a_1+a_2+a_3&=&i\\a_2+a_3+b_1&=&i\\a_3+a_1+b_2&=&i\\a_1+a_2+b_3&=&i.
 \end{array}\right .\end{eqnarray}

 Then $a_j=b_j,$ for all $j,$ and $b_j=0$, $a_j=0$ if
 $Hom_{\llz}(N_j, Cok w)=0$ for all $j,$ that is, $Cok w\cong \oplus mE_{2}.$
So we have
$r_{Cok w\cong \oplus mE_{2}}=\lr{T_0+T_1}.$

 Thus we obtain
\begin{eqnarray}\label{f:6.8}{\qquad\qquad\qquad}{\bf Z}_{\oplus mE_{2},S,X\bigoplus \oplus
(m-1)E_{2}}=P^{m-1}.\end{eqnarray}

As to the triple  $(\oplus mE_{2},S,P(6)\oplus Y_1\bigoplus \oplus
(m-1)E_{2})$ of modules, by comparing dimension vectors of the
middle term and both ends in exact sequence, similar to $(26),$ we
obtain that $N_i$ is still the test module of a pair  $(P(6)\oplus
Y_1\bigoplus \oplus (m-1)E_{2},S)$ modules.

So, we also have
$$\begin{array}{l}
{\bf Z}_{\oplus mE_2,S,P(6)\oplus Y_1\bigoplus
\oplus (m-1)E_2}=P^{m-1}.
\end{array}$$

Thus $$\begin{array}{l}{\bf Z}_{\oplus mE_{2},S,X\bigoplus \oplus
(m-1)E_{2}}, {\bf Z}_{\oplus mE_2,S,P(6)\oplus Y_1\bigoplus \oplus
(m-1)E_2}\end{array}$$ have the purity property.

 Now, from Fig.7.2, we have$$\begin{array}{l}{p_3|_{\overline{\co}_{\oplus
mE_{2}}\times \co_S}}^{-1}(X\bigoplus \oplus (m-1)E_{2})\\=(Hom(S,X\bigoplus \oplus
(m-1)E_{2})-\{0\})/k^*=P^{dim_kHom(S,X\bigoplus \oplus
(m-1)E_{2})-1}.\end{array}$$

Similarly,
$$\begin{array}{l}{p_3|_{\overline{\co}_{\oplus
mE_{1}}\times \co_S}}^{-1}(P(1)\oplus
P(5)\oplus P(7)\oplus (m-1)E_{1})\\
{\qquad\qquad}=p_3^{-1}(\tau^{-1}(P(1)\oplus P(5)\oplus
P(7))\bigoplus\oplus (m-1)E_{2} )\\
{\qquad\qquad}=p_3^{-1}((Y_1\oplus Y_2\oplus Y_3)\bigoplus\oplus
(m-1)E_{2} )\\
{\qquad\qquad}=P^{dim Hom(X,(Y_1\oplus Y_2\oplus Y_3)\bigoplus\oplus
(m-1)E_{2})-1}.
\end{array}$$

Thus the statement is true if $m=0$ or $r=0$.

We now assume that $m\neq 0,r\neq 0.$ By $(A)$ and Mayer-Vietoris
sequence, we also only need to consider the case $a=m,b=r.$

Consider the short exact sequence
\begin{eqnarray}\label{f:6.9}{\qquad\qquad\qquad}0\longrightarrow S\longrightarrow P\stackrel{g}\longrightarrow \oplus
mE_{1}\bigoplus rE_{2} \longrightarrow 0.\end{eqnarray}

Applying $Hom(\oplus mE_{1}\bigoplus\oplus rE_{2},~)$  to $(\ref{f:6.9}),$
we obtain
$$\begin{array}{l}
0\longrightarrow Hom(\oplus mE_{1}\bigoplus\oplus
rE_{2},S)\longrightarrow Hom(\oplus mE_{1}\bigoplus\oplus
rE_{2},P)\longrightarrow Hom(\oplus mE_{1}\bigoplus\oplus
rE_{2},\oplus mE_{1}\bigoplus\oplus
rE_{2})\\
\stackrel{g^*}\longrightarrow Ext^1(\oplus mE_{1}\bigoplus\oplus
rE_{2},S)\longrightarrow Ext^1(\oplus mE_{1}\bigoplus\oplus
rE_{2},P)\longrightarrow Ext^1(\oplus mE_{1}\bigoplus\oplus
rE_{2},S)\longrightarrow 0
\end{array}$$

Because $g^*$ is an injective, we get $m^2+r^2\leq m+2r.$ Thus $m=1;
r=1, 2.$

When $m=r=1,$  we consider the short exact sequence

 \begin{eqnarray}\label{f:6.10}{\qquad\qquad}0\longrightarrow S\longrightarrow P\longrightarrow E_{1}
 \oplus E_{2}\longrightarrow 0.\end{eqnarray}
Since $\lr{\dz,\udim P}=\lr{\dz,\udim S}=-3,$ $P$ has at most three indecomposable objects. Otherwise, $\lr{\dz,\udim P}\leq -4.$

Suppose that $P\cong P_1\oplus P_2\oplus P_3$ with $P_i$ being
indecomposable for $i=1,2,3.$ By the $AR-$quiver of
$\widetilde{E}_{6},$ we have
$$P\cong P(1)\oplus \tau^{-1}(P(1))\oplus \tau^{-2}(P(1)),\text{~or~}
P(5)\oplus \tau^{-1}(P(5))\oplus
\tau^{-2}(P(5)),\text{~or~}P(7)\oplus \tau^{-1}(P(7))\oplus
\tau^{-2}(P(7)).$$

Without loss of generality, we may assume that $P\cong P(1)\oplus
\tau^{-1}(P(1))\oplus \tau^{-2}(P(1))$. Because there is no
epimorphism from $\tau^{-1}(P(1))$ to $E_{2},$ by $Hom(P(1)\oplus
\tau^{-2}(P(1)),E_{2})=0,$ there is no exact sequence
$$0\longrightarrow S\longrightarrow P(1)\oplus \tau^{-1}(P_1)\oplus \tau^{-2}(P_1)
\longrightarrow E_{1}\oplus E_{2}\longrightarrow 0$$

 Suppose that $P$ is an indecomposable extension of
$S$ by $E_{1}\oplus E_{2}.$ Let $X$ be a module in
$\overline{\co}_{E_{1}\oplus E_{2}}\setminus{\co}_{E_{1}\oplus
E_{2}}$ , and we have
 the following short exact sequence
\begin{eqnarray}\label{f:6.11}{\qquad\qquad}0\longrightarrow S\longrightarrow P\longrightarrow X
\longrightarrow 0.\end{eqnarray}

Applying $Hom(~,S)$ (resp. $Hom(X,~)$) to $(\ref{f:6.11})$, we have
$$\begin{array}{l}0\longrightarrow Hom(X,S)\longrightarrow Hom(P,S)\longrightarrow Hom(S,S)\\
\longrightarrow Ext^1(X,S)\longrightarrow Ext^1(P,S)\longrightarrow
Ext^1(S,S)\longrightarrow 0\end{array}$$ and
$$\begin{array}{l}0\longrightarrow Hom(X,S)\longrightarrow Hom(X,P)
\longrightarrow Hom(X,X)\\
\longrightarrow Ext^1(X,S)\longrightarrow Ext^1(X,P)\longrightarrow
Ext^1(X,X)\longrightarrow 0.\end{array}$$

Now, it follows from $Hom(P,S)=0,Hom(X,P)=0$ and $Ext(X,S)=k^3$ that
$ dim End(X)\leq 3.$
By $X\in\overline{\co}_{E_{1}\oplus E_{2}}\setminus{\co}_{E_{1}\oplus
E_{2}}$, we obtain
 $dim End(X)=3.$

 If $X\cong X_1\oplus X_2$ and $X_i\in\ct_1,i=1,2,$  then we have $X\cong E_{1}\oplus E_{2},
 $ a contradiction.

If $X\cong X_1\oplus X_2$ and $X_1\in \cp_{prep},$ $X_2\in\ct_1$
(resp. $X_1\in \ct_1$ ,$X_2\in\cp_{prei}$), then by
Proposition~~\ref{p:5.1.2}
 we have that $E_{1}\oplus E_{2}$ is an
extension $X_1$ by $X_2.$ Since $\ct_1$ is an abelian subcategory of $\Lambda-$ module category,
this is a contradiction.

If $X\cong X_1\oplus X_2$ and $X_1\in\cp_{prep}, X_2\in\cp_{prei},$ by $dim End(X)=3,$
then $dim Hom(X_1,X_2)=1.$ Suppose $X_1=\tau^{-a}P(j),$ then we have
$$1=dim Hom(\tau^{-a}P(j),X_2)=dim Hom(P(j),\tau^{a}X_2).$$
It implies that $X_1=\tau^{-a}P(2), \tau^{-a}P(4),\text{ or
}\tau^{-a}P(6).$ According to the $AR-$quiver of $\widetilde{E}_6,$
we have $0\leq a\leq 2$ and $X_2=\tau^{2-a}I(2),\tau^{2-a}I(4) \text{ or }\tau^{2-a} I(6).$

If $X\cong X_1\oplus X_2\oplus X_3,$ we have $Hom(X_i,X_j)=0$ for
all $i\neq j.$ Suppose that there is one regular term in $X_i,
i=,1,2,3$, and the other two terms are pre-projective or
pre-injective, then, similar to the above discussions, we get
$X\notin\overline{\co}_{E_{1}\oplus E_{2}}.$ Suppose $X_3$ (resp.
$X_1$) is pre-injective (resp. pre-projective), and the other two
terms are pre-projective (resp. pre-injective). By
$X\in\overline{\co}_{E_{1}\oplus E_{2}}$ and $dim End(X)=3, dim
End(E_{1}\oplus E_{2})=2$, we have
\begin{eqnarray}\label{f:6.12}{\qquad\qquad}0\longrightarrow X_1\oplus X_2\longrightarrow
E_{1}\oplus E_{2}\longrightarrow X_3\longrightarrow0.\end{eqnarray}

Applying $Hom(~,X_3)$ to $(\ref{f:6.12}),$ we obtain
$$\begin{array}{l}0\longrightarrow Hom(X_3,X_3)\longrightarrow Hom(E_{1}\oplus E_{2},X_3)
\longrightarrow Hom(X_1\oplus X_2,X_3)\\
\longrightarrow Ext^1(X_3,X_3)\longrightarrow Ext^1(E_{1}\oplus E
_{2},X_3)\longrightarrow Ext^1(X_1\oplus X_2,X_3)\longrightarrow
0.\end{array}$$

Because  $Hom(X_1\oplus X_2,X_3)=0,$  we can deduce that $dim
Hom(E_{1}\oplus E_{2},X_3)=1.$

Applying $Hom(E_{1}\oplus E_{2},~)$ to $(\ref{f:6.12}),$ we also obtain $dim
Hom(E_{1}\oplus E_{2},X_3)\geq 2$ by $Hom(E_{1}\oplus
E_{2},X_1\oplus X_2)=0.$ This gives a contradiction.

Let $X_1$ (resp. $X_2,X_3$) be the pre-projective (resp. regular,
pre-injective) indecomposable module such that $X_2\cong E_{1}$ or
$X_2\cong E_{2}$. Suppose $X_2\cong E_{1},$ then it follows from
$(\ref{f:6.12})$ that
$$0\longrightarrow X_1\longrightarrow
E_{2}\longrightarrow X_3\longrightarrow0.$$ Thus $X_1\cong
\tau^{-1}P(1)$, or $X_1\cong \tau^{-1}P(5)$, or $X_1\cong
\tau^{-1}P(7).$

Suppose $X_2\cong E_{2},$ then it follows from $(\ref{f:6.12})$ that
$$0\longrightarrow X_1\longrightarrow
E_{1}\longrightarrow X_3\longrightarrow0.$$ Thus $X_3\cong I(1)$,or $\tau^2I(1),$or
$X_3\cong I(5),$ or $\tau^2I(5),$ or $X_3\cong I(7),$ or $\tau^2I(7).$

We, thus, only need to consider the following two cases:

(I) Let $\;   X\cong X_1\oplus X_2, X_1\in\cp_{prep}, X_2\in
\cp_{prei},\text{ and } Hom(X_1,X_2)=1.$ Then, $X\cong
\tau^{2-a}I(2)\oplus\tau^{-a} P(2),$ where $0\leq a\leq 2$.
Without loss of generality, we may only consider $a=2,$ that is, $X\cong I(2)\oplus\tau^{-2} P(2).$  Since $dim Ext^1(\tau^{-2}P(2),P(1))=1,$ we
fixed an exact sequence as follows
$$0\longrightarrow P(1) \stackrel{p}\longrightarrow P\stackrel{q}\longrightarrow\tau^{-2}P(2)\longrightarrow 0.$$
Since
 $$dim Hom(P(1), I(2))=dim Hom(S,P(1))=1,$$
 we choose $u_0\in
Hom(S,P(1))$ (resp. $v_0\in Hom(P(1), I(2))$) such that $\{u_0\}$
(resp. $\{v_0\}$) is the basis of $Hom(S,P(1))$ (resp. $Hom(P(1),
I(2))$).

For any $(u,v)\in W (S,P(1); I(2)),$ we have following commutative
diagram

\begin{equation*}
\begin{array}{c c c c c c}

0\longrightarrow &S & \stackrel{u}{\longrightarrow} &
P(1)&\stackrel{v}\longrightarrow & I(2) \longrightarrow 0
\\
& \parallel &  &p\downarrow&&\downarrow \\
0  \longrightarrow & S & \stackrel{pu}{\longrightarrow }& P
&\stackrel{\pi}\longrightarrow & Cok (pu)\longrightarrow 0.
\end{array}
\end{equation*}

Based on the Snake Lemma , we get $Cok(pu)\stackrel{w}\cong
I(2)\oplus\tau^{-2}P(2).$
 Thus $(pu,w\pi)\in W (S,
I(2)\oplus\tau^{-2}P(2); P)$ by $Hom(I(2),P)=0.$

Conversely, let $(f,g)\in W (S, I(2)\oplus\tau^{-2}P(2); P),$ by the
projectivity of $P(1),$ there exist morphisms $u=au_0, p'=bp$ such
that the following diagram commutes

\begin{equation*}
\begin{array}{c c c c c c}

0\longrightarrow &S & \stackrel{u}{\longrightarrow} &
P(1)&\stackrel{v_0}\longrightarrow & I(2) \longrightarrow 0
\\
& \parallel &  &p'\downarrow&&\left[\begin{array}{l}1\\0\end{array}\right]\downarrow \\
0  \longrightarrow & S & \stackrel{f}{\longrightarrow }& P
&\stackrel{g}\longrightarrow &
I(2)\oplus\tau^{-2}P(2)\longrightarrow 0.
\end{array}
\end{equation*}

Thus $(bu,b^{-1}v_0)\in W (S,I(2);P(1))$.  We deduce that
$$V(S,I(2)\oplus\tau^{-2}P(2);P)=V(S,I(2);P(1))=\{ 1 \text{~point~}\}.$$

(II) Let $X\cong X_1\oplus X_2\oplus X_3,
X_1\in\cp_{prep},X_2\in\cp_{regular}, X_3\in \cp_{prei},\text{ and }
Hom(X_i,X_j)=0, i\neq j.$ Based on above consider, X has nine kinds of possibility. Without loss of generality, we may let $X\cong \tau^{-1} P(1)\oplus E_{1}
\oplus \tau I(1)$.

Since
$$Hom(P,\tau^{-1} P(1))=Hom(P,\tau^{-1} P(5))=Hom(P,\tau^{-1}
P(7))=0$$ and
$$
Ext^1(I(1),S)=Ext^1(I(5),S)=Ext^1(I(7),S)=0,$$

there is no exact sequence of the form (\ref{f:6.11}).

From the discussions  above, we obtain
$$\begin{array}{l}{p_3|_{\overline{\co}_{
E_{1}\bigoplus E_{2}}\times
\co_S}}^{-1}(P)\\
{\qquad\qquad\qquad}=V(S,E_{1}\bigoplus E_{2};P)\text{ \d{$\cup$}}
\text{ \d{$\cup$}}_{a=0}^2V(S,\tau^{2-a}I(2)\oplus\tau^{-a}P(2);P)\\
{\qquad\qquad\qquad}\text{ \d{$\cup$}}\text{ \d{$\cup$}}_{a=0}^2V(S,\tau^{-a}P(4)\oplus\tau^{2-a} I(4);
P)\text{ \d{$\cup$} }\text{ \d{$\cup$}}_{a=0}^2V(S,\tau^{-a}P(6)\oplus\tau^{2-a}I(6); P).\end{array}$$

We first point out that $P$ is the extension of $S$ by
 $E_{1}\oplus E_{2}.$ In fact, $P\cong
 M(x),$ and
$E_{1}\oplus E_{2}\cong M(y),$ where $x_{12}=\left[\begin{array}{l}
0\\1
\end{array}\right],x_{23}=\left[\begin{array}{ll}
1&0\\0&0\\0&1\\0&0
\end{array}\right],x_{43}=\left[\begin{array}{ll}
0&0\\1&0\\0&0\\0&1
\end{array}\right],x_{54}=\left[\begin{array}{l}
0\\1 \end{array}\right], \\
x_{63}=\left[\begin{array}{ll}
0&1\\1&1\\1&1\\0&1
\end{array}\right], x_{76}=\left[\begin{array}{l}
0\\1 \end{array}\right] ;y_{12}=\left[\begin{array}{l} 1\\0
\end{array}\right],y_{23}=\left[\begin{array}{ll}1& 0\\0&0\\0&1
\end{array}\right],y_{43}=\left[\begin{array}{ll} 0&0\\1&0\\0&1
\end{array}\right],y_{54}=\left[\begin{array}{l} 1\\0
\end{array}\right],\\
 y_{63}=\left[\begin{array}{ll} 1&0\\1&0\\0&1
\end{array}\right], y_{76}=\left[\begin{array}{l} 1\\0
\end{array}\right].$ It is easy to see that $(f_i)_i$ is an
epimorphism from $P$ to $E_{1}\oplus E_{2},$ where $f_1=1,
f_2=\left[\begin{array}{ll} 1&1\\1&0
\end{array}\right], f_3=\left[\begin{array}{llll}
1&0&1&0\\0&1&0&1\\1&-1&0&0
\end{array}\right], f_4=\left[\begin{array}{ll} 1&1\\-1&0
\end{array}\right], f_5=1, \\f_6=\left[\begin{array}{ll} 1&2\\-1&0
\end{array}\right], f_7=2.$

Let $X$ be the set of epimorphisms from $P$ to $E_{1}\oplus E_{2}.$
Then
$$V(S,E_{1}\oplus E_{2};P)=X/{Aut(E_{1}\oplus E_{2})}. $$
Here the action of $Aut_{\llz}(E_{1}\oplus E_{2})$ on $X$ is defined
by
$$h\circ f=hf.$$

Since
$$\begin{array}{l}
Hom_{\llz}(P,E_{1}\oplus E_{2})\\
{\qquad\qquad}=\{(f_i)_i| f_1=a, f_2=\left[\begin{array}{ll}
b_1&a\\b_2&0
\end{array}\right], f_3=\left[\begin{array}{llll}
b_1&0&a&0\\0&a&0&b_1\\b_2&-b_2&0&0
\end{array}\right], f_4=\left[\begin{array}{ll} a&b_1\\-b_2&0
\end{array}\right], f_5=b_1, \\
{\qquad\qquad}f_6=\left[\begin{array}{ll} a&a+b_1\\-b_2&0
\end{array}\right], f_7=a+b_1; a, b_1, b_2\in k\},\\
End_{\llz}(E_{1}\oplus E_{2})\\
{\qquad\qquad}=\{(h_i)_i|h_1=c, h_2=\left[\begin{array}{ll} c&0\\0&d
\end{array}\right],h_3=\left[\begin{array}{lll}
c&0&0\\0&c&0\\0&0&d
\end{array}\right], h_4=\left[\begin{array}{ll} c&0\\0&d
\end{array}\right], h_5=c, \\
{\qquad\qquad}h_6=\left[\begin{array}{ll} c&0\\0&d
\end{array}\right], h_7=c; c, d \in k\},

\end{array}$$

 we have
$$\begin{array}{l}
X=\{(f_i)_i|(f_i)_i\in Hom_{\llz}(P,E_{1}\oplus E_{2}), a\neq 0, b_1\neq 0, b_2\neq 0, a+b_1\neq 0 \},\\
Aut_{\llz}(E_{1}\oplus E_{2})=\{(h_i)_i|(h_i)_i\in
End_{\llz}(E_{1}\oplus E_{2}), c\neq 0, d\neq 0\}.
\end{array}$$

Hence $X/{Aut(E_{1}\oplus E_{2})}$ has the points correspond to the
orbits of $(f_1=1, f_2=\left[\begin{array}{ll} b_1/a&1\\1&0
\end{array}\right], \\f_3=\left[\begin{array}{llll}
b_1/a&0&1&0\\0&1&0&b_1/a\\1&-1&0&0
\end{array}\right], f_4=\left[\begin{array}{ll} 1&b_1/a\\-1&0
\end{array}\right], f_5=b_1/a,
f_6=\left[\begin{array}{ll} 1&1+b_1/a\\-1&0
\end{array}\right], \\
f_7=1+b_1/a); b_1/a\neq 0,-1,$ that is,
$$V(S,E_{1}\oplus E_{2};P)=\bba^1\setminus \{2 \text{ point }\}.$$

Now we give the specific characterization of these two points. Let
$z_{12}=\left[\begin{array}{l} 1\\0
\end{array}\right], z_{23}=\left[\begin{array}{ll} 1&0\\0&0\\0&1
\end{array}\right], z_{43}=\left[\begin{array}{l} 0\\1\\-1
\end{array}\right], z_{54}=0, x_{63}=\left[\begin{array}{ll} 1&0\\1&0\\0&1
\end{array}\right], z_{76}=\left[\begin{array}{l} 1\\0
\end{array}\right].$ We then have  $\tau^{-2}(P(4))=M(z),$ and then $\tau^{-2}(P(4))\oplus I(4)=M(w),$
where $w_{12}=\left[\begin{array}{l} 1\\0
\end{array}\right], w_{23}=\left[\begin{array}{ll} 1&0\\0&0\\0&1
\end{array}\right], w_{43}=\left[\begin{array}{ll} 0&0\\1&0\\-1&0
\end{array}\right], w_{54}=\left[\begin{array}{l} 0\\1
\end{array}\right], w_{63}=\left[\begin{array}{ll} 1&0\\1&0\\0&1
\end{array}\right], w_{76}=\left[\begin{array}{l} 1\\0
\end{array}\right].$ Let $f=(f_i)_i$ be a morphism corresponds to $b_1=0$ in $X/{Aut(E_{1}\oplus E_{2})}.$
It is easy to see that $f$ is an epimorphism from $P$ to $\tau^{-2}(P(4))\oplus I(4)$. Thus this $f$ is an one point in two points above. Similarly, we can get another point. So, we prove $$\begin{array}{l}{p_3|_{\overline{\co}_{
E_{1}\bigoplus E_{2}}\times
\co_S}}^{-1}(P)\\
{\qquad\qquad\qquad}=V(S,E_{1}\bigoplus E_{2};P)\text{ \d{$\cup$}}
\text{ \d{$\cup$}}_{a=0}^2V(S,\tau^{2-a}I(2)\oplus\tau^{-a}P(2);P)\\
{\qquad\qquad\qquad}\text{ \d{$\cup$}}\text{ \d{$\cup$}}_{a=0}^2V(S,\tau^{-a}P(4)\oplus\tau^{2-a} I(4);
P)\text{ \d{$\cup$} }\text{ \d{$\cup$}}_{a=0}^2V(S,\tau^{-a}P(6)\oplus\tau^{2-a}I(6); P)\\
{\qquad\qquad\qquad}=\bba^1\cup\{7 \text{ point}\}.\end{array}$$

By the proof of Lemma 16.13 in \cite{KW}, we have ${p_3|_{\overline{\co}_{E_{1}\oplus E_{2}}\times
\co_S}}^{-1}(P)$ has
purity property.

 Suppose that
$P\cong P_1\oplus P_2$ with $P_i$ is indecomposable for $i=1,2.$ By
the $AR-$quiver of $\widetilde{E}_{6}$ and \cite{Z}, we have
$$\begin{array}{l}
P\cong P(1)\oplus \tau^{-2}(P(2)),\text{~or~} P(5)\oplus
\tau^{-2}(P(4)),\text{~or~}P(7)\oplus \tau^{-2}(P(6));\\
P\cong \tau^{-2}P(1)\oplus \tau^{-1}(P(2)),\text{~or~}
\tau^{-2}P(5)\oplus \tau^{-1}(P(4)),\text{~or~}\tau^{-2}P(7)\oplus
\tau^{-1}(P(6)).
\end{array}$$

Without loss of generality, we assume that $P\cong P(1)\oplus
\tau^{-2}(P(2)).$ We claim that $P$ is the nontrivial extension of
$S$ by $E_{1}\oplus E_{2}.$

 Let
$x_{12}=0, x_{23}=\left[\begin{array}{l} 1\\1\\1
\end{array}\right], x_{43}=\left[\begin{array}{ll} 0&0\\1&0\\0&1
\end{array}\right], x_{54}=\left[\begin{array}{l} 1\\0
\end{array}\right], x_{63}=\left[\begin{array}{ll} 1&0\\0&0\\0&1
\end{array}\right],x_{76}=\left[\begin{array}{l} 1\\0
\end{array}\right].$ We then have  $\tau^{-2}(P(2))=M(x).$

It is easy to see that $f=(f_i)_i$ (resp.
$\vartheta=(\vartheta_i)_i$ ) is the epimorphism from $
\tau^{-2}(P(2))\oplus P(1)$ (resp.$\tau^{-2}(P(2))$ to $E_{1}\oplus
E_{2}$ (resp.$E_{2}$ ),where $f_1=1, f_2=\left[\begin{array}{ll}
1&0\\1&1
\end{array}\right], f_3=\left[\begin{array}{llll}
0&0&1&0\\1&0&0&1\\1&-1&0&0
\end{array}\right], f_4=\left[\begin{array}{ll} 0&1\\-1&0
\end{array}\right], f_5=1, f_6=\left[\begin{array}{ll}0&1\\ 1&0
\end{array}\right], f_7=1$ (resp. $\vartheta_1=0, \vartheta_2=2, \vartheta_3=\left[\begin{array}{lll}
1&1&0
\end{array}\right], \vartheta_4=\left[\begin{array}{ll} 1&0
\end{array}\right], \vartheta_5=0, \vartheta_6=\left[\begin{array}{ll} 1&0
\end{array}\right], \vartheta_7=0$). So we get that  $f$ is the epimorphism from $P(1)\oplus \tau^{-1}(P(2))$ to
$E_{1}\oplus E_{2}.$ The claim is proved.

Suppose that there is  an $X$ satisfying $E_{1}\oplus E_{2}\ncong X$
and $X\in\overline{\co}_{E_{1}\oplus E_{2}}$ such that there is a
short exact sequence
\begin{eqnarray}\label{f:6.13}{\qquad\qquad}0\longrightarrow S\longrightarrow  P(1)\oplus
\tau^{-2}(P(2))\longrightarrow X \longrightarrow 0.\end{eqnarray}

Applying $Hom(X,~)$ to $(\ref{f:6.13}),$ we get
$$\begin{array}{l}0\longrightarrow Hom(X,S)\longrightarrow Hom(X,P(1)\oplus
\tau^{-2}(P(2)))
\longrightarrow Hom(X,X)\\
\longrightarrow Ext^1(X,S)\longrightarrow Ext^1(X,P(1)\oplus
\tau^{-2}(P(2)))\longrightarrow Ext^1(X,X)\longrightarrow
0.\end{array}$$
 Since $X$ not contains $S$ as a summand, we get $Hom (X,S)=0£¬$ and $$dim Ext^1(X,S)=-\lr{\delta,\udim S}=3.$$

   If $Hom(X,P(1)\oplus \tau^{-2}(P(2)))=0.$ then we have $dim
End(X)=3.$ Similar to the case of the indecomposable $P$ , by AR-quier of $\widetilde{E}_6,$ we deduce
that there does not exist an $X$ such that $(\ref{f:6.13})$ holds.

Now let $Hom(X,P(1)\oplus \tau^{-2}(P(2)))\neq 0.$ There is an
indecomposable summand $\tau^{-a}P(j)$ of $X$ such that
$Hom(\tau^{-a}P(j), P(1)\oplus \tau^{-2}(P(2)))\neq 0.$ Set
$X_1=\tau^{-a}P(j)$ and $X=X_1\oplus X_2.$

If $X_1=P(1),$ by $Hom(P(1),X_2)=0$ , we have
$$0\longrightarrow S\longrightarrow \tau^{-2}(P(2))\longrightarrow X_2\longrightarrow 0,$$
and
$$\begin{array}{l}{p_3|_{\overline{\co}_{X_2}\times
\co_S}}^{-1}(\tau^{-2}(P(2)))\cong {p_3|_{\overline{\co}_{P(1)\oplus
X_2}\times \co_S}}^{-1}(P(1)\oplus\tau^{-2}(P(2))).
\end{array}$$

If $X_1\ncong P(1),$ $X_1$ must be non-projective.
 According to
the $AR-$ quiver of $\widetilde{E}_6,$ we have $X_1=\tau^{-a}P(j);
j=1,2,3,4,5,6,7; a=1,2.$

Suppose $a\neq 2,j\neq 2.$ For any $f\in Hom(P(1),\tau^{-a}P(j) ),$
we know that $f$ is not an epimorphism. It follows from $ Hom(\tau^{-2}P(2),\tau^{-a}P(j))\neq 0 , \udim \tau^{-a}P(j))\prec \delta,$ and the AR-quiver of $\widetilde{E}_6$ that $a= 2,j= 2,$ and $\udim
X_2= \udim I(2).$ Because $Ext^1(I(1),S)=0,$ it follows from
$(\ref{f:6.13})$ that $X_2=I(2),$ that is, the following short exact
sequence holds

 $$0\longrightarrow S\longrightarrow P(1)\oplus\tau^{-2}P(2)\longrightarrow
 I(2)\oplus\tau^{-2}(P(2))\longrightarrow 0.$$

From the discussions above, we obtain
$$\begin{array}{l}{p_3|_{\overline{\co}_{E_{1}\oplus E_{2}}\times
\co_S}}^{-1}(P(1)\oplus\tau^{-2}P(2))\\
{\qquad\qquad\qquad}=V(S,E_{1}\oplus
E_{2};P(1)\oplus\tau^{-2}P(2))\text{ \d{$\cup$} }
V(S,I(2)\oplus\tau^{-2}P(2);P(1)\oplus\tau^{-2}P(2))\\
{\qquad\qquad\qquad\qquad}\text{ \d{$\cup$}
}_{[X_2]}{p_3|_{\overline{\co}_{P(1)\oplus X_2\times
\co_S}}}^{-1}(P(1)\oplus\tau^{-2}(P(2))),\end{array}$$ where
\d{$\cup$} denotes the disjoint union.

Obviously,
$$V(S,I(2)\oplus\tau^{-2}P(2);P(1)\oplus\tau^{-2}P(2))=\{1 \text{
point }\}.$$

Moreover, $${p_3|_{\overline{\co}_{P(1)\oplus X_2\times
\co_S}}}^{-1}(P(1)\oplus\tau^{-2}(P(2)))$$ has the purity property
by \cite{L4}.

 Let
$z_{12}=\left[\begin{array}{l} 0\\1
\end{array}\right], z_{23}=\left[\begin{array}{ll}
1&0\\1&0\\1&0\\0&1
\end{array}\right], z_{43}=\left[\begin{array}{ll}
0&0\\1&0\\0&1\\0&0
\end{array}\right], z_{54}=\left[\begin{array}{l} 1\\0
\end{array}\right], z_{63}=\left[\begin{array}{ll}
1&0\\0&0\\0&1\\0&0
\end{array}\right],\\
z_{76}=\left[\begin{array}{l} 1\\0
\end{array}\right].$ Then $\tau^{-2}P(2)\oplus P(1)\cong M(z).$

Since
$$\begin{array}{l}
Hom_{\llz}(\tau^{-2}P(2)\oplus P(1),E_{1}\oplus E_{2})\\
{\qquad\qquad}=\{(f_i)_i| f_1=a, f_2=\left[\begin{array}{ll}
b_1&a\\b_2&0
\end{array}\right], f_3=\left[\begin{array}{llll}
b_1&0&0&a\\b_1&-b_1&0&0\\0&0&b_2&0
\end{array}\right], f_4=\left[\begin{array}{ll} -b_1&0\\0&b_2
\end{array}\right], f_5=-b_1, \\
{\qquad\qquad}f_6=\left[\begin{array}{ll} b_1&0\\0&b_2
\end{array}\right], f_7=b_1; a, b_1, b_2\in k\},\\

\end{array}$$

let $X$ be the set of epimorphisms from $\tau^{-2}P(2)\oplus P(1)$
to $E_{1}\oplus E_{2}.$ Then we have
$$\begin{array}{l}
X=\{(f_i)_i|(f_i)_i\in Hom_{\llz}(P,E_{1}\oplus E_{2}), a\neq 0,
b_1\neq 0, b_2\neq 0 \} .
\end{array}$$

Hence $X/{Aut(E_{1}\oplus E_{2})}$ has points correspond to the
orbits of $(f_1=a/b_1, f_2=\left[\begin{array}{ll} 1&a/b_1\\1&0
\end{array}\right], \\f_3=\left[\begin{array}{llll}
1&0&0&a/b_1\\1&-1&0&0\\0&0&1&0
\end{array}\right], f_4=\left[\begin{array}{ll} -1&0\\0&1
\end{array}\right], f_5=-1,
f_6=\left[\begin{array}{ll} 1&0\\0&1
\end{array}\right], \\
f_7=1); a/b_1\neq 0.$

Let $D$ be an one point corresponds to $a/b_1=0$ in the morphism above. Then,
$$V(S,E_{1}\oplus E_{2};P)=\bba^1\setminus \{D\}.$$
It is easy to see that this point is the point just corresponds to an epimorphism from $P(1)\oplus\tau^{-2}P(2)$ to $I(2)\oplus\tau^{-2}P(2).$

Thus $$\begin{array}{l}{p_3|_{\overline{\co}_{E_{1}\oplus E_{2}}\times
\co_S}}^{-1}(P(1)\oplus\tau^{-2}P(2))\\
{\qquad\qquad\qquad}=\bba^1\text{ \d{$\cup$}
}_{[X_2]}{p_3|_{\overline{\co}_{P(1)\oplus X_2\times
\co_S}}}^{-1}(P(1)\oplus\tau^{-2}(P(2))),\end{array}$$ where
\d{$\cup$} denotes the disjoint union.

By the proof of Lemma 16.13 in \cite{KW},
we have ${p_3|_{\overline{\co}_{E_{1}\oplus E_{2}}\times
\co_S}}^{-1}(P(1)\oplus\tau^{-2}P(2))$ has
the purity property.

When $m=1, r=2,$ let $P$ be the extension of $S$ by
$E_{1}\bigoplus\oplus 2E_{2}.$ Since $\lr{\dz, \udim P}=\lr{\dz,
\udim S}=-3,$ we can deduce that $P$ has $3$ indecomposable direct
summands at most.

If $P$ is also the extension of $S$ by $X$ with
$X\in\overline{\co}_{E_{1}\bigoplus\oplus 2E_{2}}\backslash
{\co}_{E_{1}\bigoplus\oplus 2E_{2}},$ that is,

\begin{eqnarray}\label{f:6.14}{\qquad\qquad}0\longrightarrow S\longrightarrow  P\longrightarrow X \longrightarrow 0.\end{eqnarray}
Since $\lr{\dz, \udim P_0}<0$ for each pre-projective $P_0,$ and $\lr{\dz, \udim P-\udim S}=0,$
we deduce $P$ does not contain $S$ as a direct summand. It implies $Hom(P,S)=0.$

Applying $Hom(~,S)$ and $Hom(X,~)$ to $(\ref{f:6.14}),$ we get $Hom(X,S)=0$ and $dim
End(X)\leq 5+dim Hom(X,P).$

Suppose $P$ is an indecomposable. Because of the pre-projective component is directed part,
it follows that $ Hom(X,P)=0.$ Thus $dim End(X)\leq 5.$ But $dim End(E_{1}\bigoplus\oplus 2E_{2})=5,$ so
there does not exist an
 $X\in\overline{\co}_{E_{1}\bigoplus\oplus
2E_{2}}\backslash {\co}_{E_{1}\bigoplus\oplus 2E_{2}}$ such that
$(\ref{f:6.14})$ holds. Then, similar to the $m=r=1$ case,
$${p_3|_{\overline{\co}_{E_{1}\bigoplus\oplus
2E_{2}}\times \co_S}}^{-1}(P)=V(S, E_{1}\bigoplus\oplus 2E_{2};P)$$
has the purity property.

Suppose that $P$ can be decomposed into  two indecomposable objects,
say, $P\cong P_1\oplus P_2$ such that
$$\begin{array}{l}P\cong
\tau^{-1}P(1)\oplus \tau^{-3}P(2)\text{ or } \tau^{-1}P(5)\oplus
\tau^{-3}P(4)\text{ or }\tau^{-1}P(7)\oplus
\tau^{-3}P(6)\\
P\cong \tau^{-2}P(2)\oplus \tau^{-3}P(1)\text{ or }
\tau^{-2}P(4)\oplus \tau^{-3}P(5)\text{ or } \tau^{-1}P(6)\oplus
\tau^{-3}P(7).\end{array}$$ Without loss generality, we may assume
that $P\cong \tau^{-1}P(1)\oplus \tau^{-3}P(2).$

If $Hom(X,P)=0,$ it follows from $dim End(X)\leq 5$ that $X\notin
\overline{\co}_{E_{1}\bigoplus\oplus 2E_{2}}.$

If $Hom(X,P)\neq 0,$ similar to the case of $m=r=1,$  we get $X\cong
\tau^{-1}P(1)\oplus X_2$ and $Hom(X_2,P )=0$ or $X\cong X_1\oplus
\tau^{-3}P(2).$

When $X\cong X_1\oplus \tau^{-3}P(2),$ we have $ X_1=\tau I(5)\oplus
\tau I(7).$ Because of $Hom(\tau^{-1}P(1),\tau^{-3}P(2))=0,$ we obtained
$$g^{\tau^{-1}P(1)}_{\tau I(5)\oplus \tau
I(7),S}=g^{\tau^{-1}P(1)\oplus\tau^{-3}P(2)}_{\tau I(5)\oplus \tau
I(7)\bigoplus\tau^{-3}P(2),S}=1,$$ that is,
$$V(S,\tau I(5)\oplus \tau
I(7)\bigoplus\tau^{-3}P(2);\tau^{-1}P(1)\oplus\tau^{-3}P(2))=\{1
\text{ point }\}.$$

When $X\cong \tau^{-1}P(1)\oplus X_2,$ we then have $dim
End(X)=5+dim Hom (X,P)=6.$ Thus $X_2$ must be decomposable, and $dim
End(X_2)= 4$.

Suppose $X_2\cong X_{21}\oplus X_{22}$ and
$X_{21}\in\cp_{prep},X_{22}\in\cp_{prei}.$ It follows from $dim
End(X_2)= 4$ that $dim End(X_{22})+dim Hom(X_{21},X_{22})\leq 3. $

On the other hand, by Proposition~\ref{p:5.1.2}, we have
$$0\longrightarrow \tau^{-1}P(1)\oplus X_{21}\longrightarrow E_{1}\bigoplus\oplus 2E_{2}
\longrightarrow X_{22}\longrightarrow 0.$$

Applying $Hom(~,X_{22})$  (resp. $Hom(E_{1}\bigoplus\oplus
2E_{2},~)$) to the exact sequence above,
 we get
$dim Hom_{\llz}(E_{1}\bigoplus\oplus 2E_{2},X_{22})\leq 4$ (resp.
$dim Hom(E_{1}\bigoplus\oplus 2E_{2},X_{22})\geq 5$) by
$Hom(\tau^{-1}P(1), X_2)=k$ (resp. $Hom(E_{1}\bigoplus\oplus
2E_{2},\tau^{-1}P(1)\oplus X_{21})=0$).

It implies that $X_2$ contains a summand $\widetilde{X}_2$ with $\widetilde{X}_2\in\ct_1.$

Assume  that  $X_2=E_{1}\oplus X_2',$ then there does not exist any
pre-projective direct summand $M$ of $X_2'$ such that $Hom(M,P)=0$
and there exists an epimorphism $f\in Hom(P,M) $ by the $AR-$ quiver
of $\widetilde{E}_6$.

If $X_2'$ does not contain a summand of regular module, it follows
from Proposition~\ref{p:5.1.2} that
$$0\longrightarrow \tau^{-1}P(1)\longrightarrow \oplus 2 E_{2}\longrightarrow X_2'\longrightarrow 0. $$

Applying $Hom( , \oplus 2 E_{2}),$ we get $4=dim End(\oplus 2
E_{2})\leq 2 =dim Hom(\tau^{-1}P(1),\oplus 2 E_{2} )=dim Hom(P(1),\oplus 2 E_{1} ).$  Therefore,
$X\cong \tau^{-1}P(1)\oplus E_{1}\oplus E_{2}\oplus \tau I(1). $

Assume  that  $X_2=E_{2}\oplus X_2'.$ When $X_2'$ contains a
pre-projective summand $M$ such that $Hom(M,P)=0$ and $\exists f\in
Hom(P,M) $ with $f$ to be an epimorphism, we can deduce $M\cong
\tau^{-3}P(1)$ by the $AR-$ quiver of $\widetilde{E}_6.$ Thus
$X\cong\tau^{-1}P(1)\oplus E_{2}\oplus\tau^{-3}P(1)\oplus I(5)\oplus
I(7) $ and $dim End(X)\geq 7,$ which gives a contradiction.

When $X_2'$ has only pre-injective summands, we get
$$0\longrightarrow \tau^{-1}P(1)\longrightarrow E_{1}\oplus  E_{2}
\longrightarrow X_2'\longrightarrow 0. $$

Applying $Hom(~,E_{1}\oplus E_{2})$  to the exact sequence above, we
get a contradiction $2 =dim End(E_{1}\oplus E_{2})\leq 1=dim
Hom(\tau^{-1}P(1),E_{1}\oplus E_{2}).$

From the discussions above, we get $X\cong \tau^{-1}P(1)\oplus
E_{1}\oplus E_{2}\oplus \tau I(1).$

Thus
$$\begin{array}{l}
{p_3|_{\overline{\co}_{ E_{1}\bigoplus\oplus 2E_{2}}\times
\co_S}}^{-1}(\tau^{-1}P(1)\oplus \tau^{-3}P(2))=V(S,
E_{1}\bigoplus\oplus 2E_{2};\tau^{-1}P(1)\oplus
\tau^{-3}P(2))\\
{\qquad\qquad\qquad\qquad\qquad\qquad}\text{\d{$\cup$}} V(S,\tau
I(5)\oplus \tau I(7);\tau^{-1}P(1)
\text{\d{$\cup$}}({p_3|_{\overline{\co}_{ E_{1}\bigoplus
E_{2}}\times \co_{P(2)}}})^{-1}(\tau^{-3}P(2))
 .\end{array}$$

By $m=r=1,$ we deduce that $({p_3|_{\overline{\co}_{ E_{1}\oplus
E_{2}}\times \co_{P(2)}}})^{-1}(\tau^{-3}P(2))$ has the purity
property. Thus ${p_3|_{\overline{\co}_{ E_{1}\bigoplus\oplus
2E_{2}}\times \co_S}}^{-1}(\tau^{-1}P(1)\oplus \tau^{-3}P(2))$ has
also the purity property.

Suppose $P=P_1\oplus P_2\oplus P_3.$ Without loss of generality, let
$P=\tau^{-1}P(1)\oplus\tau^{-2}P(1)\tau^{-3}P(1). $ Because $S_1$ is
a summand of  $top(E_{1}/Im(f))$ for any $f\in Hom(P_1\oplus
P_2\oplus P_3,E_{1}),$ we obtain that $P_1\oplus P_2\oplus P_3$ is
not an extension of $S$ by $E_{1}\bigoplus\oplus 2E_{2}.$ The proof
for the Case $3$ is complete.

Case $4: Q=\widetilde{E}_7.$

 We only consider $\ct_2$ in 6.3. The other cases can be proved
 similarly.
Let $E_i', i=1,2,3,$ be simple objects in the full subcategory
$\ct_2$ of $\llz-\text{mod}$ corresponding to dimension vectors
listed in 6.3.

Since $$dim_k Ext^1(E_1',S)=dim_k Ext^1(E_2',S)=1, dim_k
Ext^1(E_3',S)=2,$$ we have
$$dim_k Ext^1(\oplus mE_1',S)=dim_k
Ext^1(\oplus mE_2',S)=m, dim_k Ext^1(\oplus mE_3',S)=2m.$$

In addition, we know that
$$\begin{array}{l}
\theta:{\qquad\qquad}0\longrightarrow
S\stackrel{\left[\begin{array}{l} 1\\1\\-1
\end{array}\right]}\longrightarrow P(1)\oplus P(7)\oplus
P(8)\longrightarrow E_1'\longrightarrow 0,\\
\theta':{\qquad\qquad}0\longrightarrow
S\stackrel{\left[\begin{array}{l} 1\\-1
\end{array}\right]}\longrightarrow P(2)\oplus
P(6)\longrightarrow E_2'\longrightarrow 0
\end{array}$$
is a basis of $Ext^1(E_1', S)$ and $Ext^1(E_2', S),$ respectively.

Obviously, $P(1)\bigoplus P(7)\bigoplus P(8)\bigoplus\oplus (m-1)E_1'$
(resp. $P(2)\bigoplus P(6)\bigoplus\oplus (m-1)E_2'$) is a non-trivial
extension of $S$ by $\oplus mE_1'$ (resp. $\oplus mE_2'$).

Based on the claim above, we may assume that $P\bigoplus(\oplus
(m-a)E_1'\oplus\oplus (p-b)E_2'\oplus\oplus (r-c)E_3')$ is a
non-trivial extension of $S$ by $\oplus mE_1'\bigoplus\oplus
pE_2'\bigoplus\oplus rE_3',$ for some $a,b,c,$ where $P$ is
pre-projective.

Up to isomorphism, $(P(1)\oplus P(7)\oplus P(8))\bigoplus\oplus
(m-1)E_1'$ (resp. $(P(2)\oplus P(6))\bigoplus\oplus (m-1)E_2'$) is the
unique non-trivial extension of $S$ by $\oplus mE_1'$ (resp. $\oplus
mE_2'$).

Since
$$\begin{array}{l}
{p_3}^{-1}|_{\overline{\co}_{\oplus m E_1'}\times \co_S}(P(1)\bigoplus
P(7)\bigoplus P(8)\bigoplus\oplus(m-1)E_1')=P^{2m}, \\{p_3|}^{-1}_{\overline{\co}_{\oplus m
E_2'}\times \co_S}(P(2)\bigoplus P(6)\bigoplus\oplus(m-1)E_2')=P^{m},
\end{array}$$ we can deduce that they have the purity property.

 Because $E_3'$ is non-sincere, we can reduce the study of $Ext^1_{\llz}(\oplus
mE_3',S)$ to the study of $Ext^1_{{\llz}'}(\oplus mE_3',S),$ where
${\llz}'$ is the path algebra of Dynkin quiver $D_4$ determined by
vertices $3,4,5$ and $8.$

By  Case 3, we get

$$\left\{\begin{array}{l}
{\bf Z}_{\oplus mE_3',S,X_1\bigoplus \oplus
(m-1)E_3'}=(\bba^{m+1}-\cup_{j=1}^3\mathscr{L}(r_{N_j}))/k^*,\\
{\bf Z}_{\oplus mE_3',S,(P(8)\oplus Y_1)\bigoplus \oplus
(m-1)E_3'}=(\bba^{m+1}-\cup_{j=1}^3\mathscr{L}(r_{N_j}))/k^*,
\end{array}\right.$$ and
$$\left\{\begin{array}{l}p_3^{-1}|{\overline{\co}_{\oplus m E_3'\times \co_S }}(X_1\bigoplus \oplus
(m-1)E_3')\\
{\qquad\qquad\qquad}=P^{dim_kHom(S,X_1\bigoplus \oplus
(m-1)E_3')-1},\\
p_3^{-1}|{\overline{\co}_{\oplus m E_3'\times \co_S }}(P(8)\oplus
Y_1\bigoplus \oplus
(m-1)E_3')\\
{\qquad\qquad\qquad}=P^{dim_kHom(S,P(8)\oplus Y_1\bigoplus \oplus
(m-1)E_3')-1}.\end{array}\right.$$

Thus $p_3^{-1}|{\overline{\co}_{\oplus m E_3'\times \co_S
}}(X_1\bigoplus \oplus (m-1)E_3'))$ and
$p_3^{-1}|{\overline{\co}_{\oplus m E_3'\times \co_S }}(P(8)\oplus
Y_1\bigoplus \oplus (m-1)E_3')$ have the purity property.

Similarly, the statement is true if $mp\neq 0, r=0$ or $pr\neq 0,
m=0$ or $rm\neq 0, p=0.$

We now assume that $mpr\neq 0.$ Because of $(A)$ and Mayer-Vietoris
sequence, we only need to consider $a=m, b=p, c=r.$

Consider the short exact sequence
\begin{eqnarray}\label{f:6.15}{\qquad\qquad\qquad}0\longrightarrow S\longrightarrow P\longrightarrow \oplus
mE_1'\oplus\oplus pE_2'\oplus\oplus rE_3'\longrightarrow 0.\end{eqnarray}

Applying $H(\oplus mE_1'\bigoplus\oplus pE_2'\bigoplus\oplus rE_3',~
),i=1,2,3$ to $(\ref{f:6.15}),$ we get
$$m^2+p^2+r^2\leq m+p+2r.$$
Thus $m=p=1,1\leq r\leq 2.$

When $m=p=r=1.$

Consider the short exact sequence
$$0\longrightarrow S\longrightarrow P\longrightarrow E_1'\oplus E_2'\oplus E_3'
\longrightarrow 0. $$

 Because $\lr{\dz,\udim P}=-4,$ we know that $P$ contains four indecomposable objects at most.

Suppose that $P$ is indecomposable. Let
$X\in\overline{\co}_{E_1'\oplus E_2'\oplus E_3'}$ such that $P$ is
an extension of $S$ by $X.$ Then we have the following exact
sequence

\begin{eqnarray}\label{f:6.16}{\qquad\qquad\qquad}0\longrightarrow S\longrightarrow
P\longrightarrow X \longrightarrow 0. \end{eqnarray}

Applying $Hom(~,S)$ (resp. $Hom(X,~)$) to $(\ref{f:6.16}),$ we get from  $
Hom(P,S)=0$ and $Hom(X,P)=0$ that
$$Hom(X,S)=0, Ext^1(X,S)=k^4, dim End(X)=4. $$

If $X\cong X_1\oplus X_2$ with $X_i$ indecomposable for all $i=1,2$,
we have $X_1\in\cp_{prep} ,X_2\in\cp_{prei}$ and $Hom(X_1,X_2)=2.$
Thus $X_1\cong \tau^{-3}P(3) \text{ or } \tau^{-3}P(5)$, and
$X\cong\tau^{-3}P(3)\oplus I(3) \text{ or } \tau^{-3}P(5)\oplus
I(5).$

Without loss of generality, we may assume that
$X\cong\tau^{-3}P(3)\oplus I(3).$

Fix one exact sequence
$$0\longrightarrow P(1)\longrightarrow P\longrightarrow \tau^{-3}P(3)\longrightarrow 0.$$

Because $dim Hom(P(1),I(3))=1,dim Hom(S,P(1))=1,$ and in the same
way as in $(I)$ and Case $3$, we get
$$V(S,\tau^{-3}P(3)\oplus
I(3);P)\cong V(S,I(3);P(1))=\{ 1 \text{ point }\}.$$

Suppose $X\cong X_1\oplus X_2\oplus X_3$ with $X_i$ indecomposable
for all $i=1,2,3.$ Then assume that there is no regular module in
$X_i.$ Without loss of generality, let $X_1, X_2\in\cp_{prep}$ and
$X_3\in \cp_{prei}.$ Because of Proposition~\ref{p:5.1.2} and $dim
End(X)=4$, we have
$$0\longrightarrow X_1\oplus X_2\longrightarrow E_1'
\oplus E_2'\oplus E _3'\longrightarrow X_3\longrightarrow 0.$$

Applying $Hom(~,X_3)$ (resp. $Hom(E_1' \oplus E_2'\oplus E_3',~)$)
to the sequence above, we get
$$\begin{array}{l}dim _k Hom_{\llz}(E_1'
\oplus E_2'\oplus E_3', X_3)\leq 2,\\
dim _k Hom_{\llz}(E_1' \oplus E_2'\oplus E_3', X_3)\geq
3.\end{array}$$ This gives a contradiction. Thus
 $X_1\in\cp_{prep} ,X_2\in\ct_2,
X_3\in\cp_{prei}.$

When $X_2\cong E_1',$ since $Hom(X,P)=0,$ it follows from the
$AR-$quiver of $\widetilde{E}_7$ that $X_1\cong \tau^{-5}P(1)$ or
$X_1\cong \tau^{-5}P(7).$ Thus

$$X\cong \tau^{-5}P(1)\oplus E_1'\oplus\tau^2 I(1) \text{ or } X\cong
\tau^{-5}P(7)\oplus E_1'\oplus\tau^2 I(7). $$

By $End(\tau^{-5}P(1)\oplus E_1'\oplus\tau^2
I(1))=End(\tau^{-5}P(7)\oplus E_1'\oplus\tau^2 I(7))=k^3,$ we get a
contradiction.

When $X_2\cong E_2',$ we get that $X_1$ is isomorphic to
$\tau^{-3}P(8)$ or $\tau^{-i}P(j)$ for $i=3,4,5,7; j=1,7.$

Suppose $X_1\cong\tau^{-3}P(8),$ and $P$ is the extension of $S$ by
$\tau^{-3}P(8)\oplus E_2'\oplus I(8).$ Since
$$\begin{array}{l}dim_k Hom_{\llz}(P,\tau^{-3}P(8))=dim_k Hom_{\llz}(P,E_2')=1, dim_k Hom_{\llz}(P,I(8))=2,\\
dim_k Hom_{\llz}(\tau^{-3}P(8),I(8))=1,
Aut_{\llz}(\tau^{-3}P(8)\oplus E_2'\oplus I(8))=\\
\{\left[\begin{array}{lll}
a&0&0\\0&b&0\\d&0&c
\end{array}\right]| a\in Aut(\tau^{-3}P(8)),b\in Aut(E_2'),c\in Aut(I(8)),d\in Hom(\tau^{-3}P(8),I(8))\},\end{array}$$
we may deduce  $Hom(P,\tau^{-3}P(8)\oplus E_2'\oplus I(8))$ have two orbits under action on $Aut_{\llz}(\tau^{-3}P(8)\oplus E_2'\oplus I(8)),$ One of the orbit corresponds to the morphism which factor through $\tau^{-3}P(8).$ Thus
$$V(S,\tau^{-3}P(8)\oplus E_2'\oplus I(8);
P)=\{ two~~ points\}.$$

Suppose $X_1\cong \tau^{-3}P(1)$ or $X_1\cong \tau^{-3}P(7)$ , then
we know that $X$ has four indecomposable summands at least.

Suppose $X_1\cong \tau^{-4}P(1)$ or $X_1\cong \tau^{-4}P(7),$ then
we have $End(X)=3.$

Suppose $X_1\cong \tau^{-5}P(1)$ or $X_1\cong \tau^{-5}P(7),$ then
  $P$ has the direct summand $I(1)$ or $I(7).$

Suppose $X_1\cong \tau^{-7}P(1)$ or $X_1\cong \tau^{-7}P(7),$ we
then have $X\cong \tau^{-7}P(1)\oplus E_2'\oplus I(1)$ or $X\cong
\tau^{-7}P(7)\oplus E_2'\oplus I(7).$ Since $Ext^1_{\llz}(I(1),
S)=Ext^1_{\llz}(I(7), S)=0,$ we get that $P$ is decomposable.

 Therefore, in this case, $P$ is the only possible
extension of $S$ by $\tau^{-3}P(8)\oplus E_2'\oplus I(8).$

When $X_2\cong E_3',$ we get
$$\begin{array}{l}X\cong
\tau^{-3}P(2)\oplus E_3'\oplus I(2),\text{ or } \tau^{-3}P(1)\oplus
E_3'\oplus \tau^4I(1), \text{ or } \tau^{-6}P(1)\oplus E_3'\oplus
\tau I(1),\\X\cong \tau^{-3}P(6)\oplus E_3'\oplus I(6), \text{ or }
\tau^{-3}P(7)\oplus E_3'\oplus \tau^4I(7), \text{ or }
\tau^{-6}P(7)\oplus E_3'\oplus \tau I(7).\end{array}$$

Suppose $X\cong\tau^{-3}P(2)\oplus E_3'\oplus I(2),\text{ or
}\tau^{-3}P(6)\oplus E_3'\oplus I(6).$ Because
$Ext^{1}(I(2),S))=Ext^{1}(I(6),S))=0$, $I(2)$ or $I(6)$ is a direct
summands of $P.$

Suppose $$X\cong\tau^{-3}P(1)\oplus E_3'\oplus \tau^4I(1), \text{ or
}  \tau^{-3}P(7)\oplus E_3'\oplus \tau^4I(7), $$ then we get $dim
End(X)=3.$

Suppose $$X\cong\tau^{-6}P(1)\oplus E_3'\oplus \tau I(1), \text{ or
} \tau^{-6}P(7)\oplus E_3'\oplus \tau I(7).$$ Then it follows from
$Ext^1(\tau I(1), S)=Ext^1(\tau I(7), S)=0$ that $P$ is
decomposable.

Thus $P$ is not an extension of $S$ by $X_1\oplus E_3'\oplus X_3$
either.

If $X\cong X_1\oplus X_2\oplus X_3 \oplus X_4$ with $X_i$
indecomposable for all $i=1,2,3,4$, then we have that
$\{X_4,X_3,X_2,X_1\}$ is an orthogonal exceptional sequence. Thus
$\fk{C}(X_4,X_3,X_2,X_1)$ is isomorphic to $k\widetilde{A}_3-mod.$
Because $X\in\overline{\co}_{E_1'\oplus E_2'\oplus E_3'}$ and $dim
End (X)=4,$ we see that $\{X_1,X_2,X_3,X_4\}$ are simple objects in
$\fk{C}(X_4,X_3,X_2,X_1)$, and $X_1\cong \tau^{-3}P(7),X_2\cong
E_2',X_3\cong E_3',X_4\cong\ I(7)$ or $X_1\cong
\tau^{-3}P(1),X_2\cong E_2',X_3\cong E_3',X_4\cong\ I(1)$. Therefore
$P$ must be decomposable. This also gives a contradiction.

From above, we get
$$\begin{array}{l}(p_3|_{\overline{\co}_{E_1'\oplus E_2'\oplus E_3'}
\times\co_S})^{-1}(P)\\
{\qquad}=V(S,E_1'\oplus E_2'\oplus E_3';P)
\text{ \d{$\cup$} }V(S,\tau^{-3}P(3)\oplus I(3);P) \\
{\qquad\quad}\text{ \d{$\cup$} }V(S,\tau^{-3}P(5)\oplus
I(5);P)\text{ \d{$\cup$} }V(S, \tau^{-3} P(8)\oplus E_2'\oplus I(8);
P).\end{array}$$

Since
$$\begin{array}{l}dim_k Hom_{\llz}(P,E_1')=dim_k Hom_{\llz}(P,E_3')=1, dim_k Hom_{\llz}(P,E_2')=0,\\
Aut_{\llz}(E_1'\oplus E_2'\oplus E_3')=\{\left[\begin{array}{lll}
a&0&0\\0&b&0\\0&0&c
\end{array}\right]| a\in Aut(E_1'),b\in Aut(E_2'),c\in Aut(E_3')\},\end{array}$$
we get
 $$V(S,E'_1\oplus E_2'\oplus E_3';P)=\{1~~ point\}.$$ Thereby $(p_3|_{\overline{\co}_{E_1'\oplus
E_2'\oplus E_3'} \times\co_S})^{-1}(P)$ has the purity property.

Suppose $P\cong P_1\oplus P_2$ with $P_i$ indecomposable. Without
loss of generality, we may assume that $P\cong P(7)\oplus
\tau^{-3}P(5).$

Let $X\in\overline{\co}_{E_1'\oplus E_2'\oplus E_3'}\setminus
\co_{E_1'\oplus E_2'\oplus E_3'}$ such that $P(7)\oplus
\tau^{-3}P(5)$ is an extension of $S$ by $X.$

If $Hom(X,P(7)\oplus \tau^{-3}P(5))=0,$ then similar to the case $P$
indecomposable, we can deduce that $P(7)\oplus \tau^{-3}P(5)$ is not
an extension of $S$ by $X.$

If $Hom(X,P(7)\oplus \tau^{-3}P(5))\neq 0,$ we then get $X\cong
P(7)\oplus X_2$ or $X\cong \tau^{-3}P(5)\oplus I(5)$ by the
$AR-$quiver of $\widetilde{E}_7.$

 When $X\cong P(7)\oplus X_2,$ it is easy to see that $V(S,X;P(7)\oplus \tau^{-3}P(5))=V(S,X_2;
\tau^{-3}P(5)).$

When $X\cong \tau^{-3}P(5)\oplus I(5),$ we get $V(S,X;P(7)\oplus
\tau^{-3}P(5))=V(S,I(5); P(7))=\{ 1 \text{ point }\}.$

Thus
$$\begin{array}{l}(p_3|_{\overline{\co}_{E_1'\oplus E_2'\oplus E_3'}
\times\co_S})^{-1}(P(7)\oplus \tau^{-3}P(5))\\
{\qquad} = V(S,E_1'\oplus E_2'\oplus E_3';P(7)\oplus \tau^{-3}P(5))
 \text{ \d{$\cup$} }_{[X_2]}(p_3|_{\overline{\co}_{X_2}
\times\co_S})^{-1}(\tau^{-3}P(5))\text{ \d{$\cup$} }
V(S,I(5);P(7)).\end{array}$$

Let
$x_{12}=\left[\begin{array}{l} 1\\0
\end{array}\right], x_{23}=\left[\begin{array}{ll}
1&0\\0&1\\0&0
\end{array}\right], x_{34}=\left[\begin{array}{lll}
1&0&0\\0&1&0\\0&0&1\\0&0&0
\end{array}\right], x_{54}=\left[\begin{array}{ll} 0&0\\0&0\\1&0\\0&1
\end{array}\right], x_{65}=
\left[\begin{array}{l}
0\\1
\end{array}\right],\\x_{76}=0,
x_{84}=\left[\begin{array}{ll} 1&1\\0&1\\1&0\\1&1
\end{array}\right],$ and $y_{12}=\left[\begin{array}{l} 1\\0
\end{array}\right], y_{23}=\left[\begin{array}{ll}
1&0\\0&1\\0&0
\end{array}\right], y_{34}=\left[\begin{array}{lll}
1&0&0\\0&0&0\\0&1&0\\0&0&1
\end{array}\right], \\y_{54}=\left[\begin{array}{lll} 0&0&0\\1&0&0\\0&1&0\\0&0&1
\end{array}\right], y_{65}=
\left[\begin{array}{ll}
1&0\\0&1\\0&0
\end{array}\right],y_{76}=\left[\begin{array}{l}
1\\0
\end{array}\right],
y_{84}=\left[\begin{array}{ll} 1&0\\1&0\\0&0\\0&1
\end{array}\right],$ then $\tau^{-3}P(5)=M(x),\\E'_1\oplus E'_2\oplus E'_3=M(y).$

 Since
$$\begin{array}{l}
 Hom_{\llz}(\tau^{-3}P(5)\oplus P(7),E'_{1}\oplus E'_{2}\oplus E'_{3})\\
{\qquad\qquad}=\{(f_i)_i| f_1=a, f_2=\left[\begin{array}{ll}
a&0\\0&b
\end{array}\right], f_3=\left[\begin{array}{lll}
a&0&0\\0&b&b\\0&0&c
\end{array}\right], f_4=\left[\begin{array}{lllll} a&0&0&0&0\\0&0&0&a&d\\0&b&b&-b&0\\0&0&c&0&0
\end{array}\right], \\f_5=\left[\begin{array}{lll}
0&a&d\\b&-b&0\\c&0&0
\end{array}\right],
f_6=\left[\begin{array}{ll} a&d\\-b&0
\end{array}\right], f_7=d,f_8=\left[\begin{array}{ll} a&a\\c&0
\end{array}\right]; a, b, c,d\in k\},\\

\end{array}$$ and $$\begin{array}{l}
Aut_{\llz}(E'_{1}\oplus E'_{2}\oplus E'_{3})\\
{\qquad\qquad}=\{(h_i)_i| h_1=a, h_2=\left[\begin{array}{ll}
a&0\\0&b
\end{array}\right], h_3=\left[\begin{array}{lll}
a&0&0\\0&b&0\\0&0&c
\end{array}\right], h_4=\left[\begin{array}{llll} a&0&0&0\\0&a&0&0\\0&0&b&0\\0&0&0&c
\end{array}\right], h_5=\left[\begin{array}{lll}
a&0&0\\0&b&0\\0&0&c
\end{array}\right], \\
{\qquad\qquad}h_6=\left[\begin{array}{ll} a&0\\0&b
\end{array}\right], h_7=a,h_8=\left[\begin{array}{ll} a&0\\0&c
\end{array}\right]; a, b, c\in k, abc\neq 0\},\\
\end{array}$$ let $eHom_{\llz}(\tau^{-3}P(5)\oplus P(7),E'_{1}\oplus E'_{2}\oplus E'_{3})$ be the set of epimorphism from $\tau^{-3}P(5)\oplus P(7)$ to $E'_{1}\oplus E'_{2}\oplus E'_{3}$, we then have $$\begin{array}{l}eHom_{\llz}(\tau^{-3}P(5)\oplus P(7),E'_{1}\oplus E'_{2}\oplus E'_{3})\\=\{(f_i)_i|(f_i)_i\in Hom_{\llz}(\tau^{-3}P(5)\oplus P(7),E'_{1}\oplus E'_{2}\oplus E'_{3}),abcd\neq 0\},\text{ and}\\
eHom_{\llz}(\tau^{-3}P(5)\oplus P(7),E'_{1}\oplus E'_{2}\oplus E'_{3})/Aut_{\llz}(E'_{1}\oplus E'_{2}\oplus E'_{3})\\
=\{(\widetilde{f}_i)_i| \widetilde{f}_1=1, \widetilde{f}_2=\left[\begin{array}{ll}
1&0\\0&1
\end{array}\right], \widetilde{f}_3=\left[\begin{array}{lll}
1&0&0\\0&1&1\\0&0&1
\end{array}\right], \widetilde{f}_4=\left[\begin{array}{lllll} 1&0&0&0&0\\0&0&0&1&d/a\\0&1&1&-1&0\\0&0&1&0&0
\end{array}\right], \widetilde{f}_5=\left[\begin{array}{lll}
0&1&d/a\\1&-1&0\\1&0&0
\end{array}\right], \\
{\qquad\qquad}\widetilde{f}_6=\left[\begin{array}{ll} 1&d/a\\-1&0
\end{array}\right], \widetilde{f}_7=d/a,\widetilde{f}_8=\left[\begin{array}{ll} 1&1\\1&0
\end{array}\right]; a,d\in k,d/a\neq0.\}
\end{array}$$ Thus we get $$V(S,E_1'\oplus E_2'\oplus E_3',P(7)\oplus \tau^{-3}P(5))=\bba^1\setminus \{1 \text{ point}\}.$$

By \cite{L4}, we know that $(p_3|_{\overline{\co}_{X_2}
\times\co_S})^{-1}(\tau^{-3}P(5))$ has the purity property. Since
$$(p_3|_{\overline{\co}_{E_1'\oplus E_2'\oplus
E_3'} \times\co_S})^{-1}(P(7)\oplus \tau^{-3}P(5))=\bba^1\text{ \d{$\cup$} }_{[X_2]}(p_3|_{\overline{\co}_{X_2}
\times\co_S})^{-1}(\tau^{-3}P(5))$$
we can deduce that $(p_3|_{\overline{\co}_{E_1'\oplus E_2'\oplus
E_3'} \times\co_S})^{-1}(P(7)\oplus \tau^{-3}P(5)) $ has the purity
property.

Suppose that $P\cong P_1\oplus P_2\oplus P_3$ with $P_i$
indecomposable for all $i=1,2,3.$ Without loss of generality, we
 may assume that $P\cong P(7)\oplus \tau^{-1}P(7)\oplus \tau^{-3}P(6).$

Let $X\in\overline{\co}_{E_1'\oplus E_2'\oplus E_3'}\setminus
\co_{E_1'\oplus E_2'\oplus E_3'}$ such that $P(7)\oplus
\tau^{-1}P(7)\oplus \tau^{-3}P(6)$ is an extension of $S$ by $X.$

If $Hom(X,P(7)\oplus \tau^{-1}P(7)\oplus \tau^{-3}P(6))=0,$ then
similar to the case $P$ indecomposable, we can deduce that
$P(7)\oplus \tau^{-1}P(7)\oplus \tau^{-3}P(6)$ is not an extension
of $S$ by $X.$

If $Hom(X,P(7)\oplus\tau^{-1}P(7)\oplus \tau^{-3}P(6))\neq 0,$ we
then know that $X$  contains a direct summand $ P(7),\text{ or }
\tau^{-1}P(7), \text{ or }\tau^{-2}(P(7)),\text{ or
}\tau^{-2}(P(6)),\text{ or }\tau^{-3}(P(6)).$

When $X\cong \tau^{-2}(P(7))\oplus X_2$ with $X_2$ indecomposable,
we have
$$Hom(\tau^{-2}(P(7)),
X_2)=Hom((P(7)), \tau^2X_2)=0.$$ Thus $dim End(X)=2,$ which gives a
contradiction.

When $X\cong \tau^{-2}(P(7))\oplus X_2\oplus X_3$ with $X_2,X_3$
being indecomposable modules, we then have $X_2\in\ct_2.$

Suppose $X\cong \tau^{-2}(P(7))\oplus E_1'\oplus X_3.$ Then
$End(X)=3,$ it implies $\tau^{-2}(P(7))\oplus E_1'\oplus
X_3\notin\overline{\co}_{E_1'\oplus E_2'\oplus E_3'}\setminus
\co_{E_1'\oplus E_2'\oplus E_3'}.$

Suppose $X\cong \tau^{-2}(P(7))\oplus E_2'\oplus X_3$ (resp. $X\cong
\tau^{-2}(P(7))\oplus E_3'\oplus X_3.$)  Then $I(1)$ (resp. $\tau
I(7)$) is a direct summand of $X_3$. It implies that $I(1)$ (resp.
$\tau I(7)$) to be a direct summand of $P.$ This also contradicts to
the assumption.

Similarly, we may show that $P$ is not an extension of $S$ by $X$ if
$\tau^{-1}P(7)$ or $\tau^{-2}P(6)$ or $\tau^{-3}P(6)$ is a direct
summand of $X.$

Let
$x_{12}=\left[\begin{array}{l} 1\\0
\end{array}\right], x_{23}=\left[\begin{array}{ll}
1&0\\0&1
\end{array}\right], x_{34}=\left[\begin{array}{ll}
1&0\\0&1\\0&0
\end{array}\right], x_{54}=\left[\begin{array}{ll} 0&0\\1&0\\0&1
\end{array}\right], x_{65}=
\left[\begin{array}{l}
0\\1
\end{array}\right],\\x_{76}=0,
x_{84}=\left[\begin{array}{l} 1\\1\\1
\end{array}\right],$
$y_{12}=0, y_{23}=0, y_{34}=0, y_{54}=1, y_{65}=
1,y_{76}=1,
y_{84}=0,$  and
$z_{12}=0, z_{23}=0, z_{34}=1, z_{54}=0, z_{65}=
0,z_{76}=0,
z_{84}=1.$ We then have $\tau^{-3}P(6)=M(x), P(7)=M(y), \tau^- P(7)=M(z),$ and $P(7)\oplus\tau^{-}P(7)\oplus\tau^{-3}P(6)=M(w),$ where $w_{12}=\left[\begin{array}{l} 1\\0
\end{array}\right], w_{23}=\left[\begin{array}{ll}
1&0\\0&1\\0&0
\end{array}\right], w_{34}=\left[\begin{array}{lll}
1&0&0\\0&1&0\\0&0&0\\0&0&0\\0&0&1
\end{array}\right], w_{54}=\left[\begin{array}{lll} 0&0&0\\1&0&0\\0&1&0\\0&0&1\\0&0&0
\end{array}\right], w_{65}=
\left[\begin{array}{ll}
0&0\\1&0\\0&1
\end{array}\right],w_{76}=\left[\begin{array}{l}
0\\1
\end{array}\right],
\\w_{84}=\left[\begin{array}{ll} 1&0\\1&0\\1&0\\0&0\\0&1
\end{array}\right].$

Because
$$\begin{array}{l}
Hom_{\llz}( P(7)\oplus\tau^- P(7)\tau^{-3}P(6),E'_{1}\oplus E'_{2}\oplus E'_{3})\\
{\qquad\qquad}=\{(f_i)_i| f_1=a, f_2=\left[\begin{array}{ll}
a&0\\0&b
\end{array}\right], f_3=\left[\begin{array}{lll}
a&0&0\\0&b&0\\0&0&c
\end{array}\right], f_4=\left[\begin{array}{lllll} a&0&0&0&0\\0&0&a&d&0\\0&b&-b&0&0\\0&0&0&0&c
\end{array}\right], f_5=\left[\begin{array}{lll}
0&a&d\\b&-b&0\\0&0&0
\end{array}\right], \\
{\qquad\qquad}f_6=\left[\begin{array}{ll} a&d\\-b&0
\end{array}\right], f_7=d,f_8=\left[\begin{array}{ll} a&0\\0&c
\end{array}\right]; a, b, c,d\in k\},\\
\end{array}$$ and $f_5$ is not a surjection.
 It implies there is no an epimorphism from $P(7)\oplus\tau^- P(7)\tau^{-3}P(6)$ to $E'_{1}\oplus E'_{2}\oplus E'_{3}.$
Thus
$V(S,E_1'\oplus E_2'\oplus E_3';P_1\oplus P_2\oplus P_3)=\emptyset.$

Suppose $P\cong P_1\oplus P_2\oplus P_3\oplus P_4$ with $P_i$
indecomposable for all $i=1,2,3,4.$ Without loss of generality, we
may assume that $P\cong P(7)\oplus \tau^{-1}P(7)\oplus
\tau^{-2}P(7)\oplus \tau^{-3}P(7).$

Because $Hom(P(7)\oplus \tau^{-1}P(7)\oplus \tau^{-3}P(7),E_2')=0$
and there is not an epimorphism from $\tau^{-2}P(7)$ to $E_2',$ we
have that $P(7)\oplus \tau^{-1}P(7)\oplus \tau^{-2}P(7)\oplus
\tau^{-3}P(7)$ is not an extension  of $S$ by $E_1'\oplus E_2'\oplus
E_3'.$

 When $m=p=1, r=2,$ let $P$ be the extension of $S$ by $E_1'
 \oplus E_2'\oplus 2E_3'.$ Because $\lr{\dz,\udim S}=\lr{\dz,\udim
 P}=-4,$ we can deduce that $P$ has four indecomposable direct summands  at most.

 Assume that $P$ is also the extension of $S$ by $X$ with $X\in \overline{\co}_
 {E_1'
 \oplus E_2'\oplus 2E_3'}\setminus{\co}_
 {E_1'
 \oplus E_2'\oplus 2E_3'}.$ Then there is an
 exact sequence
\begin{eqnarray}\label{f:6.17}{\qquad\qquad}0\longrightarrow S\longrightarrow P\longrightarrow X\longrightarrow 0.\end{eqnarray}

Because $Hom(P,S)=0,$ we deduce that $Hom(X,S)=0, Ext^1(X,S)=k^6.$
Applying $Hom(X,~)$ to (\ref{f:6.17}), we get $dim End(X)\leq 6+dim
Hom(X,P).$

Suppose $P$ is indecomposable. We then have $dim End(X)\leq 6.$ Thus
$X\notin\overline{\co}_
 {E_1'
 \oplus E_2'\oplus 2E_3'}. $ This gives a
 contradiction.

Let $x_{12}=\left[\begin{array}{l} 1\\0
\end{array}\right], x_{23}=\left[\begin{array}{ll}
1&0\\0&1\\0&0\\0&0
\end{array}\right], x_{34}=\left[\begin{array}{llll}
1&0&0&0\\0&1&0&0\\0&0&1&0\\0&0&0&1\\0&0&0&0\\0&0&0&0
\end{array}\right], x_{54}=\left[\begin{array}{llll} 0&0&0&0\\0&0&0&0\\1&0&0&0\\0&1&0&0\\0&0&1&0\\0&0&0&1
\end{array}\right], \\x_{65}=
\left[\begin{array}{ll}
0&0\\0&0\\1&0\\0&1
\end{array}\right],x_{76}=\left[\begin{array}{l}
0\\1
\end{array}\right],
x_{84}=\left[\begin{array}{lll} 1&1&0\\1&0&0\\0&1&0\\0&0&1\\1&1&1\\1&0&1
\end{array}\right],$ and $y_{12}=\left[\begin{array}{l} 1\\0
\end{array}\right], y_{23}=\left[\begin{array}{ll}
1&0\\0&1\\0&0\\0&0
\end{array}\right], \\y_{34}=\left[\begin{array}{llll}
1&0&0&0\\0&0&0&0\\0&1&0&0\\0&0&1&0\\0&0&0&1
\end{array}\right], y_{54}=\left[\begin{array}{llll} 0&0&0&0\\1&0&0&0\\0&1&0&0\\0&0&1&0\\0&0&0&1
\end{array}\right], y_{65}=
\left[\begin{array}{ll}
1&0\\0&1\\0&0\\0&0
\end{array}\right],y_{76}=\left[\begin{array}{l}
1\\0
\end{array}\right],
y_{84}=\left[\begin{array}{lll} 1&0&0\\1&0&0\\0&0&0\\0&1&0\\0&0&1
\end{array}\right],$ we then have $P=M(x), E'_1\oplus E'_2\oplus 2E'_3=M(y).$

Since
$$\begin{array}{l}
 Hom_{\llz}(P,E'_{1}\oplus E'_{2}\oplus 2E'_{3})
=\{(f_i)_i| f_1=x_1, f_2=\left[\begin{array}{ll}
x_1&-x_1\\0&x_2
\end{array}\right], f_3=\left[\begin{array}{llll}
x_1&-x_1&0&0\\0&-x_2&-x_2&-x_2\\0&0&x_3&x_4\\0&0&x_5&x_6
\end{array}\right], \\f_4=\left[\begin{array}{llllll} x_1&-x_1&0&0&0&0\\0&0&0&0&x_1&-x_1\\0&-x_2&-x_2&-x_2&x_2&0\\0&0&x_3&x_4&0&0\\0&0&x_5&x_6&0&0
\end{array}\right], f_5=\left[\begin{array}{llll}
0&0&x_1&-x_1\\-x_2&-x_2&x_2&0\\x_3&x_4&0&0\\x_5&x_6&0&0
\end{array}\right],
f_6=\left[\begin{array}{ll} x_1&-x_1\\x_2&0
\end{array}\right], f_7=-x_1,\\f_8=\left[\begin{array}{lll} 0&x_1&0\\0&x_3&x_4\\0&x_5&x_6
\end{array}\right]; x_1, x_2, x_3,x_4,x_5,x_6\in k\},\\

\end{array}$$ and $$\begin{array}{l}
Aut_{\llz}(E'_{1}\oplus E'_{2}\oplus 2E'_{3})\\
=\{(h_i)_i| h_1=a, h_2=\left[\begin{array}{ll}
a&0\\0&b
\end{array}\right], h_3=\left[\begin{array}{llll}
a&0&0&0\\0&b&0&0\\0&0&c&e\\0&0&f&d
\end{array}\right], h_4=\left[\begin{array}{lllll} a&0&0&0&0\\0&a&0&0&0\\0&0&b&0&0\\0&0&0&c&e\\0&0&0&f&d
\end{array}\right], h_5=\left[\begin{array}{llll}
a&0&0&0\\0&b&0&0\\0&0&c&e\\0&0&f&d
\end{array}\right], \\
h_6=\left[\begin{array}{ll} a&0\\0&b
\end{array}\right], h_7=a,h_8=\left[\begin{array}{lll} a&0&0\\0&c&e\\0&f&d
\end{array}\right]; a,b,c,d,e,f\in k, ab\neq 0,cd-ef\neq 0\},\\
\end{array}$$ let $eHom_{\llz}(P,E'_{1}\oplus E'_{2}\oplus 2E'_{3})$ be the set of epimorphism from $ P$ to $E'_{1}\oplus E'_{2}\oplus 2E'_{3}$, we then have $$\begin{array}{l}eHom_{\llz}(P,E'_{1}\oplus E'_{2}\oplus 2E'_{3})=\{(f_i)_i|(f_i)_i\in Hom_{\llz}(P,E'_{1}\oplus E'_{2}\oplus 2E'_{3}),x_1x_2\neq 0,x_3x_6-x_4x_5\neq 0\},\end{array}$$ and
$$\begin{array}{l}eHom_{\llz}( P,E'_{1}\oplus E'_{2}\oplus 2E'_{3})/Aut_{\llz}(E'_{1}\oplus E'_{2}\oplus 2E'_{3})
=\{(\widetilde{f}_i)_i| \widetilde{ f}_1=1, \widetilde{f}_2=\left[\begin{array}{ll}
1&-1\\0&1
\end{array}\right], \\\widetilde{f}_3=\left[\begin{array}{llll}
1&-1&0&0\\0&-1&-1&-1\\0&0&1&0\\0&0&0&1
\end{array}\right], \widetilde{f}_4=\left[\begin{array}{llllll} 1&-1&0&0&0&0\\0&0&0&0&1&-1\\0&-1&-1&-1&1&0\\0&0&1&0&0&0\\0&0&0&1&0&0
\end{array}\right], \widetilde{f}_5=\left[\begin{array}{llll}
0&0&1&-1\\-1&-1&1&0\\1&0&0&0\\0&1&0&0
\end{array}\right],\\
\widetilde{f}_6=\left[\begin{array}{ll} 1&-1\\1&0
\end{array}\right], \widetilde{f}_7=-1,\widetilde{f}_8=\left[\begin{array}{lll} 0&1&0\\0&1&0\\0&0&1
\end{array}\right]\}
\end{array}$$ Thus we get $V(S,E_1'\oplus E_2'\oplus 2E_3',P)=\{1 \text{ point}\},$ and
 $$\begin{array}{l}(p_3|_{\overline{\co}_{E_1'\oplus E_2'\oplus 2E_3'}
\times\co_S})^{-1}(P) = V(S,E_1'\oplus E_2'\oplus
2E_3';P)\end{array}$$ has the purity property.

Suppose $P\cong P_1\oplus P_2$ with $P_i$ indecomposable. Without
loss of generality, we may set $P\cong\tau^{-1}P(7)\oplus
\tau^{-4}P(5).$

If $Hom(X,P)=0,$ it follows that $dim End(X)\leq 6$ and
$X\notin\overline{\co}_
 {E_1'
 \oplus E_2'\oplus 2E_3'}.$ This gives a
 contradiction.

 Thus $Hom(X,P)\neq 0,$ and we get $\tau^{-4}P(5)\vdash X$ or $\tau^{-1}P(7)\vdash X.$

 When $X\cong \tau^{-4}P(5)\oplus X',$ we deduce that $X\cong
 I(8)\oplus\tau^2I(1)\oplus\tau^{-4}P(5).$ Thus
 $$V(S,I(8)\oplus\tau^2I(1)\oplus\tau^{-4}P(5);\tau^{-1}P(7)\oplus
\tau^{-4}P(5))=V(S,I(8)\oplus\tau^2I(1));\tau^{-1}P(7))=\{ 1 \text{
point }\}.$$

When $X\cong \tau^{-1}P(7)\oplus X_2,$ we get $ dim Hom(X, P)=1$ and
$ dim End(X)=7.$ We now set $X\cong \tau^{-1}P(7)\oplus X_2.$ It
follows from $dim Hom (\tau^{-1}P(7),X_2)=1$ that $dim End(X_2)=5$
and $X_2$ is decomposable.

Assume that $X_2\cong X_{21}\oplus X_{22}$ with
$X_{21}\in\cp_{prep}, X_{22}\in\cp_{prei}.$ By $dim End(X_2)=5,$ we
get $dim End(X_{22})+dim Hom(X_{21},X_{22})\leq 4.$ Based on
Proposition~\ref{p:5.1.2}, we have
$$0\longrightarrow \tau^{-1}P(7)\oplus X_{21}\longrightarrow
E_1'
 \oplus E _2'\oplus 2E_3'\longrightarrow X_{22}\longrightarrow 0.$$

Applying $Hom(~,X_{22})$ (resp. $Hom(E_1'
 \oplus E_2'\oplus 2E_3',~)$) to the
 sequence above, we get $$dim Hom(E_1'
 \oplus E_2'\oplus 2E_3',X_{22})\leq 5 \text{ and }dim Hom(E_1'
 \oplus E_2'\oplus 2E_3',X_{22})\geq 6. $$

 This gives a contradiction.  Thus $X_2$ must contain one regular direct summand at least.

 According to the $AR-$quiver of $\widetilde{E}_7,$ we can deduce
 that
 $X\cong\tau^{-1}P(7)\bigoplus
 E_1'
 \oplus E_2'\oplus E_3'\bigoplus \tau^2I(7).$

 Thus
$$\begin{array}{l}
{p_3|_{\overline{\co}_{\oplus E_1'\bigoplus\oplus
E_2'\bigoplus\oplus2 E_3'}\times \co_S}}^{-1}(\tau^{-1}P(7)\oplus
\tau^{-4}P(5))=V(S, E_1'\bigoplus\oplus E_2'\bigoplus\oplus2
E_3';\tau^{-1}P(7)\oplus \tau^{-4}P(5))\\
{\qquad\qquad\qquad\qquad\qquad\qquad}\text{\d{$\cup$}}V(S,I(8)\oplus\tau^2I(1));\tau^{-1}P(7))\text{\d{$\cup$}}
({p_3|_{\overline{\co}_{E_1'\oplus E_2'\oplus E_3'}\times
\co_{P(5)}}})^{-1}(\tau^{-4}P(5))
 .\end{array}$$

By $m=p=r=1,$ we deduce that $({p_3|_{\overline{\co}_{E_1'\oplus
E_2'\oplus E_3'}\times \co_{P(5)}}})^{-1}(\tau^{-4}P(5))$ has the
purity property. So does ${p_3|_{\overline{\co}_{ E_1'\bigoplus
E_2'\bigoplus\oplus2 E_3'}\times \co_S}}^{-1}(\tau^{-1}P(7)\oplus
\tau^{-4}P(5)).$

Suppose $P\cong P_1\oplus P_2\oplus P_3$ with $P_i$ indecomposable.
Without loss of generality, we may set
$P\cong\tau^{-1}P(7)\oplus\tau^{-2}P(7)\oplus\tau^{-4}P(6).$

For any $f\in Hom(\tau^{-2}P(7),E_2'),$ we have $S(6)\vdash
E_2'/Im(f).$

By $Hom(\tau^{-1}P(7)\oplus\tau^{-4}P(6),E_2')=0,$ we deduce that
$\tau^{-1}P(7)\oplus\tau^{-2}P(7)\oplus\tau^{-4}P(6)$ is not an
extension of $S$ by $E_1'\bigoplus E_2'\bigoplus\oplus 2E_3'.$

Similarly, suppose $P\cong P_1\oplus P_2\oplus P_3\oplus P_4$ with
$P_i$ indecomposable, then we may show that there is not an
epimorphism from $P$ to $E_1'\bigoplus E_2'\bigoplus\oplus 2E_3'.$
Thus $P$ is not an extension of $S$ by $E_{1}'\bigoplus
E_2'\bigoplus\oplus 2E_3'.$

The proof is complete.  \end{proof}

\mk \nd{\bf Remark} 7.1.1.
 Let $P$ be a non-trivial extension of $\oplus tS$ by the regular
 semi-simple objects in $\ct_i$ for some $i.$

 Since $$\begin{array}{l}
 {\bf Z}_{S,\oplus (t-1)S;\oplus tS}=\mathbb{P}^1,\\
\mathbb{P}^1\times({p_3|_{\overline{\co}_M*\co_{\oplus
tS}}})^{-1}(P) \\
{\qquad\qquad}=\text{\d{$\cup$}}_{P'\oplus L, L\vdash
M}({p_3|_{\co_M*\co_{S}}})^{-1}(P'\oplus L) \times
({p_3|_{\overline{\co}_{P'\oplus
L}*\co_{\oplus(t-1)S}}})^{-1}(P)\\
{\qquad\qquad\quad}\text{\d{$\cup$}}
({p_3|_{\overline{\co}_{P''\oplus L'\oplus I',P''\oplus L'\oplus
I'\in \overline{\co}_M, L'\vdash M}*\co_{\oplus
tS}}})^{-1}(P)\times\mathbb{P}^1,\end{array}$$

 we have by induction on $t$ and $\udim M$ that Proposition~\ref{p:6.2} is true when $S$
is replaced by $\oplus m S.$

 \begin{proposition}\label{p:6.4} Let $i$ be a sink, and $P,P'$ be
  pre-projective modules of $k Q$ (except $Q=\widetilde{E}_8$), and let $M$ be a regular semi-simple
module in $\ct_j$ for some $j,1\leqslant j\leqslant l.$

  Then
   ${p_3|_{\overline{\co}_M*\co_P}}^{-1}(P'\oplus L)$ has the purity
  property,
   where  $L$ is a
 submodule of $M$ in $\ct_j.$\end{proposition}

\begin{proof} Let $P=\oplus a_1P_1\bigoplus\oplus
a_2P_2\bigoplus\oplus\cdots \bigoplus\oplus a_tP_t$ where
$Ext^(P_u,P_v)=0$ for $u<v$ and $P_u$ is indecomposable.

 First, suppose that $P=\oplus a_1P_1.$ Applying the reflection functor $\sz_i$, we can turn the
question into the case of $P=\oplus a_1 P(i).$   By Proposition~\ref{p:6.2}
and Remark~7.1.1, the statement is true.

By the properties of pre-projective components,
 the statement above is true if $M$ is replaced by $P'\oplus
M$, where $P'$ is pre-projective.

 Next, suppose that $P\cong \oplus a_1P_1\bigoplus\oplus
a_2P_2\bigoplus\oplus\cdots \bigoplus\oplus a_tP_t $ and $t\geq 2.$
Without loss of generality, we may only consider the case $t=2.$

Since
$$\begin{array}{l}
{p_3|_{\overline{\co}_M*\co_{\oplus aP_1\bigoplus\oplus
aP_2}}}^{-1}(P'\oplus L)=\text{\d{$\cup$}}_{P''\oplus L',L'\vdash
M}({p_3|_{\overline{\co}_{M}\times \co_{a_1P_1} }})^{-1}(P''\oplus
L')\times ({p_3|_{\overline{\co}_{P''\oplus L}\times \co_{a_2P_2}
}})^{-1}(P'\oplus L),
\end{array}$$
it follows from the case $t=1$ that
${p_3|_{\overline{\co}_M*\co_{\oplus aP_1\bigoplus\oplus
aP_2}}}^{-1}(P'\oplus L)$ has the purity property. The proof is
complete.\end{proof}

Dually, we have the following statement

\begin{proposition}\label{p:6.5} Let $i$ be a source, and $I,I'$ be
 pre-injective modules of $k Q$(except $Q=\widetilde{E}_8$), and let $M$ be a regular semi-simple
module in $\ct_j$ for some $j,1\leqslant j\leqslant l.$

  Then
   ${p_3|_{\co_I*\overline{\co}_M}}^{-1}(I'\oplus L)$ has the purity
  property,   where  $L$ is a
 submodule of $M$ in $\ct_j$ . \end{proposition}

\begin{proposition}\label{p:6.6}  Let $i$ be a sink, and $P,P'$
  pre-projective modules of $k Q$ ( except $Q=\widetilde{E}_8$ ), and let $M$ be a regular
module in $\ct_j$ for some $j,1\leqslant j\leqslant l.$

  Then
   ${p_3|_{\overline{\co}_M*\co_P}}^{-1}(P'\oplus L)$ has the purity
  property,
   where  $L$ is a
 submodule of $M$ in $\ct_j$ .\end{proposition}

\begin{proof} By using induction on $\udim M,$ we prove that
${p_3|_{\overline{\co_M}*\co_P}}^{-1}(P'\oplus L)$ have purity
property.

The case that $M$ is a semi-simple object in $\ct_j$ is proved in
Proposition~\ref{p:6.4}.

 Assume that $M$ is not semi-simple in $\ct_i$. Let $M'$ be the
 direct sum of indecomposable summands of $M$ with maximal length in
 the full subcategory $\ct_i$ of $\llz-mod.$ Set $M_2=soc_{\ct_i}(M').$

 According to Proposition~2.5 in \cite{GJ}, there is a regular module $M_1$ in
 $\ct_i$ such that $g^M_{M_1M_2}=1.$

 Since
$$\left\{\begin{array}{l}\co_{M_1}*\co_{M_2}=p_3p_2p_1^{-1}(\co_{M_1}\times \co_{M_2}),\\
  {\bf Z}_{M,M_1,M_2}=p_2p_1^{-1}(\co_{M_1}\times \co_{M_2})\cap p_3^{-1}(M)=\{1 \text{ point }\},\end{array}\right.$$
  and $\co_{M_1}*\co_{M_2}$ has only
finitely many orbits,
  we thus get
 \begin{eqnarray}\label{f:6.18}{\qquad\qquad}\left\{\begin{array}{l}dim \co_{M_1}*\co_{M_2}=dim \co_M,\\
 \overline{\co}_{M_1}*\overline{\co}_{M_2}=\overline{\co}_M.\end{array}\right.\end{eqnarray}

Thus
$$\begin{array}{l}{p_3|_{\overline{\co}_M*\co_P}}^{-1}(P'\oplus L)
={p_3|_{\overline{\co}_{M_1}*\overline{\co}_{M_2}*\co_P}}^{-1}(P'\oplus
L)\\
{\qquad\qquad\qquad\qquad\quad}=\text{\d{$\cup$}}_{P''\oplus
L',L'\vdash M_2}{p_3|_{\overline{\co}_{M_1}*\co_{P''\oplus
L'}}}^{-1}(P'\oplus L)\times
{p_3|_{\overline{\co}_{M_2}*\co_P}}^{-1}(P''\oplus L').
\end{array}$$

Since $M_2$ is semi-simple in the full subcategory, the disjoint
union above makes sense.  By Proposition~\ref{p:6.4}, we have that
${p_3|_{\overline{\co}_{M_2}*\co_P}}^{-1}(P''\oplus L')$ has the
purity property.

In addition, it follows from the induction hypothesis and the proof
of proposition~2.5 in \cite{GJ} that
${p_3|_{\overline{\co}_{M_1}*\co_{P''\oplus L'}}}^{-1}(P'\oplus L)$
has the purity property.

Thus the proof is complete.\end{proof}




\section{Proof of Theorem~\ref{t:5.1.5}}

\subsection{}The aim of this subsection is to prove Theorem~\ref{t:5.1.5}.

Let $Q$ be a tame quiver not of type $\widetilde{E}_8$.

%

\begin{lemma}\label{l:7.1.1} Let $M$ and $N$ be modules
belonging to different nonhomogeneous tube $\ct_i$ and $\ct_j$ in $k
Q-mod$. Assume that $\overline{\co}_M$ and $\overline{\co}_N$
have the purity property.  Then $\overline{\co}_{M\oplus N}$ has the purity property.
\end{lemma}


\begin{proof}
 Set $\az=\udim M,\bz=\udim N,$ and
$X=\overline{\co}_{M\oplus N}.$ Let
$$P=\left\{\left(\begin{array}{ll}
A&C\\
0&B
\end{array}\right)\left|\left(\begin{array}{ll}
A&C\\
0&B
\end{array}\right)\in GL_{\az+\bz}, A\in GL_{\az},B\in GL_{\bz}\right. \right\}.$$
Then $P$ is a parabolic subgroup of $GL_{\az+\bz}.$
It is well known that the closure $Y$ of
 $GL_{\az+\bz}/P$ has the purity property.

Consider the natural projection
$$p: GL_{\az+\bz}/(\aut_{\llz}(M)\times \aut_{\llz}(N))\longrightarrow GL_{\az+\bz}/P. $$
It is easy to see that the fibre of $p$ is $P/(\aut(M)\times
\aut(N)).$ Thus we have a long exact sequence
$$\ldots\longrightarrow H^i_c(Z,\overline{\bbq}_l)\longrightarrow H^i_c(X,\overline{\bbq}_l)
\longrightarrow H^i_c(Y,\overline{\bbq}_l)\longrightarrow\ldots,$$
where $Z$ is the closure of $P/(\aut_{\llz}(M)\times \aut_{\llz}(N))$.

 By the definition of
$\aut_{\llz}(M)$ and $\aut_{\llz}(N),$ we have
$\aut_{\llz}(M)\subseteq GL_{\az}$ and $ \aut_{\llz}(N)\subseteq GL_{\bz}$. Therefore,
\[P/(\aut_{\llz}(M)\times \aut_{\llz}(N))\cong
GL_{\az}/\aut_{\llz}(M)\times GL_{\bz}/\aut_{\llz}(N)\times\Pi_{i\in
I}k^{\az_i\bz_i}.\]
It follows that
the closure $Z$ of $P/(\aut_{\llz}(M)\times \aut_{\llz}(N))$
has the purity property, since by assumption the closures  of $GL_{\az}/\aut_{\llz}(M)$ and
$GL_{\bz}/aut_{\llz}(N)$ both have the purity property. Now the desired result follows from the above long exact sequence.
\end{proof}

Note that $\Hom_{\llz}(M,N)$ is a subspace of $\Pi_{i\in
I}\Hom_k(k^{\az_i},k^{\bz_i}).$  In the same way as in
 Lemma~\ref{l:7.1.1}, we have


 \begin{lemma}\label{l:7.1.2} Let $M$ and $N$ be modules of $k Q$ such that
  $Ext_{\llz}^1(M,N)=0, Hom_{\llz}(N,M)=0$. Assume that
$\overline{\co}_M$ and $\overline{\co}_N$
both  have the purity property.  Then $\overline{\co}_{M\oplus N}$ has the purity property.
\end{lemma}


\begin{proposition}\label{p:7.1.3}  Let $M$ be a regular $k Q$-module in a nonhomogeneous tube $\ct_i$ for $1\leqslant
i\leqslant l.$
  Then $\overline{\co}_{M}$ has the purity property.
\end{proposition}

\begin{proof} We use induction on $\udim M.$

The case that $M$ is a semi-simple object in $\ct_i$ is proved in
Lemma~\ref{l:5.2.1}-Lemma~\ref{l:5.2.5}.

 Assume that  $M$ is not semi-simple in $\ct_i$. Let $M'$ be the
 direct sum of indecomposable summands of $M$ with maximal length in
 the full subcategory $\ct_i$ of $\llz-mod.$ Set $M_2=soc_{\ct_i}(M').$

 According to Proposition~2.5 in \cite{GJ}, there is a regular module $M_1$ in
 $\ct_i$ such that $g^M_{M_1M_2}=1.$ It follows from $(\ref{f:6.18})$ that
 $ \overline{\co}_{M_1}*\overline{\co}_{M_2}=\overline{\co}_M.$

By induction hypothesis and Lemma~\ref{l:5.2.1}-Lemma~\ref{l:5.2.5}, we obtain
$\overline{\co}_{M_1}$ and $\overline{\co}_{M_2}$ have the purity
property.

On the other hand, because there are finitely many orbits in
$\overline{\co}_{M_2},$  for any $P\oplus R\oplus I\in
\overline{\co}_{M_2},$ by proposition~\ref{p:6.6} and the proof of
Proposition~2.5 in \cite{GJ},
$({p_3|_{\overline{\co}_{M_1}*\co_{P\oplus R}}})^{-1}(P'\oplus R')$
has the purity property. So does
$({p_3|_{\overline{\co}_{M_1}*\co_{P\oplus R\oplus
I}}})^{-1}(P'\oplus R'\oplus I').$

Thus
$$({p_3|_{\overline{\co}_{M_1}*\overline{\co}_{M_2}}})^{-1}(P'\oplus R'\oplus
I')=\text{\d{$\cup$}}_{P\oplus R\oplus I\in \overline{\co}_{
M_2}}({p_3|_{\overline{\co}_{M_1}*\co_{P\oplus R\oplus
I}}})^{-1}(P'\oplus R'\oplus I').$$

Consider the proper morphism
$$p_3:p_2p_1^{-1}(\overline{\co}_{M_1}\times \overline{\co}_{M_2})\longrightarrow
\overline{\co}_{M_1}* \overline{\co}_{M_2}.$$

From the discussions above, we obtain that $\overline{\co}_M$ has
the purity property.
\end{proof}

We now prove Theorem~\ref{t:5.1.5}.
\begin{proof}[Proof of Theorem~\ref{t:5.1.5}]By \cite{L4} and \cite{L5},
$\overline{\co}_P$ (resp. $\overline{\cn}_{{\bf w},3}$
$\overline{\co}_I$) has the purity property. It now follows from
Proposition~\ref{p:7.1.3} that $\overline{\co}_M$ has the purity property.

According to Lemma~7.1.2, $X=\overline{\co}_{P,M,{\bf w},I}$ has
purity property,  the desired conclusion follows from above.
\end{proof}

\bigskip

\section{ Application of $\ch^s(\llz)$}

\bigskip

\subsection{} In \cite{H}, A.Hubery has proved the existence of Hall
polynomials for tame quivers for Segre classes. In this subsection,
by using the extension algebras of singular Ringel-Hall algebras, we
give a simple and direct proof for the existence of Hall polynomials
for tame quivers.

\mk\nd  Let $Q$ be a affine quiver . We define a new algebra
$\ch^{s,x_1,x_2,\cdots,x_t}$ which is generated by $\{u_i,u_{[M]},
u_{[N]}: i\in I, M\in\ct_j, N\in\oplus_{i=1}^t\mathscr{H}_{x_i}
,1\leq j\leq l, \},t\in\bbn,$ where $x_t$ are $\bbf_q$ rational
points in $ \mathbb{P}^1.$

 We now give a new decomposition of $E_{n\dz}$ as follows

   $$ E_{n\delta}=E_{n\delta,1}+E_{n\delta,2}+E_{n\delta,3} , $$

 where
 \begin{eqnarray}\label{f:8.1.1}
   E_{n\delta,1}=v^{-n \dim S_1-n \dim S_2}
   \sum_{[M],M\in\fk{C}_1\bigoplus\oplus_{x_t}\mathscr{H}_{x_t},\udim
   M=n\delta}u_{[M]},
     {\qquad\qquad\qquad\quad}\end{eqnarray}
 \begin{eqnarray}\label{f:8.1.2}
    E_{n\delta,2}=v^{-n \dim S_1-n \dim S_2}
   \sum_{\substack{[M],\udim M=n\delta,\\
   M=M_1\oplus M_2, 0\neq M_1\in \fk{C}_1\bigoplus\oplus_{x_t}\mathscr{H}_{x_t}, 0\neq M_2\in
   \fk{C}_0\backslash\oplus_{x_t}\mathscr{H}_{x_t}}}u_{[M]},\end{eqnarray}
 \begin{eqnarray}\label{f:8.1.3}
   E_{n\delta,3}=
   v^{-n \dim S_1-n \dim S_2}\sum_{[M], M\in\fk{C}_0\backslash\oplus_{x_t}\mathscr{H}_{x_t},
   \udim M=n\delta}u_{[M]}. {\qquad\qquad\qquad\quad\quad}
  \end{eqnarray}

  Let
$\textbf{w}=(w_1,\cdots,w_t)$ be a partition of $n,$ we then define
$$E_{\textbf{w}\dz,3}=E_{w_1\dz,3}*\cdots*E_{w_t\dz,3}.$$

Let ${\bf P}(n)$ be the set of all partitions of $n,$ and
$\lr{N}=v^{-dim N+dim End(N)}u_{[N]}.$ Set
$${\bf
B}'=\{\lr{P}*\lr{M}*E_{\textbf{w}\dz,3}*\lr{I}|~ P\in\cp_{prep},
M\in\oplus_{i=1}^l\ct_i\bigoplus\oplus_{x_t}\mathscr{H}_{x_t},
I\in\cp_{prei},{\bf w}\in{\bf P}(n), n\in\bbn\}.$$

Similar to Theorem ~4.1.1,   we have the following:

\begin{proposition}\label{p:8.1.1 } {\sl  The set ${\bf B}'$ is a
$\bbq(v)-$ basis of $\ch^{s,x_1,x_2,\cdots,x_t}$.} \end{proposition}

\begin{theorem}\label{t:8.1.2}  Let $Q$ be a affine quiver ,  let
$P_i$ (resp. $R_i,I_i$) be a pre-projective (resp. nonhomogeneous
regular, pre-injective) $\bbf_q Q-$module, and let $H_i\in
\oplus_{j=1}^{t}\mathscr{H}_{x_j}$ be a homogeneous regular $\bbf_q
Q-$module with $x_j$ being $\bbf_q-$ rational point in
$\mathbb{P}^1$ for $i=1,2,3; t\in \bbn.$ Let $X_i=P_i\oplus
R_i\oplus H_i\oplus I_i,i=1,2,3.$

Then there exists a Hall polynomial
$\varphi_{X_1X_2}^{X_3}(x)\in\bbq[x]$ such that
$$\varphi_{X_1X_2}^{X_3}(q)=g_{X_1X_2}^{X_3}.$$ \end{theorem}

\begin{proof} Since
$$\lr{P_1\oplus M_1\oplus I_1}*\lr{P_2\oplus M_2\oplus I_2}=
a_{12}^3(v)\lr{P_3\oplus M_3\oplus I_3}+\text{ other terms},$$ and
$a_{12}^3(v)\in \bbq (v),$ then we have
$$\begin{array}{l}g_{P_1\oplus R_1\oplus H_1\oplus I_1,P_2\oplus
R_2\oplus H_2\oplus I_2}^{P_3\oplus R_3\oplus
H_3\oplus I_3} \\
=v^{dim_{\bbf_q}End(P_3\oplus M_3\oplus
I_3)-dim_{\bbf_q}End(P_1\oplus M_1\oplus
I_1)-dim_{\bbf_q}End(P_2\oplus M_2\oplus I_2)-\lr{\udim P_1\oplus
R_1\oplus H_1\oplus I_1, \udim P_2\oplus R_2\oplus H_2\oplus I_2}
}a_{12}^3(v).\end{array}$$

Set

$$(*)\begin{array}{l}\varphi_{X_1X_2}^{X_3}(v^2)\\
=v^{dim_{\bbf_q}End(P_3\oplus M_3\oplus
I_3)-dim_{\bbf_q}End(P_1\oplus M_1\oplus
I_1)-dim_{\bbf_q}End(P_2\oplus M_2\oplus I_2)-\lr{\udim P_1\oplus
R_1\oplus H_1\oplus I_1, \udim P_2\oplus R_2\oplus H_2\oplus I_2}
}a_{12}^3(v).\end{array}$$

On the other hand, we know that
$$v^{dim_{\bbf_q}End(P_3\oplus M_3\oplus
I_3)-dim_{\bbf_q}End(P_1\oplus M_1\oplus
I_1)-dim_{\bbf_q}End(P_2\oplus M_2\oplus I_2)-\lr{\udim P_1\oplus
R_1\oplus H_1\oplus I_1, \udim P_2\oplus R_2\oplus H_2\oplus I_2}
}a_{12}^3(v)$$  takes the positive integer value while $v^2$ takes
infinite many positive integer values. Then
$\varphi_{X_1X_2}^{X_3}(v^2)$ is a polynomial of $v^2$ over $\bbq.$
Thus the proof is complete.\end{proof}

\end{document}